\documentclass[10pt, 3p]{elsarticle}

\makeatletter
\def\ps@pprintTitle{%
  \let\@oddhead\@empty
  \let\@evenhead\@empty
  \let\@oddfoot\@empty
  \let\@evenfoot\@oddfoot
}
\makeatother

\journal{Computers \& Mathematics with Applications}

\usepackage[utf8]{inputenc}

\usepackage{algorithm2e}
\SetKwInput{KwData}{Input}

\usepackage{amsmath, amsfonts, amssymb, amsthm}
\usepackage{hyperref, cleveref}
\usepackage{rotating, lscape}

\usepackage{bm}
\newtheorem{proposition}{Proposition}

\usepackage{xcolor, graphicx}
\usepackage{subfigure}

\usepackage{acronym}
\acrodef{AMR}{Adaptive Mesh Refinement}
\acrodef{FVSG}{Finite-Volume Scharfetter-Gummel}
\acrodef{EAFE}{Edge-Average Finite-Element}
\acrodef{DAR}{Diffusion--Advection--Reaction}
\acrodef{FD}{Finite Differences}
\acrodef{FE}{Finite Elements}
\acrodef{FV}{Finite Volumes}
\acrodef{PDE}{Partial Differential Equation}

\newcommand{\norm}[1]{\left\|#1\right\|}

\begin{document}

\begin{frontmatter}
    
    \title{Scalable Recovery-based Adaptation on Quadtree Meshes for Advection-Diffusion-Reaction Problems}
    
    \author{Pasquale Claudio Africa}
    \author{Simona Perotto}
    \author{Carlo de Falco}
    \address{MOX -- Dipartimento di Matematica, Politecnico di Milano, P.zza L. Da Vinci 32, 20133 -- Milano -- Italy}
    
\begin{abstract}
We propose a mesh adaptation procedure for Cartesian quadtree meshes, to discretize scalar advection-diffusion-reaction problems.
The adaptation process is 
driven by a recovery-based a posteriori estimator for the $L^2(\Omega)$-norm of the discretization error, based on suitable higher order approximations of both the solution and the associated gradient. In particular, a metric-based approach exploits the information furnished by the estimator to iteratively predict the new adapted mesh.
The new mesh adaptation algorithm is successfully assessed on different configurations, and turns out to perform well also when dealing with discontinuities in the data as well as in the presence of internal layers not 
aligned with the Cartesian directions.
A cross-comparison with a standard 
estimate--mark--refine approach and with other adaptive strategies available in the literature shows the remarkable accuracy and parallel scalability 
of the proposed approach.
\end{abstract}
    
    
\end{frontmatter}


\section{Introduction}
Quadtree-based decomposition of two-dimensional spatial domains is a well-established technique whereby 
a rectangular region is 
hierarchically subdivided into smaller rectangles in such a way that each 
refinement consists of splitting one \emph{parent} rectangle into exactly four 
\emph{children}, by halving its edges.
The quadtree structure was initially devised for the organization of spatial 
data~\cite{samet85,morton1966computer}, and was greatly successful due 
to the early development of efficient algorithms for
data sorting and retrieval~\cite{Finkel1974,samet1984quadtree}.

While, in the past, simplicial meshes have been slightly favoured over quadtrees
because
of their superior geometric flexibility, a recent substantial surge in the research
interest in methods based on Cartesian quadtree meshes has been driven by two
main factors.
On the one hand, the many advancements in the field of the embedded interface methods 
(e.g., the immersed boundary or the fictitious domain approaches), that allow to accurately 
track interior interfaces, i.e., lower dimensional manifolds within the 
computational domain where problem coefficients may exhibit discontinuities, without
constructing meshes that accurately represent the geometry of the 
interface~\cite{hasbestan2017parallel,TOWERS2018424,LEE2021109958,burmancutfem,massingcutdg}.
On the other hand, the recent development of highly scalable software 
libraries~\cite{BursteddeWilcoxGhattas11,IsaacBursteddeWilcoxEtAl15,Burstedde20,LOURENCO2020124721,pablo,p4esturl,basilisk} 
based on quadtrees (or octrees in 3D) has identified Cartesian meshes as the ideal
tool to efficiently manage extreme scale problems that require
massively parallel computer architectures.

As a consequence, in recent 
years, many instances of numerical schemes based on quadtrees have been used to discretize
partial differential equation problems stemming from diverse applications, such as multiphase
and free--surface fluid flows~\cite{LAURMAA2016190,LIANG2009221} or 
semiconductor device simulation~\cite{CHEN2003341}. In the recent literature,  
some works focus on
correctly adapting quadtree (or octree) meshes to the geometry of the domain boundary as well as to possible internal interfaces~\cite{LAURMAA2016190,raeli2018finite,WDFGMFSPCDFLF,coco2018second}. 

In this paper, we propose a new metric-driven mesh adaptation procedure to sharply approximate the solution to a scalar advection-diffusion-reaction problem.
The adaptation is driven by an a posteriori estimator for the $L^2(\Omega)$-norm of the discretization error. In particular, we adopt a recovery-based estimator.
In 1987, O.C.~Zienkiewicz and J.Z.~Zhu propose a very practical way to evaluate the $H^1(\Omega)$-seminorm of the discretization error, instrumental to generate a mesh able to capture the strong gradient areas of the solution~\cite{ZZ1987}. The idea, successively formalized in the seminal papers~\cite{zz3,Zienkiewicz1992}, is very straightforward. It consists in making explicitly computable the $H^1(\Omega)$-seminorm of the discretization error by replacing the exact gradient with the so-called recovered gradient obtained by properly averaging/projecting the gradient of the discrete solution.\\
Recovery based error estimators proved to be 
reliable and computationally efficient
in different engineering contexts (see, e.g.,~\cite{yan01,porta12,mu13,simpaty}), by offering a mathematical tool characterized by several 
good properties, such as the independence of the specific problem and of the adopted discretization, the handy implementation and the computational cheapness.\\ 
In more recent years, the idea of making computable a quantity by means of suitable recovery procedures has been exploited also to evaluate norms of the discretization error different from the $H^1(\Omega)$-seminorm. Examples can be found, for instance, in~\cite{wiberg94,MAISANO20064794,fpc22}, where the authors propose different recipes to recover the exact solution in order to make computable a certain norm of the associated discretization error.
Here, we pursue a similar goal in the case of non-conforming Cartesian, 
quadtree meshes, by generalizing the approach proposed in~\cite{MAISANO20064794}.
We build a higher-order approximation for the solution gradient by means of suitable finite difference formulas. Then, the recovered gradient is exploited to compute a higher-order approximation of the solution, in order to estimate the $L^2(\Omega)$-norm of the discretization error. In particular, 
we use one-sided formulas along the interior interfaces in order to guarantee accuracy when coefficient discontinuities, even not aligned with the Cartesian directions,  occur
within the computational domain.
Moreover, in order to improve the scalability of the adaptation procedure, the
same one-sided formulas are used in correspondence with the boundaries
of the parallel partition of the computationl domain,
thus resulting in a reduction of the inter--process communication. 
Finally, a metric-based approach is shown to further improve the method efficiency, by reducing the number of estimation--adaptation steps with respect to a standard estimate--mark--refine approach~\cite{verfurth}.

The manuscript is organized as follows: in Section~\ref{sec:formulation} we 
introduce the reference partial differential equation problem, together with the corresponding 
discrete formulation based on edge-averaged finite elements 
on quadtree meshes. Section~\ref{sec:recovery} is devoted to the recovery procedures, and proposes specific recipes to enrich both the discrete gradient and solution. The error estimator for
the $L^2(\Omega)$-norm of the discretization error is introduced in Section~\ref{sec:refinement}, together with the adopted metric-based mesh adaptation strategy able to predict, in an automatic way, the level of refinement/coarsening of the grid cells. For comparison reasons, in the same section we provide also a standard marking-based mesh adaptation algorithm.
An extensive numerical assessment is carried out in Section~\ref{sec:results}, while some concluding remarks are drawn in the last section.


\section{The reference setting
}\label{sec:formulation}
As reference context, we adopt the 
conservative form of a 
standard scalar advection-diffusion-reaction (ADR) problem, completed with mixed boundary conditions
\begin{equation}
\label{eq:adeq}
\left\{
\begin{alignedat}{2}
&-\nabla\cdot\left(\varepsilon\nabla u - \bm{\beta}u\right) + b u  = f & \quad & \text{in } \Omega\subset \mathbb{R}^2 \\
&u  = g & \quad & \text{on } \Gamma_\mathrm{D} \subseteq \partial\Omega \\
&\varepsilon \nabla u \cdot \bm{\nu}  = 0 & \quad & \text{on } \partial\Omega\setminus\Gamma_\mathrm{D},
\end{alignedat}
\right.
\end{equation}
with $\epsilon>0$ the diffusive coefficient; $\bm{\beta}$
a conservative advective field, i.e.,  
$\bm{\beta} = \nabla\psi$ for a certain potential $\psi$, and such that
$\nabla\cdot\bm{\beta} = 0$; $b\ge 0$, the reaction; $f$ the source term; $g$ the data assigned on the Dirichlet portion, $\Gamma_D$, of the domain boundary, $\partial \Omega$; $\bm{\nu}$ the unit outward normal vector to $\partial \Omega$. In principle, the problem data can be assumed discontinuous in $\Omega$.

The two next sections are devoted to the discretization of problem \eqref{eq:adeq} by resorting to a finite
element scheme. In particular,
we will distinguish between Cartesian-product and non-conforming quadtree meshes, after identifying $\Omega$ with a Cartesian domain.
We highlight that non-Cartesian domains can be handled as well, by using a penalty approach (see, for instance,~\cite{raeli2018finite} and  Section~\ref{sec:penalty} for an example).

\subsection{An edge-averaged finite element discretization on a Cartesian-product grid}
A Cartesian-product mesh is a quadrilateral grid obtained by considering the Cartesian product between a one-dimensional (1D) partition of the edges of $\Omega$ aligned with the $x$- and the $y$-direction, respectively. We denote by 
$\tau_h = \left\{\Omega^{(k)}\right\}_{k=1}^{N_\mathrm{el}}$ a family of Cartesian-product partitions of \(\Omega\)  such that
$\Omega= \cup_{k=1}^{N_\mathrm{el}} \Omega^{(k)}$, with \(N_\mathrm{el}\) the total number of mesh elements and $h = \min_k \left(\mathrm{diam}\left(\Omega^{(k)}\right)\right)$ the minimum diameter of cells $\Omega^{(k)}$.\\
We introduce the  
finite element space, \(Q_h^1(\tau_h)\), of the continuous piecewise bi-linear polynomials associated with the partition $\tau_h$, where
$Q_h^n(\tau_h)  = Q_h^{n,n}(\tau_h)$,
being
$$
Q_h^{n,m}(\tau_h)  = \Big\{ v \in \mathcal{C}^0\left(\overline{\Omega}_h\right) : v|_{\Omega^{(k)}} \in \mathbb{P}_{n,m}(\Omega^{(k)})\ \forall \Omega^{(k)} \in \tau_h\Big\}
$$
with $n$, $m\in \mathbb N$, and
\begin{equation*}
\mathbb{P}_{m,n}(\Omega^{(k)}) = \Big\{ p_{m, n}(x, y) = \sum_{0\le i \leq m, \, 0\le j \leq n} a_{ij} x^i y^j,\ \mbox{with\ } {\bf x}=(x, y) \in \Omega^{(k)} \Big\}.
\end{equation*}
From now on, we simplify the notation by dropping any explicit dependence on partition \(\tau_h\). 
After denoting by \(X=\{\bm{x}_i\}_{i=1}^N\) the set of the nodes defining the partition \(\tau_h\),
we consider a Lagrangian basis $\{\varphi_i\}_{i=1}^N$ of the space $Q_h^1$, so that  
\begin{equation}
\label{eq:phi}
\varphi_i(\bm{x}_j) = \delta_{ij}\quad \forall \bm{x}_j \in X,
\end{equation}
\(\delta_{ij}\) being the Kronecker delta, and 
the discrete solution \(u_h\) to problem \eqref{eq:adeq} can be expressed as a linear combination of such basis functions through the set of (unknown) coefficients \(\{U_i\}_{i=1}^N\), being
\begin{equation*}
u_h(\bm{x}) = \sum_{i=1}^N U_i \varphi_i(\bm{x}) \quad \forall \bm{x} \in \Omega.
\end{equation*}
In particular, we require that $u_h$ belongs to the subset $V_h\subset Q_h^1$ of the functions satisfying a homogeneous Dirichlet data 
in correspondence with $\Gamma_D$ to match the essential boundary conditions in \eqref{eq:adeq}. Finally, we observe that functions $\{ \varphi_i\}_{i=1}^N$ satisfy a partition of unity property, namely
\begin{equation}\label{PU}
\displaystyle\sum_i \varphi(\bm{x})=1 \quad \forall \bm{x} \in \Omega.
\end{equation}

To discretize problem \eqref{eq:adeq}, we adopt the \ac{FVSG} method proposed in \cite{bank1998finite}. This discretization was originally introduced by D.N. De G.~Allen and R.V. Southwell in \cite{de1955relaxation} for 1D geometries and successively generalized to two-dimensional (2D) and three-dimensional (3D) diffusion-convection equations by M.A. Zl\'{a}mal in \cite{zlamal1986finite} and to self-adjoint problems by P.A. Markowich and M.A. Zl\'{a}mal in \cite{markowich1988inverse}.\\
In particular, first the original equations are rewritten in terms of the so-called Slotboom variables, this leading to a self-adjoint problem;
then, such a problem is discretized by means of a primal-mixed formulation using the harmonic average of the diffusion coefficient over the mesh edges; finally, the problem is rewritten in terms of the primal variables by inverting the Slotboom relations, so that arithmetic overflows are prevented during computations. 
The trapezoidal quadrature rule is adopted to approximate the integrals resulting from the formulation.\\
As an example, the matrix yielded by the FVSG discretization when applied to an advection-diffusion problem (e.g., to problem \eqref{eq:adeq} with $b=0$) and associated with the generic element \(\Omega^{(k)}\) in $\tau_h$, 
whose vertices are numbered by $1$ to $4$ in a lexicographical order, 
is given by
\begin{equation}
\begin{bmatrix}
\mathcal{B}^-_{12} + \mathcal{B}^+_{31} & -\mathcal{B}^+_{12} & -\mathcal{B}^-_{31} & 0 \\
-\mathcal{B}^-_{12} & \mathcal{B}^+_{12} + \mathcal{B}^-_{24} & 0 & -\mathcal{B}^+_{24} \\
-\mathcal{B}^+_{31} & 0 & \mathcal{B}^+_{43} + \mathcal{B}^-_{31} & -\mathcal{B}^-_{43} \\
0 & -\mathcal{B}^-_{24} & -\mathcal{B}^+_{43} & \mathcal{B}^-_{43} + \mathcal{B}^+_{24}
\end{bmatrix},
\label{eq:localmatrixtp}
\end{equation}
with
\begin{equation*}
\mathcal{B}^+_{ij} = \frac{\ell_{ij} h^{(k)}_{ij}(\varepsilon) B(\psi_i - \psi_j)}{2|\Omega^{(k)}|}, \qquad
\mathcal{B}^-_{ij} = \frac{\ell_{ij} h^{(k)}_{ij}(\varepsilon) B(\psi_j - \psi_i)}{2|\Omega^{(k)}|},
\end{equation*}
for $ij\in \{12, 31, 24, 43\}$, \(\ell_{ij}\) the length of the edge, \(ij\), joining vertices \(i\) and \(j\), \(\psi_k\) the value of the potential \(\psi\) at node \(k\) for $k=i, j$, \(B(x) = x/(e^x - 1)\) the Bernoulli function, \(h^{(k)}_{ij}(f) = |\ell_{ij}|^{-1} \int_{ij} f^{-1}\) the harmonic mean operator over the edge \(ij\), $|\varrho|$ denoting the measure of a generic set $\varrho\subset \mathbb R^d$ for $d=1, 2$.

It has been checked that, under the choices above and over a Cartesian-product mesh, the local matrix \eqref{eq:localmatrixtp} ensures the fully-assembled matrix to be an M-matrix. This, in turn, guarantees that the FVSG discretization is monotone, i.e., that a discrete maximum principle holds for the numerical solution to \eqref{eq:adeq} \cite{brezzi1989two}.

\subsection{An edge-averaged finite element discretization on a non-conforming grid}\label{sec2}
Non-conforming quadrilateral meshes stem from refinement and coarsening procedures applied to a Cartesian-product mesh.\\
The quadtree data structure characterizing a Cartesian-product grid allows to easily deal with non-conforming mesh refinement and coarsening, where the former operation consists in replacing an element with four children of equal size, while the latter occurs by replacing four children with their parent~\cite{BursteddeWilcoxGhattas11}. 
Both refinement and coarsening procedures can lead to the generation of hanging nodes. As a consequence, the mesh cannot be represented as a Cartesian-product set, and the \(Q_h^1\) polynomial space turns out to be not well-defined. Thus, space \(Q_h^1\) is replaced by a new function environment, \(\tilde{Q}_h^1\), obtained by imposing in \(Q_h^1\) a set of constraints in order to account for the presence of hanging nodes -- which are not included as additional degrees of freedom (dofs) -- and to preserve the partition of unity property \eqref{PU}. 

We demand that neighboring cells may differ by at most one refinement level, thereby enforcing that only a single hanging node can exist per face. This procedure, referred to as 2:1 balance, simplifies the derivation of the discretization scheme and the recovery procedures presented below. This balance also ensures that the workload is well-balanced among parallel processors, since the computational burden becomes proportional to the number of local cells, rather than to the number of dofs~\cite{bangerth2012algorithms}, thus making the proposed approach well-suited to frequent adaptations in large-scale computing scenarios~\cite{isaac2012low}.

In the remaining part of the section, we analyze more into the details how to manage refinement and coarsening of non-conforming meshes, 
with a focus also on the definition of space $\tilde{Q}_h^1$. 
To this aim, we denote by \(\tilde{\tau}_h\) the non-conforming mesh obtained by refinement and/or  coarsening of the Cartesian-product grid \(\tau_h\), and by \(\tilde{X}\) the set of nodes defining \(\tilde{\tau}_h\).

\begin{figure}
    \centering
    \includegraphics[width=0.99\textwidth]{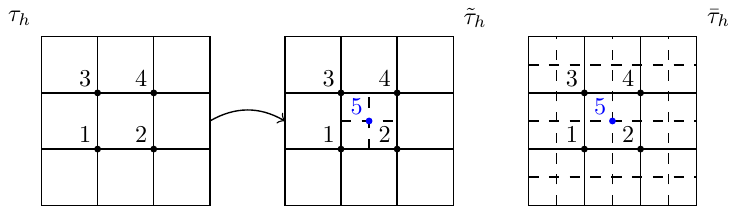}
    \caption{Example of refinement for a non-conforming mesh: initial configuration (left); refined mesh (center), where node $5$ corresponds to the newly added degree of freedom; mesh used to define the basis function \(\bar{\varphi}_5\) (right).}
    \label{fig:refine}
\end{figure}

\paragraph{Refinement of a non-conforming mesh}
When refining a mesh element in a non-conforming framework, an additional node (where an additional dof is located) has to be added (see Figure~\ref{fig:refine}, where node \(5\) denotes the added dof).

\noindent
Now, we introduce the partition \(\bar{\tau}_h\) obtained by uniformly refining \(\tau_h\) and coinciding with the minimum-size Cartesian-product mesh containing \(\tilde{\tau}_h\), and the basis function, \(\bar{\varphi}_5 \in Q_h^1\left(\bar{\tau}_h\right)\), 
associated with dof \(5\) over the uniformly refined mesh \(\bar{\tau}_h\).\\
Now, we associate functions
\begin{equation*}
\tilde{\varphi}_i = \varphi_i - \frac{1}{4}\bar{\varphi}_5 \quad i=1 \dots 4,
\end{equation*}
with the four vertices, \(1, \dots, 4\), identifying the element to be refined (see Figure~\ref{fig:refine}), and we define the set
\(\varPhi = \big\{\{\varphi_i\}_{i \neq 1, \dots, 4}, \{\tilde{\varphi}_i\}_{i=1}^4, \bar{\varphi}_5\big\}\). 
Thus, the new polynomial space associated with the non-conforming grid $\tilde{\tau}_h$ can be defined as
\begin{equation*}
\tilde{Q}_h^1 = \mathrm{span}\left(\varPhi\right).
\end{equation*}

\noindent
With an abuse of notation, from now on, we denote by $\tilde{\varphi}_i$ the generic function in $\varPhi$.
We remark that $\tilde{Q}_h^1$ is a linear subspace of $Q_h^1\left(\bar{\tau}_h\right)$ and that the basis functions in $\varPhi$ satisfy the partition of unity property, 
\begin{equation}
\label{eq:partition_unity}
\sum_{i} \tilde{\varphi}_i(\bm{x}) = 1 \quad \forall \bm{x} \in \Omega,
\end{equation}
as well as the Lagrangian property, 
\begin{equation}
\label{eq:lagrangian}
\tilde{\varphi}_i(\bm{x}_j) = \delta_{ij} \quad \forall \bm{x}_j \in \tilde{X},
\end{equation}
analogously as in \eqref{eq:phi} and \eqref{PU}.

\begin{figure}
    \centering   \includegraphics[width=0.7\textwidth]{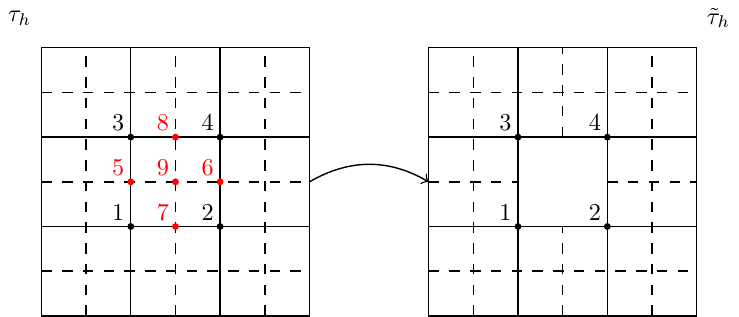}
    \caption{Example of coarsening for a non-conforming mesh: initial configuration (left) where the red nodes highlight the degrees of freedom to be removed; coarsened mesh (right).}
    \label{fig:coarsen}
\end{figure}

\paragraph{Coarsening of a non-conforming mesh}
We exemplify such an operation in Figure~\ref{fig:coarsen}, where we assume that the cell with vertices $1$, $2$, $3$, $4$ has to be coarsened by removing node $9$. This action leads to 
get rid of five dofs (the ones associated with the nodes $5$, $6$, $7$, $8$, $9$), since nodes $5$, $6$, $7$, $8$ become hanging and node \(9\) does not yet represent a mesh node.

\noindent
To characterize space $\tilde{Q}_h^1$,
we associate the new functions
\begin{equation*}
\tilde{\varphi}_i = \varphi_i + \frac{1}{4}\varphi_9 + \frac{1}{2}\sum_{j \in 
\mathcal{N}_i} \varphi_j \quad i \in 1 \dots 4,
\end{equation*}
with vertices $1$, $2$, $3$, $4$,
where \(\mathcal{N}_i\) denotes the set of the hanging nodes sharing an edge with node 
\(i\), and we define the new set 
 \(\varPhi =
\{\{\varphi_i\}_{i \neq 1, \dots, 9}, \{\tilde{\varphi}_i\}_{i=1}^4\}\). 
Thus, the new discrete space associated with the non-conforming coarsened grid turns out to be
\begin{equation*}
\tilde{Q}_h^1 = \mathrm{span}\left(\varPhi\right),
\end{equation*}
which coincides with 
a linear subspace of space \(Q_h^1(\tau_h)\). Functions in \(\varPhi\) still satisfy the partition of unity and the
Lagrangian property \eqref{eq:partition_unity} and \eqref{eq:lagrangian}, 
respectively.

\vspace*{.2cm}

In practice, a discretization based on a non-conforming mesh generated by refinement and/or coarsening of a Cartesian-product grid, is managed first by assembling the standard local matrices associated with the basis functions of the space \(Q_h^1\); then, the resulting system is extended through a set of equations that constrain the solution at each hanging node to be the arithmetic mean of the solution at the two parent vertices; finally, the constraints are eliminated by static condensation~\cite{BK07}. This technique reduces the size of the linear system by minimizing the coupling among dofs.

\section{Recovery procedures}\label{sec:recovery}
In this section, we present the recovery techniques used to define the error estimator driving the automatic mesh adaptation procedure. To this aim,
we extend the results in~\cite{MAISANO20064794} for triangular grids to quadtree meshes.\\ In particular, given a finite element solution, \(u_h \in \tilde{Q}_h^1\), and the associated gradient, \(\bm{\sigma}_h=\nabla u_h \in \tilde{Q}_h^{0,1} \times \tilde{Q}_h^{1, 0}\), first we derive a formula to recover an enriched gradient, \(\bm{\sigma}_h^* = \nabla^* u_h\in \big[\tilde{Q}_h^1\big]^2\); successively, we employ $\bm{\sigma}_h^*$ to build an enriched solution, \(u_h^* \in \tilde{Q}_h^2\), coinciding with a piecewise bi-quadratic polynomial.

\subsection{Gradient Recovery}\label{grad_rec_sec}
\begin{figure}
    \centering
\includegraphics[height=0.4\textwidth]{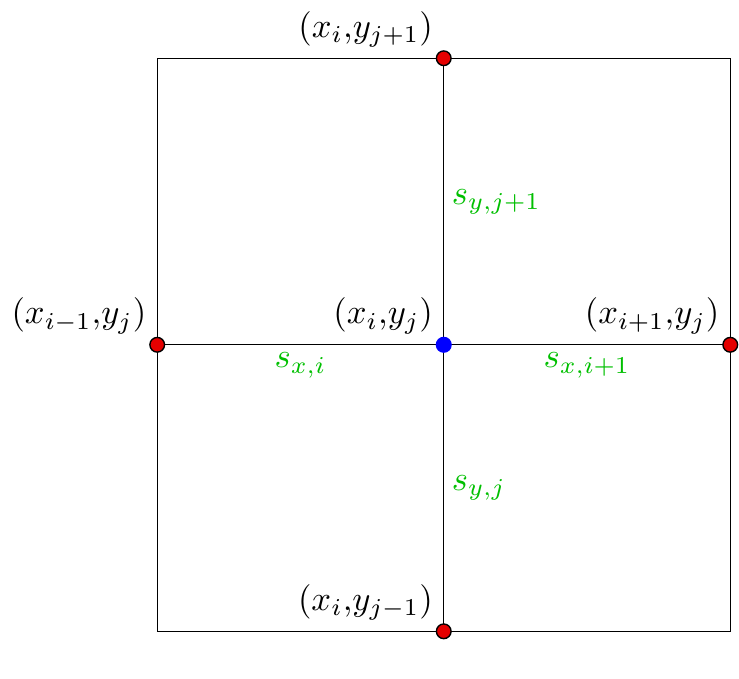}
    \quad
        \includegraphics[height=0.4\textwidth]{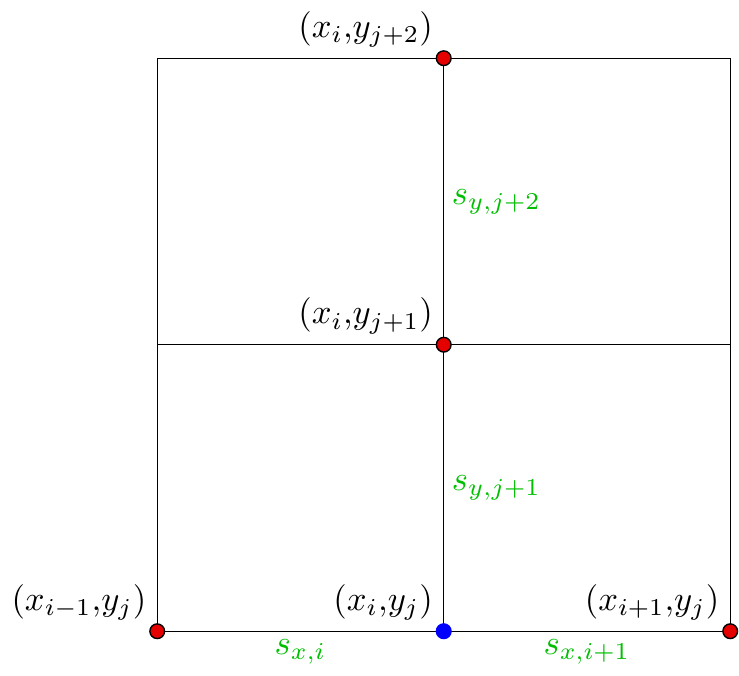}
    \caption{Gradient recovery: patch associated with the internal (left) and the boundary (right) node \((x_i, y_j)\). The mesh edges used in the recovery procedure are labeled in green.}\label{fig:recovery}
\end{figure}

To compute the recovered gradient, we resort to an averaging technique.\\ 
We denote by $\{s_{x,i}\}_{i=1}^{N_{\mathrm{s},x}}$ and $\{s_{y,j}\}_{j=1}^{N_{\mathrm{s},y}}$ the sets of the mesh edges oriented along the $x$ and $y$ axis, respectively. Now, we
consider the patch of elements, $\{\Omega^{(k)}\in \tau_h \, :\, (x_i, y_j)\in \Omega^{(k)}\}$, associated with the internal node \((x_i, y_j) \in \Omega\) (see Figure~\ref{fig:recovery} (left)).
Since the solution \(u_h\) is a bi-linear polynomial, the gradient components can be associated with the degrees of freedom of the Nédélec finite element space~\cite{nedelec1980mixed}, and can be exactly computed by using the finite difference formulas
\begin{equation}
\label{eq:sigma_star}
\begin{alignedat}{3}
\bm{\sigma}_h\big|_{s_{x,i}} & =
\begin{bmatrix}
\displaystyle \frac{u_h(x_i, y_j) - u_h(x_{i-1}, y_j)}{h_{x,i}} \\
0
\end{bmatrix},
\
\bm{\sigma}_h\big|_{s_{y,j}} & =
\begin{bmatrix}
0 \\
\displaystyle \frac{u_h(x_i, y_j) - u_h(x_i, y_{j-1})}{h_{y,j}}
\end{bmatrix},
\end{alignedat}
\end{equation}
with $h_{x,m} = x_m - x_{m-1}$ and $h_{y,n} = y_n - y_{n-1}$ the length of the generic edges \([x_{m-1}, x_m]\) and \([y_{n-1}, y_n]\). In the sequel, with an abuse of notation, we denote restrictions $\bm{\sigma}_h\big|_{s_{x,i}}$ and $\bm{\sigma}_h\big|_{s_{y,j}}$ by $\bm{\sigma}_h(s_{x,i})$ and $\bm{\sigma}_h(s_{y,j})$, respectively.\\
Then, following \cite{MAISANO20064794}, we compute the recovered gradient, $\bm{\sigma}_h^*$, at $(x_i, y_j)$ by averaging the values of the discrete gradient, $\bm{\sigma}_h$, along the $x$- and $y$-direction, respectively, with weights depending on the reciprocal edge lengths, so that we obtain
\begin{equation}
\label{eq:recovery}
\begin{array}{rcl}
\bm{\sigma}_h^*(x_i, y_j) &=
& \displaystyle \frac{1}{\displaystyle \frac{1}{h_{x,i}} + \frac{1}{h_{x,i+1}}} \left(\frac{\bm{\sigma}_h(s_{x,i})}{h_{x,i}} + \frac{\bm{\sigma}_h(s_{x,i+1})}{h_{x,i+1}}\right)  \\[7mm]
& +& \displaystyle \frac{1}{\displaystyle \frac{1}{h_{y,j}} + \frac{1}{h_{y,j+1}}} \left(\frac{\bm{\sigma}_h(s_{y,j})}{h_{y,j}} + \frac{\bm{\sigma}_h(s_{y,j+1})}{h_{y,j+1}}\right).
\end{array}
\end{equation}

A particular care has to be taken when \((x_i, y_j)\) is a mesh node on the boundary $\partial \Omega$. In such a case, the patch adopted to compute the recovered gradient becomes unilateral and formula \eqref{eq:recovery} has to be modified accordingly. For instance, with reference to Figure~\ref{fig:recovery} (right),
we have to consider a patch which is unilateral along the \(y\)-direction, while standard along the $x$-direction, \((x_i, y_j)\) coinciding with a node along the bottom boundary of $\Omega$.  Formula \eqref{eq:recovery} is consequently modified into
\begin{equation}
\label{eq:recovery_boundary}
\begin{array}{rcl}
\bm{\sigma}_h^*(x_i, y_j) &=& \displaystyle\frac{1}{\displaystyle\frac{1}{h_{x,i}} + \frac{1}{h_{x,i+1}}} \left(\frac{\bm{\sigma}_h(s_{x,i})}{h_{x,i}} + \frac{\bm{\sigma}_h(s_{x,i+1})}{h_{x,i+1}}\right)  \\[7mm]
&+& \displaystyle\frac{1}{\displaystyle\frac{1}{h_{y,j+1}} + \frac{1}{h_{y,j+2}}} \left[\bm{\sigma}_h(s_{y,j+1}) \left(\frac{1}{h_{y,j+1}} + \frac{2}{h_{y,j+2}}\right) - \frac{\bm{\sigma}_h(s_{y,j+2})}{h_{y,j+2}}\right],
\end{array}
\end{equation}
which falls back to a standard three-point forward-difference formula.
Once a value for the recovered gradient is computed at each dof in $\tilde{Q}_h^1$, we are able to define the recovered gradient \(\bm{\sigma}_h^* \in [\tilde{Q}_h^1]^2\). Notice that the presence of hanging nodes does not modify the procedure, since formulas  \eqref{eq:recovery}-\eqref{eq:recovery_boundary} employ only the tangential component of the gradient along each mesh edge, which is continuous also in the case of non-conforming meshes.

We remark that we exploit the
one-sided formulas in \eqref{eq:recovery_boundary} also at the 
interfaces between regions that are assigned to different processors 
in a distributed-memory parallel processing paradigm. 
This strategy allows us to highly reduce the amount of inter--process 
communication in the recovery procedure, thus ensuring a better scalability (see Section~\ref{sec:results} for more details).

The following property features the gradient recovery procedure just formalized:

\begin{proposition}\label{prop:gradient}
    Let \(\Omega \subset \mathbb{R}^2\) be an open, bounded domain and let \(\tau_h = \left\{\Omega^{(k)}\right\}_{k}\) be a quadtree partition of \(\Omega\). Let \(u \in \mathbb{P}_{2,2}(\Omega)\) be a bi-quadratic polynomial. Then, the recovery procedure \eqref{eq:recovery}-\eqref{eq:recovery_boundary} applied to the bi-linear interpolant \(\varPi_h^1 u \in \tilde{Q}_h^1\) of \(u\) exactly recovers the gradient \(\bm{\sigma} = \nabla u\) at the partition vertices.
    
    \begin{proof}
        With reference to Figure~\ref{fig:recovery}, we aim to prove that \(\bm{\sigma}_h^*(x_i, y_j) = \bm{\sigma}(x_i, y_j)\) for each internal node \((x_i, y_j)\), being $\bm{\sigma}_h^*$ and $\bm{\sigma}_h$ a piecewise linear and a linear polynomial, respectively. To show this, we
        proceed component-wise. Since \(u\) is a bi-quadratic polynomial, the following relations hold for the \(x\)-component \(\sigma_x\) of the exact gradient:
        \begin{equation}
        \label{eq:proof1}
        \begin{aligned}
        \sigma_x(s_{x,i}) &= \displaystyle \frac{\varPi_h^1 u(x_i, y_j) - \varPi_h^1 u(x_{i-1}, y_j)}{h_{x,i}}, \\[5mm]
        \sigma_x(s_{x,i+1}) &= \displaystyle \frac{\varPi_h^1 u(x_{i+1}, y_j) - \varPi_h^1 u(x_i, y_j)}{h_{x,i+1}}, \\[5mm]
        \sigma_x(x_i, y_j) &= \displaystyle \frac{1}{\displaystyle \frac{1}{h_{x,i}} + \frac{1}{h_{x,i+1}}} \left(\frac{\sigma_x(s_{x,i})}{h_{x,i}} + \frac{\sigma_x(s_{x,i+1})}{h_{x,i+1}}\right),
        \end{aligned}
        \end{equation}
        where the last equality follows from the linearity of \(\bm{\sigma}\).
        
        On the other hand, procedure \eqref{eq:recovery} applied to the computation of the \(x\)-component \(\sigma_{h,x}^*\) of the recovered gradient of \(\varPi_h^1 u\) yields
        \begin{equation}
        \label{eq:proof2}
        \sigma_{h,x}^*(x_i, y_j) = \displaystyle \frac{1}{\displaystyle \frac{1}{h_{x,i}} + \frac{1}{h_{x,i+1}}} \left(\frac{\sigma_{h,x}(s_{x,i})}{h_{x,i}} + \frac{\sigma_{h,x}(s_{x,i+1})}{h_{x,i+1}}\right),
        \end{equation}
        with $\sigma_{h,x}$ the $x$-component of the discrete gradient, which, by using \eqref{eq:sigma_star} and the definition of interpolation, provides \(\sigma_x(x_i, y_j) = \sigma_{h,x}^*(x_i, y_j)\).\\
        The same procedure can be repeated for the \(y\)-component of the gradient and for boundary nodes, which concludes the proof.
    \end{proof}
\end{proposition}

In particular, if \(u \in \mathbb{P}_{2,2}(\Omega)\), then \(\bm{\sigma}_h^*(s_{x,i}) = \bm{\sigma}(s_{x,i})\) and, analogously, \(\bm{\sigma}_h^*(s_{y,j}) = \bm{\sigma}(s_{y,j})\). Thus, we can state that, on each mesh edge, the tangential component of the recovered gradient, \(\bm{\sigma}_h^*\), provides an exact representation of the tangential component of the exact gradient, \(\bm{\sigma}\).

As a consequence of Proposition~\ref{prop:gradient}, we have also numerically verified that the recovered gradient is super-convergent with respect to the mesh size $h$ (we refer, e.g., to Fig.~\ref{fig:dr1_convergence}). 
In particular, it holds
\begin{equation*}
\norm{\bm{\sigma}_h^* - \bm{\sigma}}_{L^2(\Omega)} = O(h^2).
\end{equation*}

\subsection{Solution Recovery}\label{recsol_sec}
\begin{figure}
    \centering
        \includegraphics[height=0.25\textwidth]{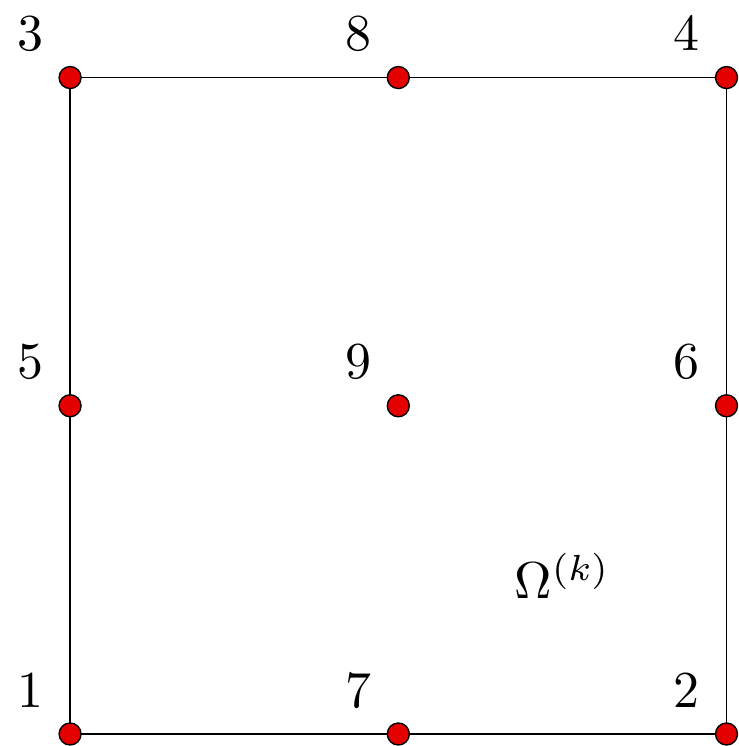}
    \quad
        \includegraphics[height=0.25\textwidth]{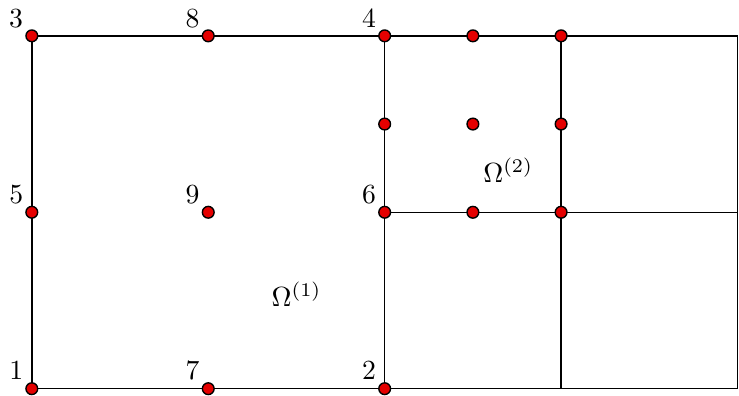}
    \caption{Solution recovery: local numbering of dofs in \(\Omega^{(k)}\) characterizing the space \(\tilde{Q}_h^2\) (left); 
    sketch of a non-conforming configuration when interested in recovering the solution at dof \(6\) (right).}\label{fig:dofs_q2}
\end{figure}

In this section we present a procedure that, starting from the recovered gradient \(\bm{\sigma}_h^*\) in \eqref{eq:recovery}-\eqref{eq:recovery_boundary}, enables us to recover an enriched solution \(u_h^* \in \tilde{Q}_h^2\). We remind that the polynomials in this space are identified by \(9\) dofs on each mesh element $\Omega^{(k)}$.\\ 
To explain the recovery procedure, we adopt the numbering of the dofs in Figure~\ref{fig:dofs_q2} (left), and we
denote by \(u_{h,i} = u_h(\bm{x}_i)\) the value taken by the discrete solution \(u_h\) at the mesh node \(\bm{x}_i\). Now, we set the recovered solution to coincide with the discrete solution at the vertices of \(\Omega^{(k)}\), i.e.,  
\begin{equation}
\label{eq:recover_vertices}
u_{h,i}^*\big|_{\Omega^{(k)}} = u_h^*\big|_{\Omega^{(k)}}(\bm{x}_i)= u_{h,i}\quad i=1, \dots, 4,
\end{equation}
while we properly integrate the recovered gradient \(\bm{\sigma}_h^*\) to compute the values of \(u_h^*\big|_{\Omega^{(k)}}\) at the dofs discriminating a \(\tilde{Q}_h^2\) from a \(\tilde{Q}_h^1\)-approximation (namely, dofs \(5, \dots, 9\)).
We exemplify the adopted recovery process on dof $5$. 
After denoting by
\begin{equation}
\label{eq:midpoint1}
\begin{aligned}
u_5^{(1)} & = u_{h,1} + \int_{y_1}^{y_5} \sigma_{h,y}^* \ \mathrm{d}y, \qquad
u_5^{(2)} & = u_{h,3} - \int_{y_5}^{y_3} \sigma_{h,y}^* \ \mathrm{d}y
\end{aligned}
\end{equation}
with $\sigma_{h,y}^*$ the $y$-component of the recovered gradient, we assign the value
\begin{equation}
\label{eq:solrec1}
u_{h,5}^*\big|_{\Omega^{(k)}} = \frac{u_5^{(1)} + u_5^{(2)}}{2}
\end{equation}
to dof \(5\). We notice that the gradient recovery procedure in 
Section~\ref{grad_rec_sec} does not guarantee, a priori, that 
\(\bm{\sigma}_h^*\) is a conservative vector field. Thus, the two 
line integrals in \eqref{eq:midpoint1} could differ each as they 
depend on the integration path.
The same approach as in \eqref{eq:midpoint1}-\eqref{eq:solrec1} is 
adopted to recover the value of the solution at the midpoint dofs 
\(6, 7, 8\). \\Finally, the value associated with the centroid of the 
element (i.e., with dof $9$) is computed as
\begin{equation}
\label{eq:solrec2}
u_{h,9}^*\big|_{\Omega^{(k)}} = \frac{1}{4} \sum_{i=1}^{4} u_9^{(i)},
\end{equation}
with
\begin{equation}
\label{eq:midpoint2}
\begin{aligned}
u_9^{(1)} & = u_{h,5} + \int_{x_5}^{x_9} \sigma_{h,x}^* \ \mathrm{d}x, \qquad
u_9^{(2)} & = u_{h,6} - \int_{x_9}^{x_6} \sigma_{h,x}^* \ \mathrm{d}x, \\
u_9^{(3)} & = u_{h,7} + \int_{y_7}^{y_9} \sigma_{h,y}^* \ \mathrm{d}y, \qquad
u_9^{(4)} & = u_{h,8} - \int_{y_9}^{y_8} \sigma_{h,y}^* \ \mathrm{d}y,
\end{aligned}
\end{equation}
and $\sigma_{h,x}^*$, $\sigma_{h,y}^*$ the $x$- and the $y$-component of the recovered gradient, respectively. All the integrals in \eqref{eq:midpoint1} and \eqref{eq:midpoint2} are exactly computed by the midpoint quadrature rule, since the functions to be integrated are linear polynomials.

In the presence of non-conforming meshes, the procedure to recover the solution at hanging nodes is properly modified. With reference to Figure~\ref{fig:dofs_q2} (right), let us assume to be interested in recovering $u_h^*$ at node $6$.
Instead of setting \(u_{h,6}^*|_{\Omega^{(2)}}=u_{h,6}=\frac{1}{2}(u_{h,2}+u_{h,4})\), we assign
\begin{equation}
\label{eq:solrec3}
u_{h,6}^*\big|_{\Omega^{(2)}}=u_{h,6}^*\big|_{\Omega^{(1)}},
\end{equation}
i.e., the value recovered at the same point as a dof of the element \(\Omega^{(1)}\). It is straightforward to verify that this choice guarantees that the recovered solution \(u_h^*\) is continuous also across non-conforming edges.

\vspace*{.2cm}

A result analogous to Proposition~\ref{prop:gradient} can be proved also for the recovered solution, as stated below.

\begin{proposition}
    Let \(\Omega \subset \mathbb{R}^2\) be an open, bounded domain and let \(\tau_h = \left\{\Omega^{(k)}\right\}_{k}\) be a quadtree partition of \(\Omega\). Let \(u \in \mathbb{P}_{2,2}(\Omega)\) be a bi-quadratic polynomial. Then, the solution recovery procedure \eqref{eq:recover_vertices}-\eqref{eq:solrec3} applied to the bi-linear interpolant \(\varPi_h^1 u \in \tilde{Q}_h^1\) of \(u\) exactly recovers the solution \(u\) at the partition vertices.
\end{proposition}

Also for the recovered solution, we have numerically verified (see Fig.~\ref{fig:dr1_convergence}) a convergence of order \(2\) to the exact solution with respect to the $L^2(\Omega)$-norm, being
\begin{equation*}
\norm{u_h^* - u}_{L^2(\Omega)} = O(h^2).
\end{equation*}

\section{The mesh adaptation procedure}\label{sec:refinement}
To carry out the mesh adaptation process, we resort to the generalization of the approach proposed in~\cite{MAISANO20064794} 
to a quadtree tassellation of the domain $\Omega$. In particular, we adopt the a posteriori error estimator 
\begin{equation*}
    \eta = \Bigg[\sum_{k=1}^{N_\mathrm{el}} \eta_k^2\Bigg]^{\frac{1}{2}}
\end{equation*}
for the $L^2$-norm, $\|u-u_h \|_{L^2(\Omega)}$, of the discretization error, with
\begin{equation}\label{local_estim}
    \eta_k = \norm{u_h^* - u_h}_{L^2\left(\Omega^{(k)}\right)}
\end{equation}
the estimator associated with the mesh element \(\Omega^{(k)}\) (notice that the enhanced gradient, \(\bm{\sigma}_h^*\), is used only to recover the enriched solution \(u_h^*\)).\\
The current goal is to employ the information provided by $\eta$ and, in particular, by $\eta_k$, to generate a computational mesh able to follow the features of the solution at hand. To this aim, we exploit an iterative procedure which, starting from an initial (possibly uniform) grid, predicts a new mesh characterized by a refinement of the elements where the solution features an irregular behaviour, while exhibiting an element coarsening where a smooth trend is detected. 

Below, we investigate two alternative strategies to drive the mesh adaptation. The former is based on a standard 
solve $\rightarrow$ estimate $\rightarrow$ mark $\rightarrow$ adapt
algorithm~\cite{verfurth}, whereas the latter exploits the concept of metric~\cite{frey}.
In either cases, the strategy we adopt consists in minimizing the mesh cardinality while ensuring a user-defined accuracy, $\mathrm{tol}$, on the $L^2(\Omega)$-norm of the discretization error (i.e., on the estimator $\eta$), together with an error equidistribution across the grid elements,
so that 
\begin{equation}\label{eq:equidist}
    \eta_k \simeq \frac{\mathrm{tol}}{\sqrt{N_\mathrm{el}}}.
\end{equation}
This check can be tuned by a constant in order to penalize the elements that most infringe criterion \eqref{eq:equidist}, both from above and from below (we refer to parameters $\delta_1$ and $\delta_2$ in Algorithm~\ref{algo:marking}).

To stop the iterative procedure, different choices can be made. We set a maximum number of iterations to keep affordable the computational effort. As an alternative, we could monitor the stagnation 
of mesh elements by checking whether the relative variation in the mesh cardinality of two successive grids is smaller than a fixed threshold.

In the next sections, we introduce the marking-based and the metric-based approaches followed to generate adapted meshes, starting from estimator $\eta$.
All the numerical tests presented hereby have been performed using \texttt{bim++}, an in-house \texttt{C++} library, with utility classes
and methods to solve ADR problems, with an interface to \texttt{p4est} for quadtree mesh support.

\subsection{A marking-based mesh adaptation algorithm}\label{marksec}
As a first approach, we resort to a standard solve $\rightarrow$ estimate $\rightarrow$ mark $\rightarrow$ adapt iterative algorithm. \\
The whole procedure is itemized in Algorithm~\ref{algo:marking}, whose input parameters coincide with 
the initial mesh, $\tau_h^{(0)}$, the tolerance, $\mathrm{tol}$, set by the user on the estimator $\eta$, and the maximum number of iterations, \(i_{\mathrm{max}}\), constraining the adaptive loop.
\\ 
After computing the discrete solution on the current mesh, the local error estimator \eqref{local_estim} is evaluated, and the mesh elements are consequently marked to be refined, coarsened or kept unaltered, according to an error equidistribution principle. 
Then, cells $\Omega^{(k)}$ are modified consistently with the received label by resorting to the local hierarchical cell refinement and coarsening operations, according to what detailed in Section~\ref{sec2}.\\
Parameters \(\delta_1\) and \(\delta_2\)
are meant to confine the refinement and the coarsening phase to the elements characterized by the largest and by the smallest values for $\eta_k$ , respectively in terms of error equidistribution. In particular, we set \(\delta_1=1.5\) and \(\delta_2=0.5\) in the numerical assessment below.
\begin{algorithm}[h!]
    \KwData{$\tau_h^{(0)}$; \(\mathrm{tol}\); \(i_{\mathrm{max}}\);}
      \(i=0\)\;
    \While{\(i \leq i_{\mathrm{max}}\)}{
        compute \(u_h\)\;
        \For{\(k = 1, \ldots, N_{\mathrm{el}}\)}{
            compute \(\eta_k\)\;
            \uIf{\(\eta_k \geq \delta_1\displaystyle \frac{\mathrm{tol}}{\sqrt{N_\mathrm{el}}}\)}{
                mark \(\Omega^{(k)}\) for refinement\;
            }
            \uElseIf{\(\eta_k \leq \delta_2\displaystyle \frac{\mathrm{tol}}{\sqrt{N_\mathrm{el}}}\)}{
                mark \(\Omega^{(k)}\) for coarsening\;
            }
            \Else{keep \(\Omega^{(k)}\)}
        }
        refine and/or coarsen the mesh\;
        \(i = i+1\)\;
    }
    \caption{Workflow for the marking-based adaptation.}
    \label{algo:marking}
\end{algorithm}

\subsection{A metric-based mesh adaptation algorithm}
The second adaptive approach coincides with a a metric-based procedure. In such a case, the error estimator is used to directly predict the spacing (i.e., the metric) characterizing the new mesh, instead of selecting the operation (refinement/coarsening/preservation) to be performed on each cell.
In more detail, according to Section~\ref{marksec}, the prediction of the new spacing is driven by a certain accuracy to be guaranteed to the error estimator, together with 
the minimization of the mesh cardinality and the equidistribution of the error (namely, of the estimator).\\
Since we deal with isotropic meshes, the prediction of the spacing characterizing the new grid reduces to the computation of the cell side, for each cell of the new mesh.
To this aim, we first scale the local error estimator in \eqref{local_estim} with respect to the area of the corresponding element, in order to confine, at least asymptotically (i.e., when the mesh is sufficiently fine), the information related to the measure of
$\Omega^{(k)}$. This leads us to define the scaled estimator
\begin{equation}\label{uno}
    \hat{\eta}_k^2 = \frac{\eta_k^2}{\big[ h_x^{(k)}\big]^2},
\end{equation}
where $h_x^{(k)}$ denotes a characteristic size of
$\Omega^{(k)}$ such as, for instance, the cell edge length along the \(x\)-direction. 
Now, we exploit the equidistribution constraint \eqref{eq:equidist} to predict the new size \(h_{x,\mathrm{new}}^{(k)}\) of the element $\Omega^{(k)}$, being
\begin{equation}\label{due}
    \hat{\eta}_k^2 \big[h_{x,\mathrm{new}}^{(k)}\big]^2 \approx \frac{\mathrm{tol}^2}{N_\mathrm{el}}.
\end{equation}
Relations \eqref{uno} and \eqref{due}, combined together, yield
\begin{equation*}
    \frac{h_{x,\mathrm{new}}^{(k)}}{h_x^{(k)}} \approx \frac{\mathrm{tol}}{\eta_k \sqrt{N_\mathrm{el}}}.
\end{equation*}
We recall that, in the case of quadtree meshes, the relation \({h_{x,\mathrm{new}}^{(k)}} / {h_x^{(k)}} = 2^{-\ell_k}\) holds, so that the number, \(\ell_k\), of refinement or coarsening steps required to reach the desired size, \(h_{x,\mathrm{new}}^{(k)}\), can be computed as
\begin{equation}\label{prima_ell}
    \ell_k = \left\lceil \log_2\left(\frac{\eta_k \sqrt{N_\mathrm{el}}}{\mathrm{tol}}\right) \right\rceil,
\end{equation}
where the choice between refinement and coarsening depends on the sign of the ratio \({h_{x,\mathrm{new}}^{(k)}} / {h_x^{(k)}}\).\\
Finally, to have a sharper control on the quality of the generated mesh, we resort to two additional parameters, $n_{\mathrm{ref}}$ and $n_{\mathrm{coarsen}}$, with \(0 \leq n_{\mathrm{ref}}, n_{\mathrm{coarsen}} \leq \ell_k\), which modify the prediction for $ \ell_k$ in \eqref{prima_ell} into the new quantity
\begin{equation*}
    \ell_k =
    \left\{
    \begin{alignedat}{3}
    	& \mathrm{max} (0, \ell_k - n_{\mathrm{ref}}) & \quad & \text{if } \ell_k \geq 0, \\
        & \mathrm{min} (0, \ell_k + n_{\mathrm{coarsen}}) & \quad & \text{if } \ell_k < 0.
    \end{alignedat}
    \right.
\end{equation*}
Parameters $n_{\mathrm{ref}}$ and $n_{\mathrm{coarsen}}$ have a role similar to the one played by \(\delta_1\) and \( \delta_2\) in the marking-based approach, and provide a bound to the predicted maximum number of refinement/coarsening steps. This becomes essential when small tolerance values are assigned or when dealing with large meshes.

\noindent
The metric-based adaptation procedure is itemized in
Algorithm~\ref{algo:metrics}.
\begin{algorithm}[t]
    \KwData{$\tau_h^{(0)}$; \(\mathrm{tol}\); \(i_{\mathrm{max}}\)\;}
 \(i=0\)\;
    \While{\(i \leq i_{\mathrm{max}}\)}{
        compute \(u_h\)\;
        \For{\(k = 1, \ldots, N_{\mathrm{el}}\)}{
            compute \(\eta_k\)\;
            compute the number \(\ell_k\) of refinement/coarsening steps\;
            refine/coarsen \(\Omega^{(k)}\) \(\ell_k\) times\;
        }
        \(i = i+1\)\;
    }
    \caption{Workflow for the metric-based adaptation.}
    \label{algo:metrics}
\end{algorithm} 

\noindent
The predictive feature represents the main advantage of a metric-based approach. In general, a metric-based approach converges to the adapted mesh after a smaller number of iterations when compared with a marking-based algorithm.
Nevertheless, as we will verify in the next section,  Algorithm~\ref{algo:metrics} requires the initial mesh to be properly selected: for instance, should it be too coarse, \(\ell_k\) would be overestimated during the first steps, thus generating overly-refined grids at intermediate iterations.

\section{Numerical results}
\label{sec:results}
In this section, we assess the performance of both the marking-based and the metric-based mesh adaptation procedures on diffusion-reaction and advection-diffusion problems. The presence of boundary and internal layers as well as of a possible discontinuity in the problem data further support the adoption of a mesh customized to the problem at hand. Finally, the computational footprint of \Cref{algo:marking,algo:metrics} are analyzed, thus confirming the efficiency and the scalability of the proposed adaptation procedures.

\subsection{Test case 1: a diffusion-reaction problem with multiple layers}\label{prob:dr1}
\begin{figure}
    \centering
    \subfigure[Adapted mesh: $i=1$.]{\includegraphics[width=0.4\textwidth]{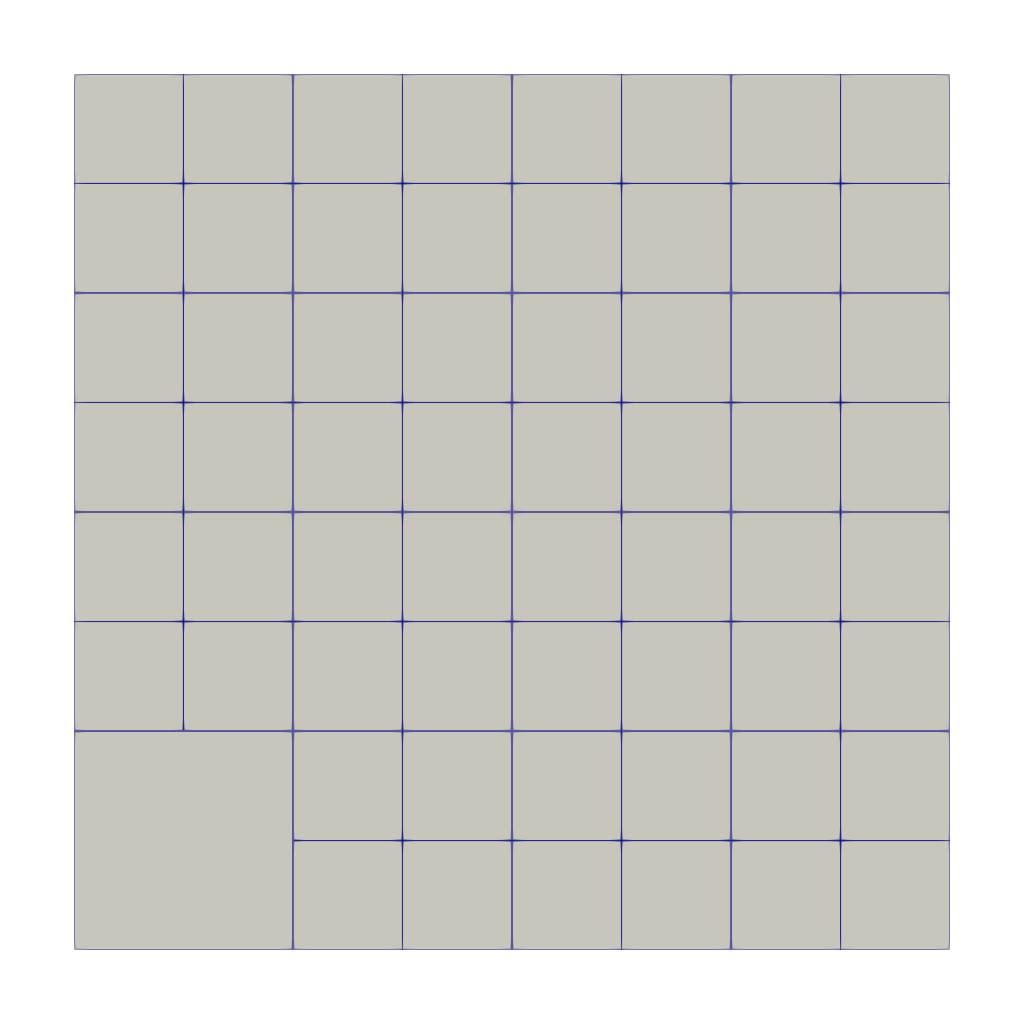}}
    \subfigure[Adapted mesh: $i=4$.]{\includegraphics[width=0.4\textwidth]{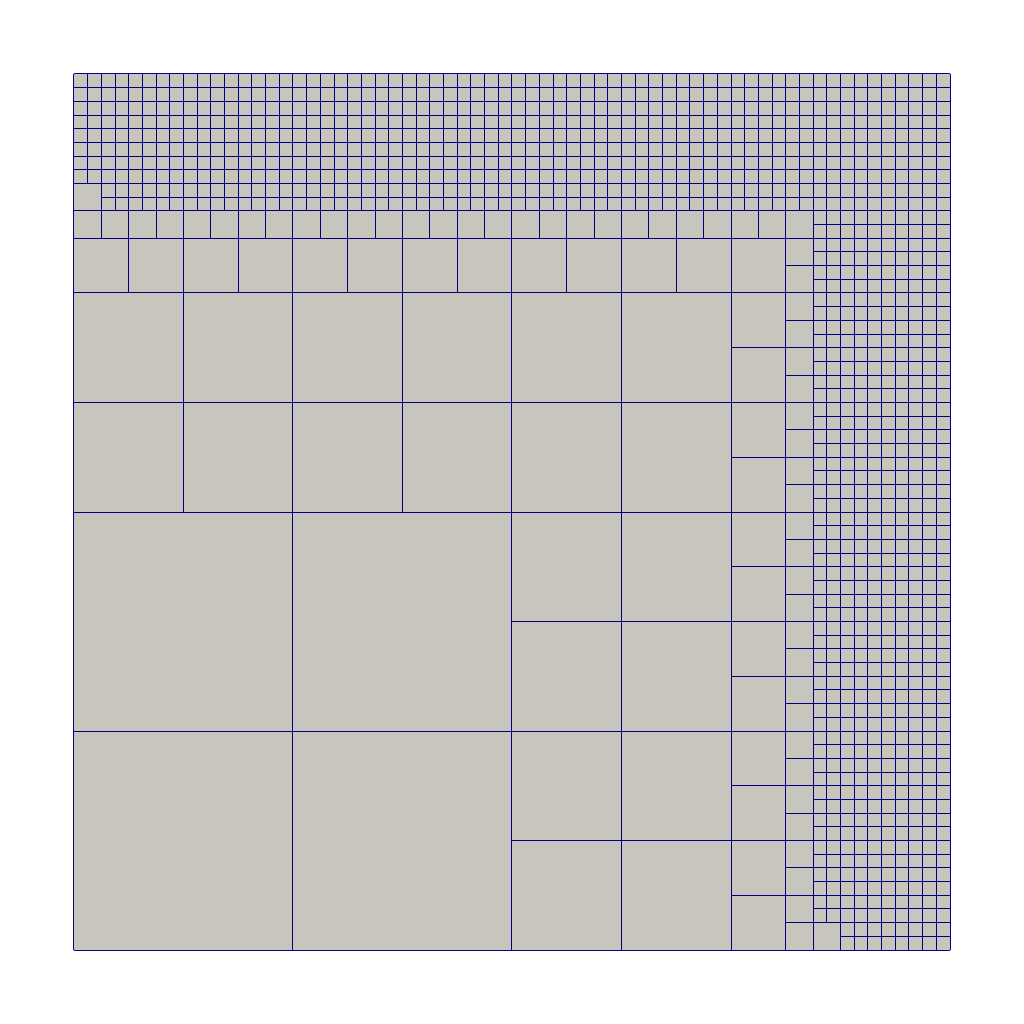}} \\
    \subfigure[Adapted mesh: $i=6$.]{\includegraphics[width=0.4\textwidth]{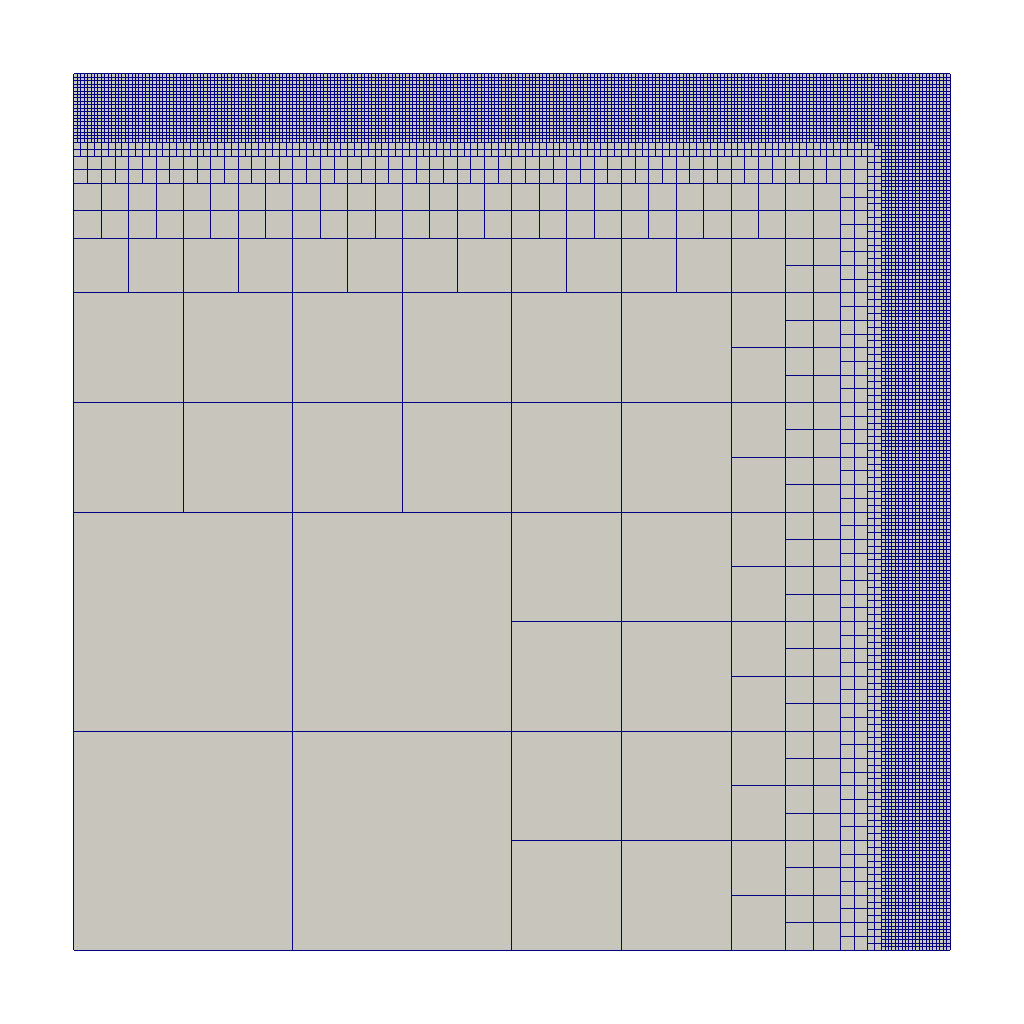}}
    \subfigure[Adapted mesh: $i=9$.]{\includegraphics[width=0.4\textwidth]{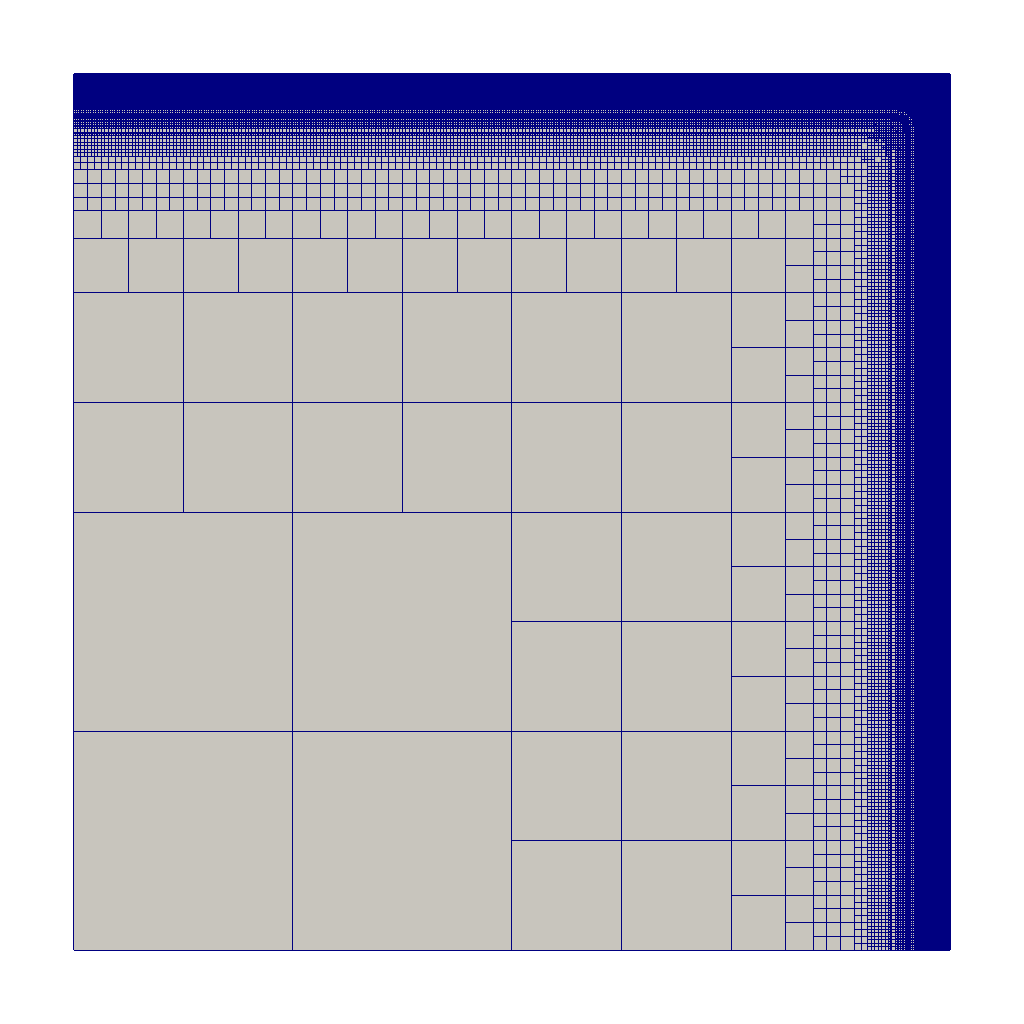}}
    \caption{Test case 1. Marking-based adaptation.}
    \label{fig:dr1_marker}
\end{figure}

\begin{figure}
    \centering
    \subfigure[Adapted mesh: $i=1$.]{\includegraphics[width=0.4\textwidth]{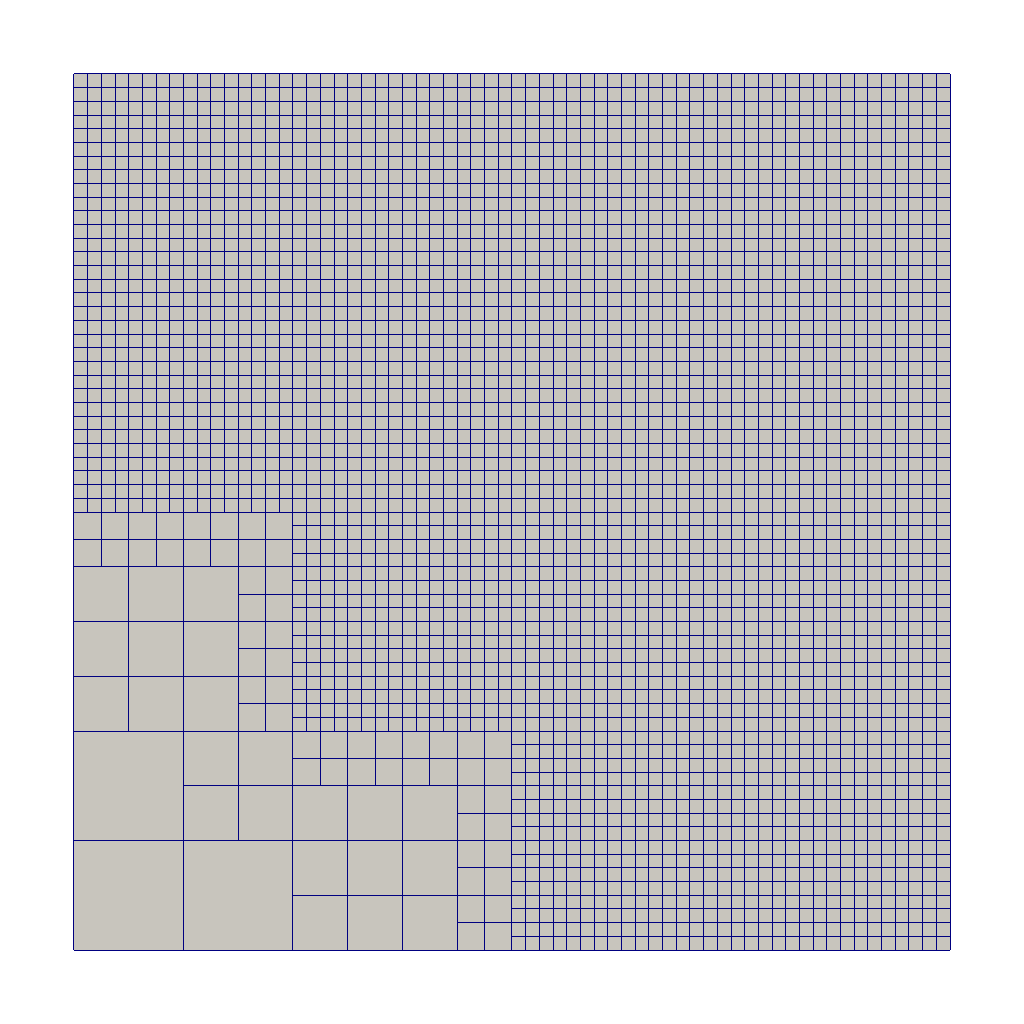}}
    \subfigure[Adapted mesh: $i=3$.]{\includegraphics[width=0.4\textwidth]{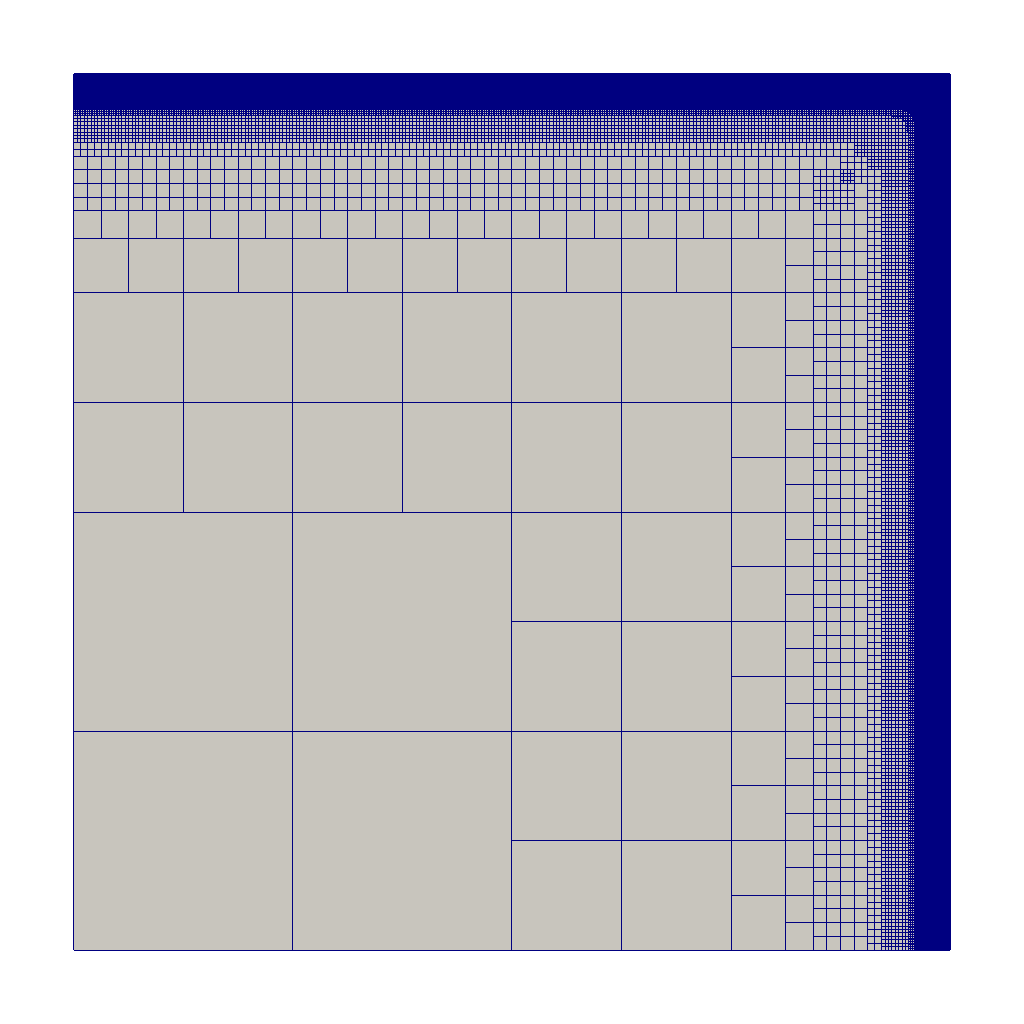}}
    \caption{Test case 1. Metric-based adaptation.}
    \label{fig:dr1_metrics}
\end{figure}

We solve problem \eqref{eq:adeq} 
in the domain \(\Omega = (0, 1)^2\)
when dealing with constant coefficients. In particular, we set \(\varepsilon = 10^{-4}\),
$\beta=0$, \(b=1\), while $f$ and $g$ are chosen such that the exact solution is 
\begin{equation*}
u_{\mathrm{ex}}(x, y) = \left(1 - \frac{\sinh(x / \sqrt{\varepsilon})}{\sinh(1 / \sqrt{\varepsilon})}\right)\left(1 - \frac{\sinh(y / \sqrt{\varepsilon})}{\sinh(1 / \sqrt{\varepsilon})}\right),
\end{equation*}
and $\Gamma_D$ coincides with the whole boundary $\partial \Omega$.
Solution $u_{\mathrm{ex}}$ exhibits a boundary layer both along the top and the right side of $\Omega$.\\
Algorithms~\ref{algo:marking} and~\ref{algo:metrics} are run starting from the same initial mesh, i.e., a uniform mesh consisting of 16 squared cells, and by setting \(i_{\mathrm{max}} = 10\) and \(\mathrm{tol} = 10^{-5}\). Algorithm~\ref{algo:marking} converges after $9$ iterations, while the metric-based procedure demands only $3$ iterations.
Figures~\ref{fig:dr1_marker} and~\ref{fig:dr1_metrics} show the evolution of the mesh throughout the two different adaptive procedures. 
The boundary layers are correctly refined by both algorithms. However, it is evident that the metric-based approach takes less iterations ($3$ adaptations are enough to capture the two layers in contrast to the $9$ iterations demanded by Algorithm~\ref{algo:marking}) and detects the two layers with a higher sharpness, as highlighted by the thinner refined areas.

\begin{figure}
    \centering
\includegraphics[width=0.45\textwidth]{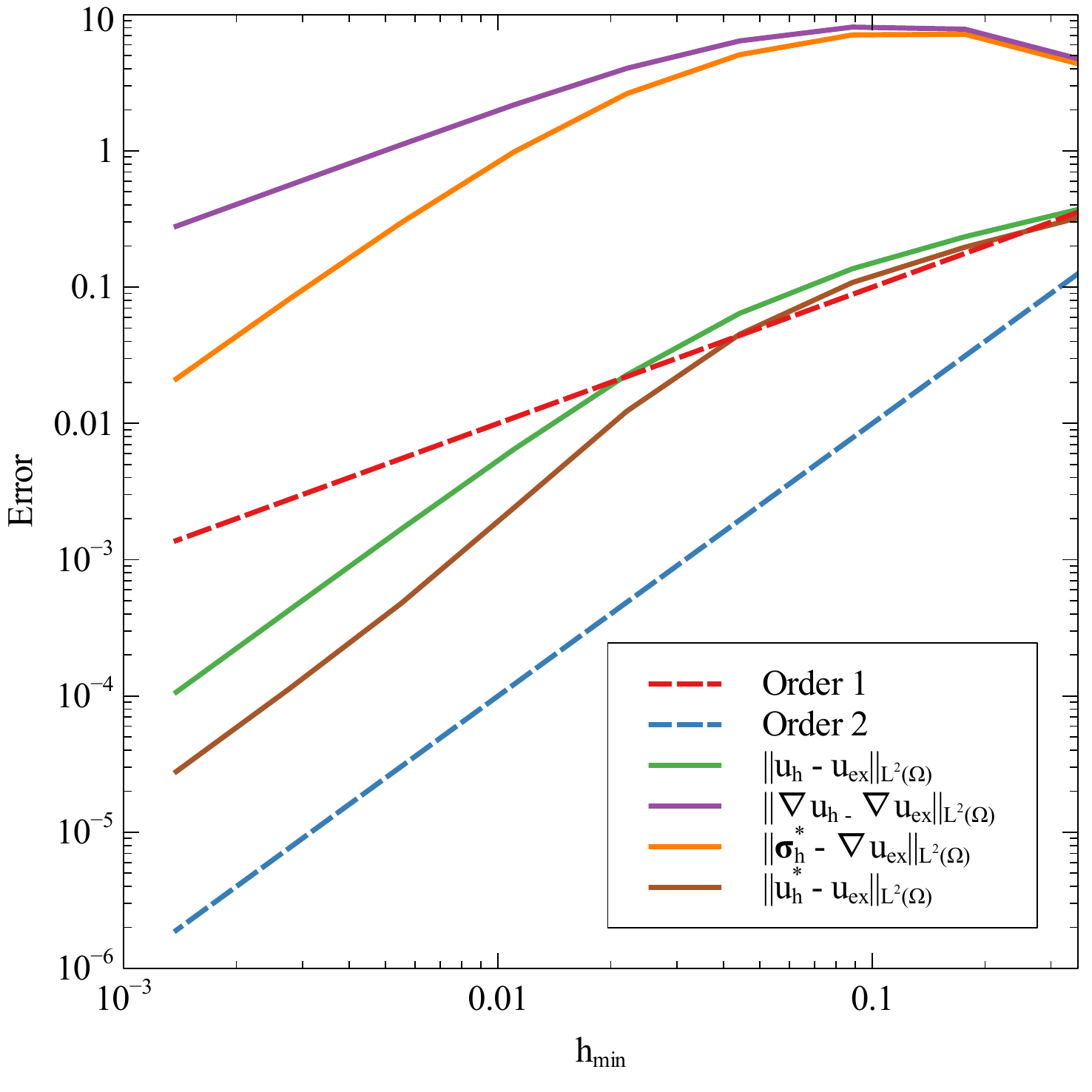}\hspace*{0.8cm}
\includegraphics[width=0.45\textwidth]{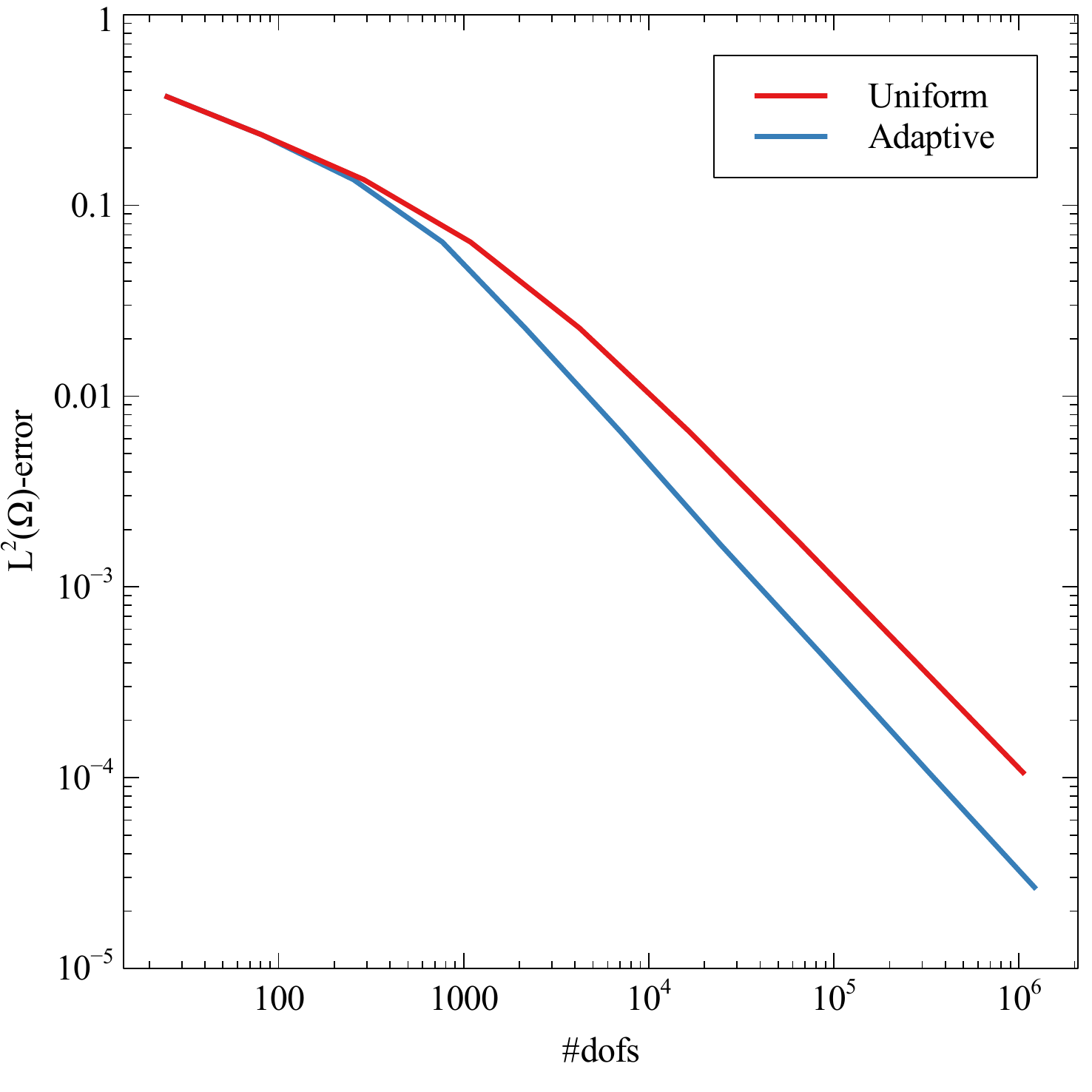}
    \caption{Test case 1. Convergence analysis of the  error associated with the discrete and the recovered  solution and gradient, as a function of the mesh size $h_{\mathrm{min}}$ (left); log-log plot of the norm \(\norm{u - u_h}_{L^2(\Omega)}\) as a function of the number of dofs on a uniform and on a metric-based adapted mesh (right).}
    \label{fig:dr1_convergence}
\end{figure}

\begin{figure}
    \centering
\includegraphics[width=0.45\textwidth]{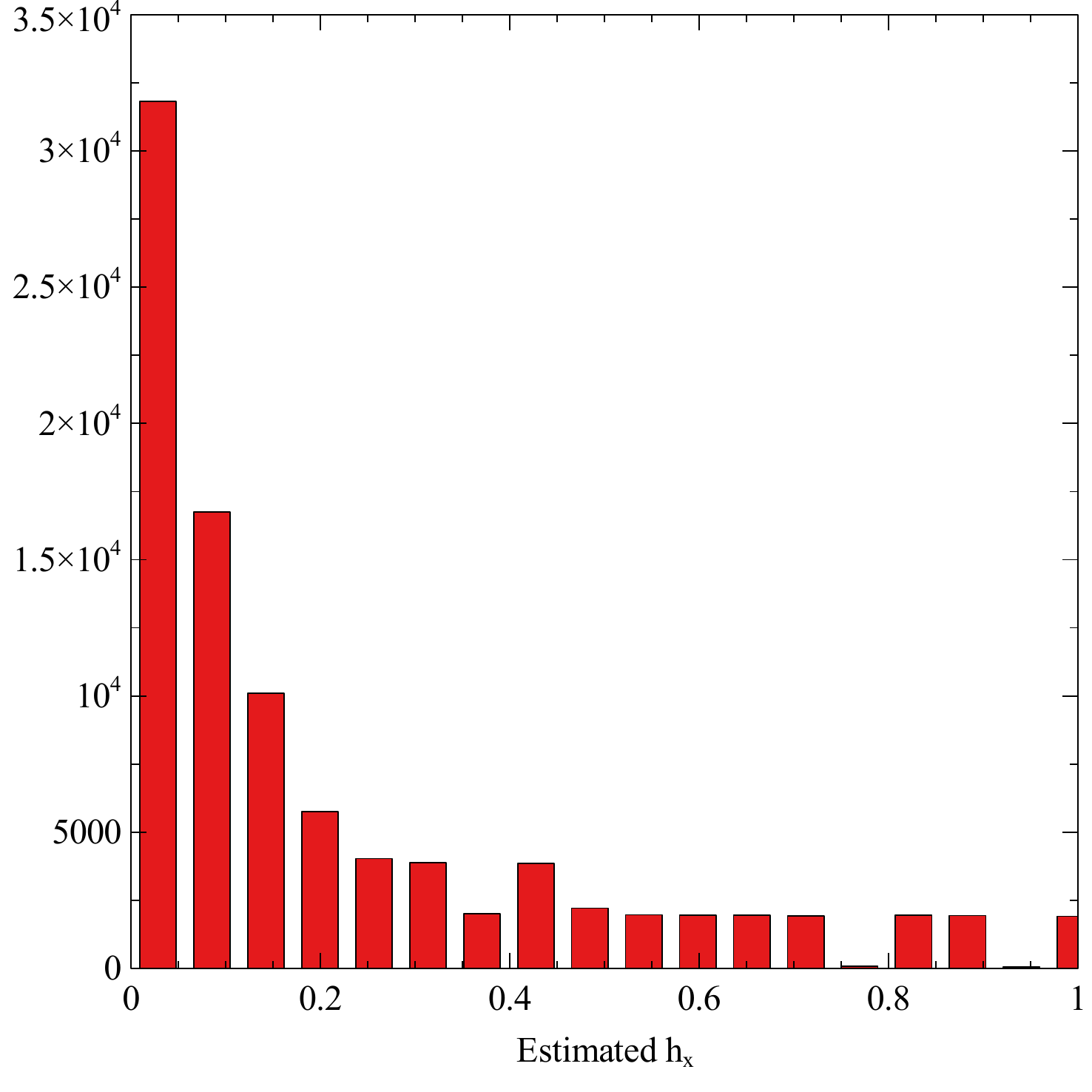}\hspace*{0.8cm}
\includegraphics[width=0.45\textwidth]{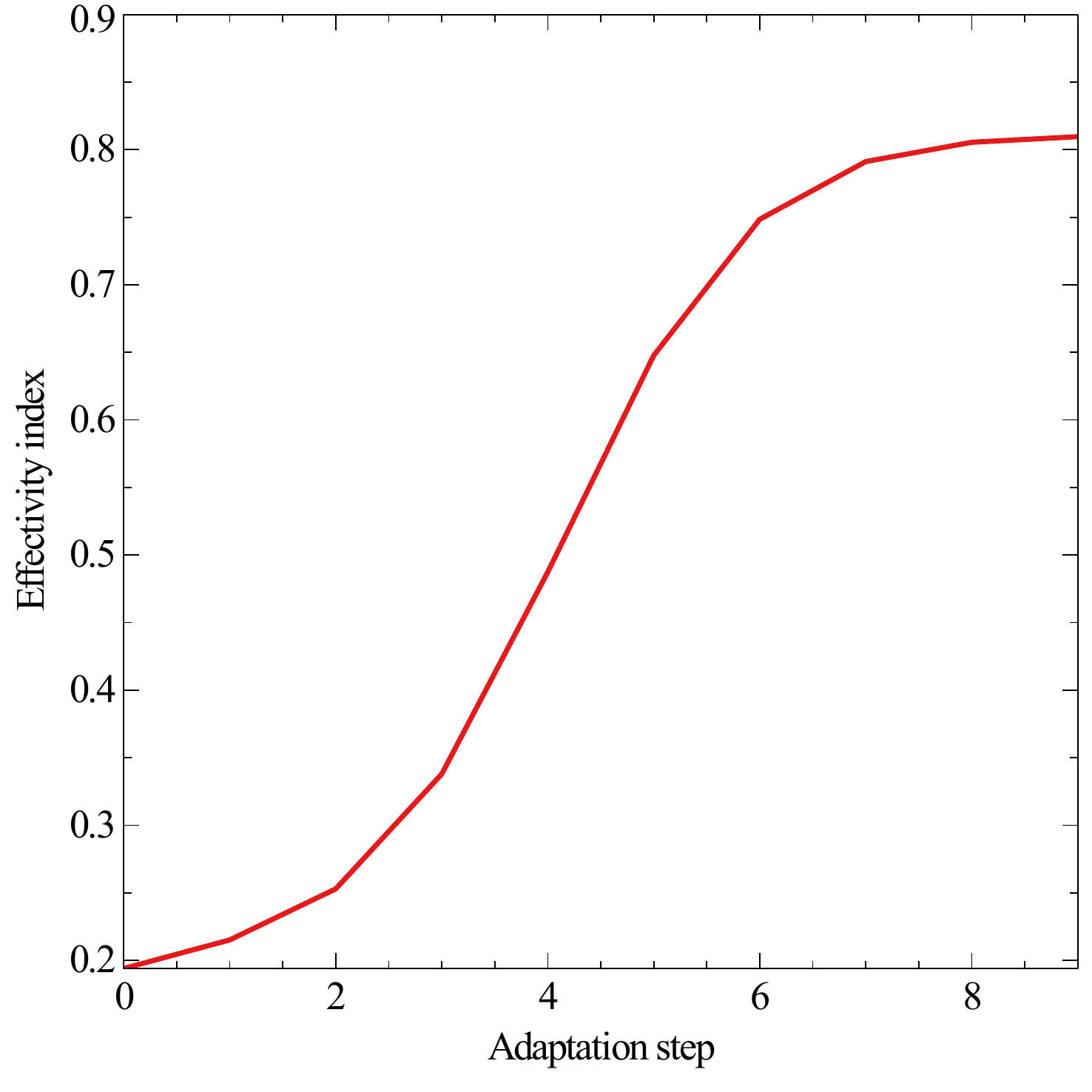}
    \caption{Test case 1. Histogram of the element size computed at the last iteration of the metric-based adaptive procedure (left); trend of the effectivity index as a function of the adaptation step (right).}
    \label{fig:dr1_effectivity}
\end{figure}

In Figure~\ref{fig:dr1_convergence} (left) we investigate the approximation property of the recovered quantities $u_h^*$ and $\bm{\sigma}_h^*$. In particular, we compare the convergence of the $L^2(\Omega)$-norm of the discretization errors, $\left(u_{\mathrm{ex}}-u_h\right)$ and $\left(\nabla u_{\mathrm{ex}} - \nabla u_h\right)$, as a function of the minimum mesh size, $h_{\mathrm{min}}$, with the trend of the quantities, $\left(u_{\mathrm{ex}}- u_h^*\right)$ and $\left(\nabla u_{\mathrm{ex}} - \bm{\sigma}_h^*\right)$, respectively. The expected order of convergence is guaranteed, at least asymptotically, by all the errors, with a superconvergent behaviour for the error associated with the recovered gradient. 
Additionally, it turns out that both the recovered gradient and solution provide a better approximation when compared with the corresponding discrete quantities, thus confirming the reliability of the recovered procedures set in Sections~\ref{grad_rec_sec} and~\ref{recsol_sec}. 
\\The right panel in Figure~\ref{fig:dr1_convergence} highlights the gain in terms of dofs ensured by the metric-based adaptive procedure when compared to a uniform refinement of $\tau_h$. For a sufficiently large number of dofs, about half of the unknowns required by  
a uniform refinement process
is demanded by the metric-based adaptive algorithm to ensure a certain accuracy on the discrete solution. In particular, in Figure~\ref{fig:dr1_effectivity} (left), we show the 
histogram of the elements size predicted by the metric at the last iteration of the metric-based adaptive procedure.

Furthermore, we verify the robustness of the error estimator $\eta$ by computing the effectivity index
\begin{equation*}
\xi = \frac{\eta}{\norm{u - u_h}_{L^2(\Omega)}},
\end{equation*}
with an ideal value equal to $1$. The value of $\xi$ 
at each adaptation step is provided in Figure~\ref{fig:dr1_effectivity} (right). We observe that $\eta$ slightly under-estimates the exact error, with a scaling factor between $0.8$ and $1$. Nevertheless, $\xi$ exhibits a stagnation trend across the adaptive iterations, which confirms the robustness of the proposed estimator.

\begin{figure}
{\includegraphics[width=0.45\textwidth]{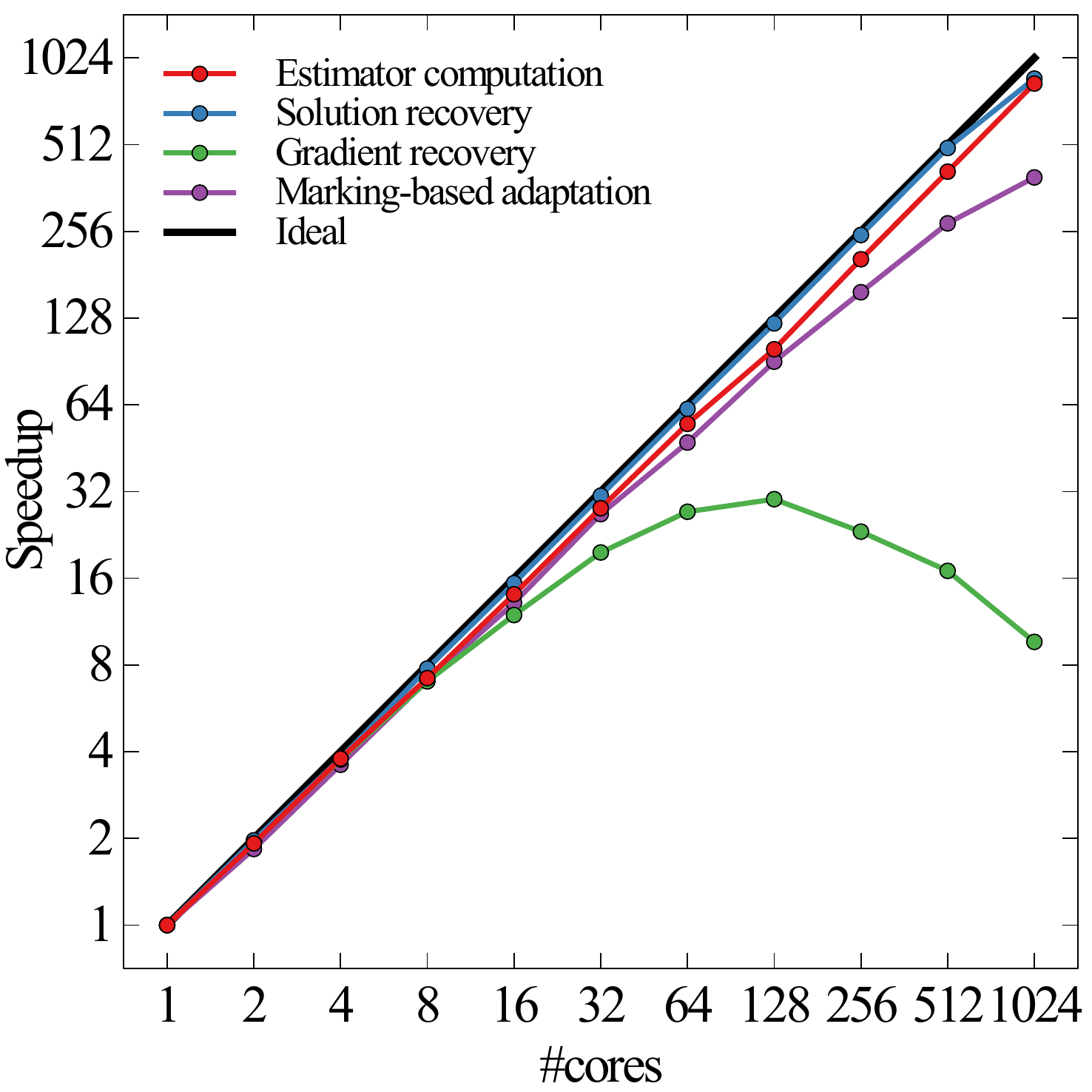}}\hspace*{0.8cm}
{\includegraphics[width=0.45\textwidth]{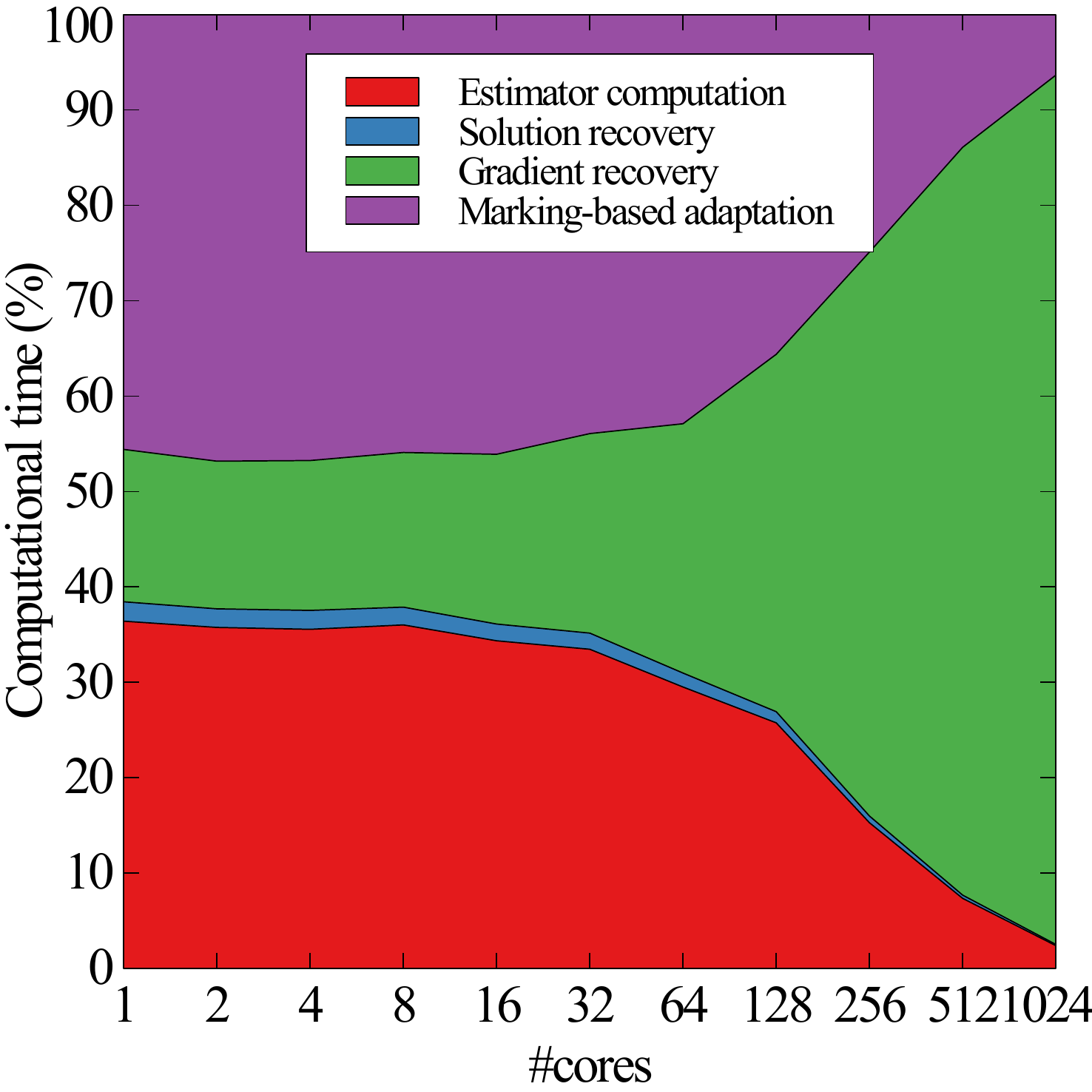}}
\caption{Test case 1. Scalabilty analysis of the
different phases of the marking-based adaptation procedure:
parallel speedup (left) and stacked barplot of the time percentage associated with the different phases of the algorithm (right) as a function of the number of parallel processes.}\label{fig:dr1_marker_speedup}
\end{figure}

In order to assess the scalability
performance of the marking-based adaptation procedure, we have run  \Cref{algo:marking} on a different number of parallel processes. The results are summarized in \Cref{fig:dr1_marker_speedup}. 
In the left panel, we show 
the speedup as a function of the number of the parallel processes, 
for the different phases of the algorithm, namely the recovery of the gradient \(\bm{\sigma}_h^*\) (as described in \Cref{grad_rec_sec}) and of the solution \(u_h^*\) (as described in \Cref{recsol_sec}), the computation of the error estimator \(\eta_k\) over each mesh element (as described in \Cref{sec:refinement}) and the mesh adaptation. 
All these phases exhibit an almost linear speedup, up to \(1024\) cores, with the only exception being the gradient recovery whose trend
deteriorates starting from \(128\) cores. This can be ascribed to an unbalanced load among processes due to the presence of the 
two very narrow boundary layers. Indeed, while the solution recovery can be performed independently on each cell, the gradient recovery involves an exchange of data among neighboring elements. 
The demanding computational footprint characterizing the recovery of the gradient is confirmed by the right panel in 
\Cref{fig:dr1_marker_speedup}, which displays the percentage of the computational time associated with each phase of the marking-based algorithm. Finally, we remark that 
the whole adaptation procedure requires a contained 
computational time (between $1$ and $50$s), thus providing an efficient tool to deal with more challenging scenarios, such as large-scale or time-dependent problems.

\subsection{Test case 2: diffusion-reaction problems with discontinuous data}\label{prob:dr3}
In this section we examine two configurations which exhibit a discontinuity, rectilinear and circular, respectively, The latter setting offers a challenging benchmark since dealing with 
a Cartesian mesh adaptation procedure. 

\subsubsection{Rectilinear discontinuity}
\begin{figure}
    \centering
{\includegraphics[width=0.45\textwidth]{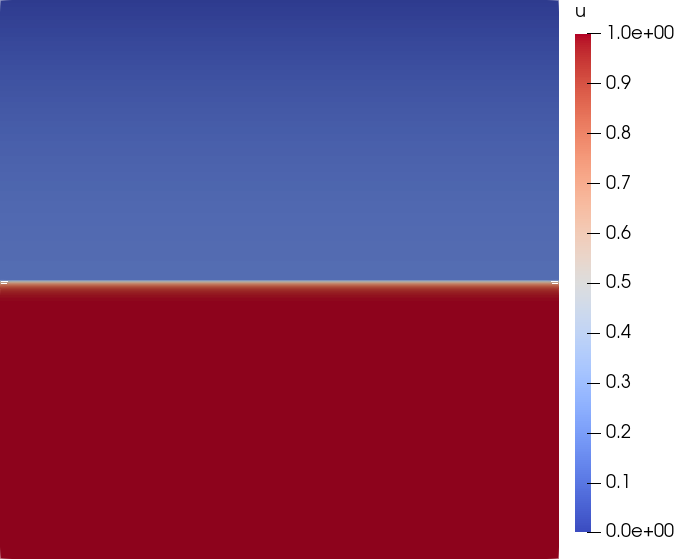}}\hspace*{0.8cm}
{\includegraphics[width=0.45\textwidth]{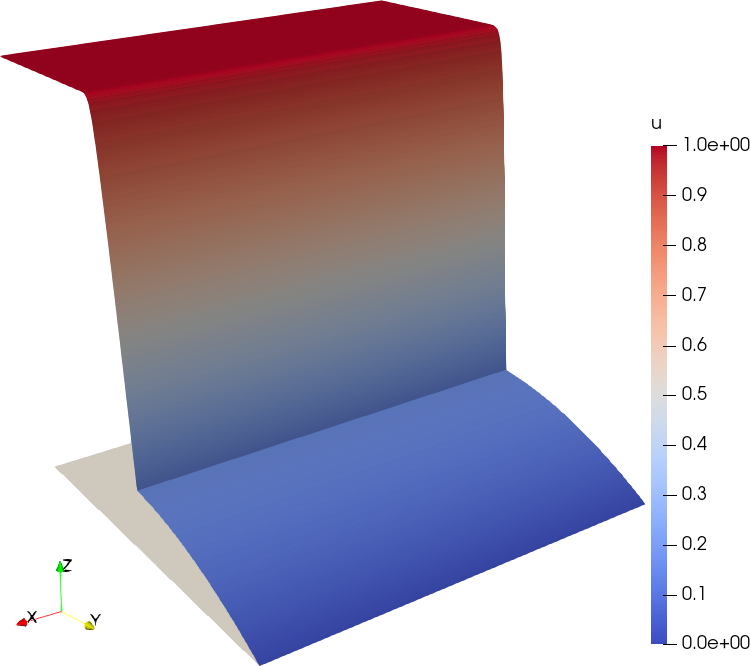}}
    \caption{Test case 2 (rectilinear discontinuity). 2D (left) and 3D (right) view of the exact solution.}
    \label{fig:dr3_sol}
\end{figure}

We consider a diffusion-reaction problem firstly proposed in \cite{de2010interior}, characterized by a discontinuous diffusion and reaction.
We solve problem \eqref{eq:adeq} in the unit square \(\Omega = (0, 1)^2\), by setting 
\begin{equation*}
\varepsilon = 
\left\{
\begin{alignedat}{2}
\varepsilon_1 & = 5 \cdot 10^{-5} & \quad & \text{in } \Omega_1 \\
\varepsilon_2 & = 10^{-1} & \quad & \text{in } \Omega_2,
\end{alignedat}
\right.
\end{equation*}
with \(\Omega_1 = \{ (x, y) \, :\, y \leq 0.5 \}\) and \(\Omega_2 = \Omega \setminus \Omega_1\), $\beta=0$,
\begin{equation*}
b = 
\left\{
\begin{aligned}
1 & \quad \text{in } \Omega_1 \\
0 & \quad \text{in } \Omega_2,
\end{aligned}
\right.
\qquad
f = 
\left\{
\begin{aligned}
1 & \quad \text{in } \Omega_1 \\
\varepsilon_2 & \quad \text{in } \Omega_2.
\end{aligned}
\right.
\end{equation*}
Concerning the boundary data, we select 
\begin{equation*}
g =
\left\{
\begin{alignedat}{3}
& 1 & \quad & \text{on } \Gamma_{D1} \\
& 0 & \quad & \text{on } \Gamma_{D2},
\end{alignedat}
\right.
\end{equation*}
with $\Gamma_{D1}=\{(x,y)\, :\, y=0\}$, 
$\Gamma_{D2}=\{(x,y)\, :\, y=1\}$.
This choice of data leads us to identify function 
\begin{equation*}
u_\mathrm{ex}(x, y) =
\left\{
\begin{aligned}
1 + 2c_1 \sinh(y / \sqrt{\varepsilon_1}) & \quad \text{in } \Omega_1 \\
-0.5(y - 1)(y + 2c_2) & \quad \text{in } \Omega_2,
\end{aligned}
\right.
\end{equation*}
with the exact solution to the problem, for
\begin{equation*}
\begin{aligned}
    c_1 & = -\frac{7\varepsilon_2}
    {8\sqrt{\varepsilon_1} \cosh(0.5 / \sqrt{\varepsilon_1}) + 16\varepsilon_2 \sinh(0.5 / \sqrt{\varepsilon_1})} \\[2mm]
    c_2 & = \frac{7\sqrt{\varepsilon_1} \cosh(0.5 / \sqrt{\varepsilon_1})}
    {4\sqrt{\varepsilon_1} \cosh(0.5 / \sqrt{\varepsilon_1}) + 8\varepsilon_2 \sinh(0.5 / \sqrt{\varepsilon_1})}.
\end{aligned}
\end{equation*}
Function $u_\mathrm{ex}$ exhibits a strictly one-dimensional dynamics along the $y$-direction, together with a jump across the line $\{ (x,y)\, :\, y=0.5\}$ (see Fig.~\ref{fig:dr3_sol}).

\begin{figure}
    \centering
    \subfigure[Adapted mesh: $i=1$.]{\includegraphics[width=0.4\textwidth]{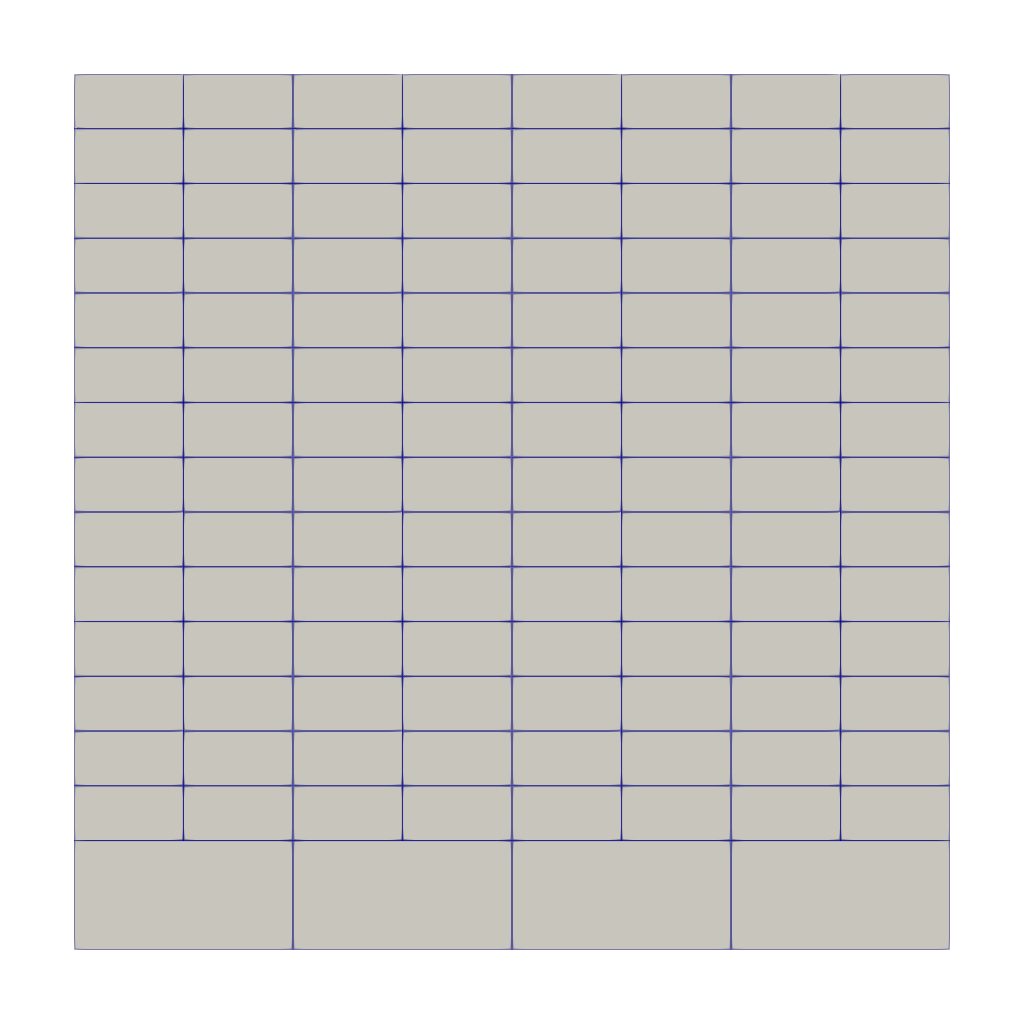}}
    \subfigure[Adapted mesh: $i= 3$.]{\includegraphics[width=0.4\textwidth]{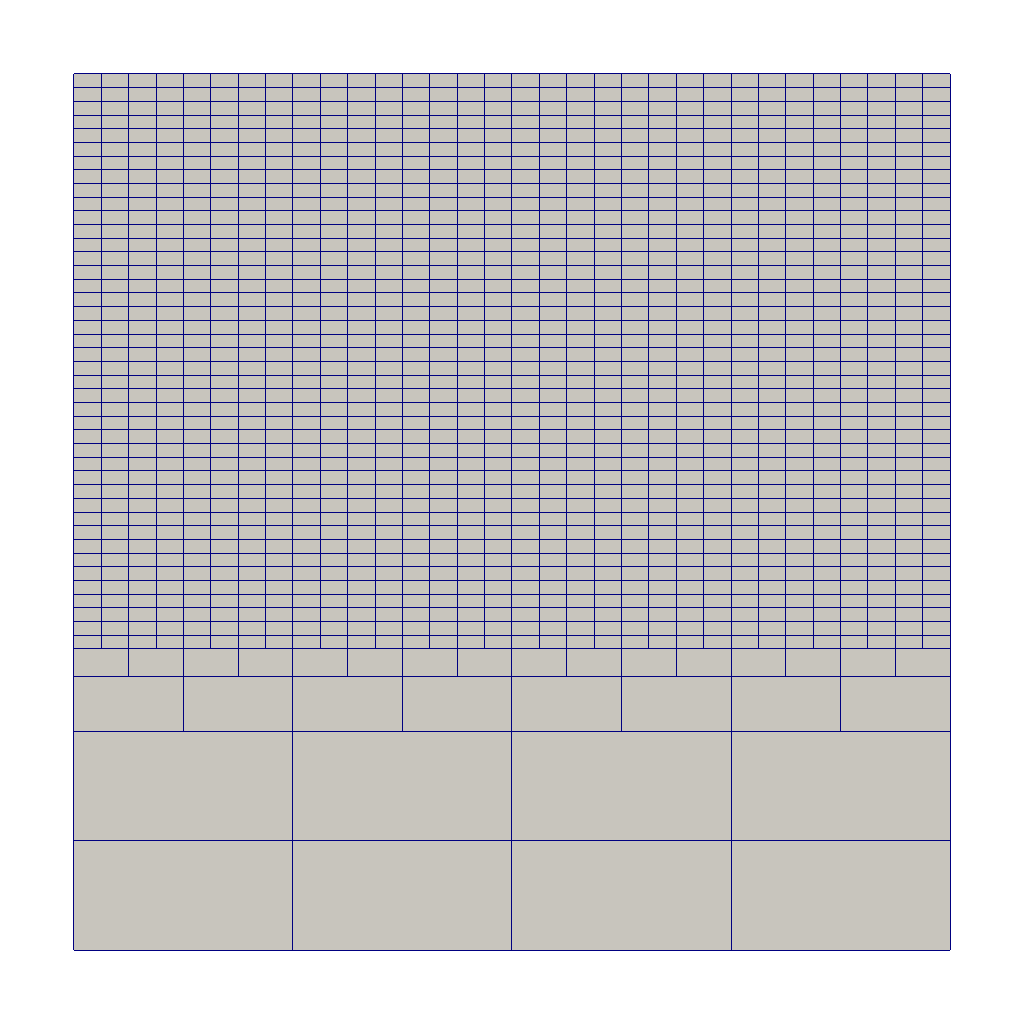}} \\
    \subfigure[Adapted mesh: $i=5$.]{\includegraphics[width=0.4\textwidth]{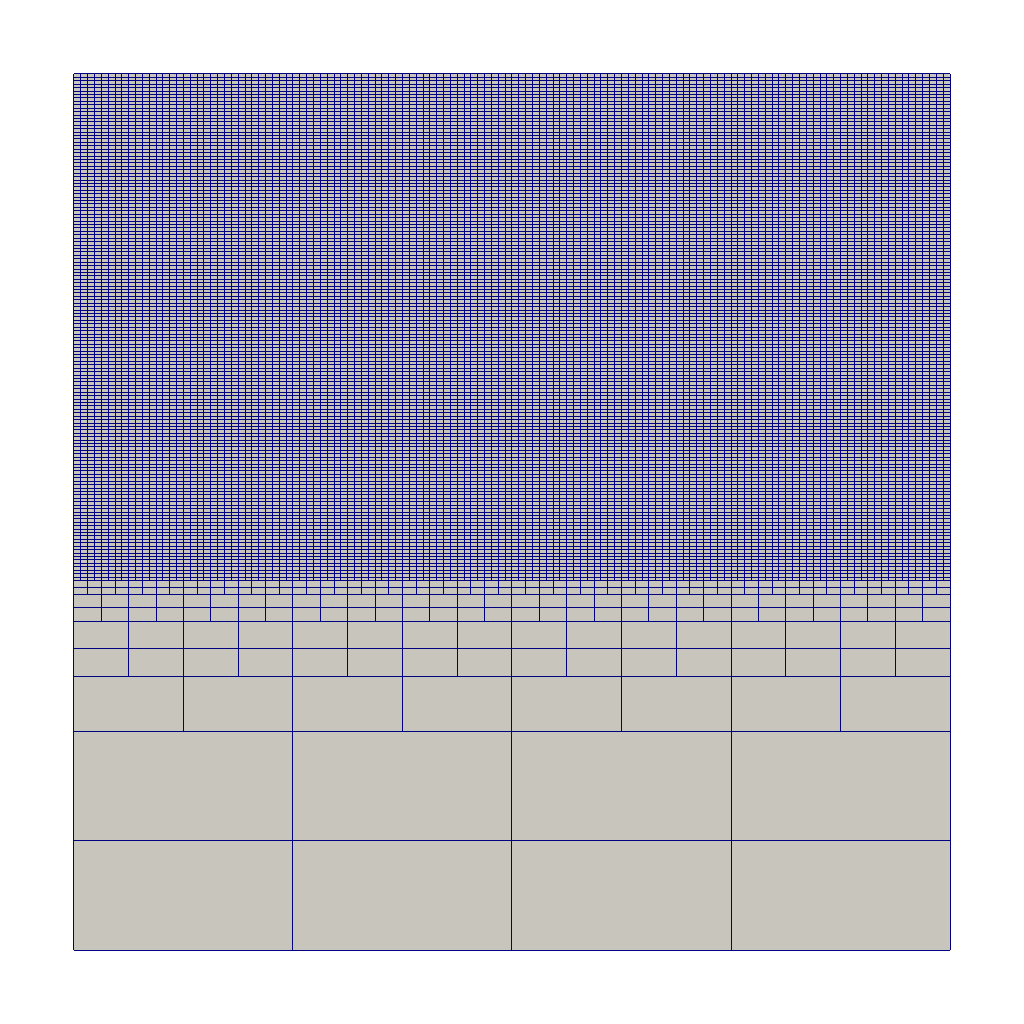}}
    \subfigure[Adapted mesh: $i= 9$.]{\includegraphics[width=0.4\textwidth]{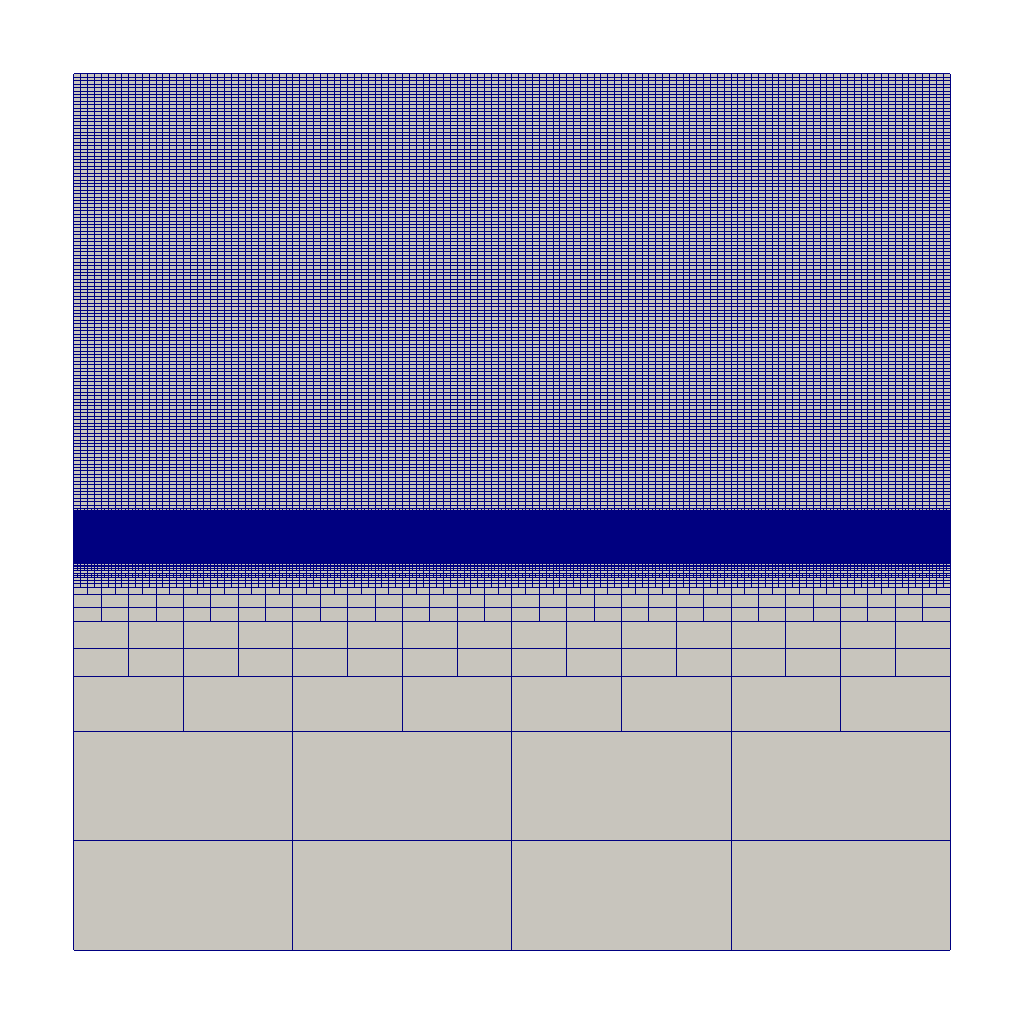}}
    \caption{Test case 2 (rectilinear discontinuity). Marking-based adaptation.}
    \label{fig:dr3_marker}
\end{figure}

\begin{figure}
    \centering
    \subfigure[Adapted mesh: $i= 1$.]{\includegraphics[width=0.4\textwidth]{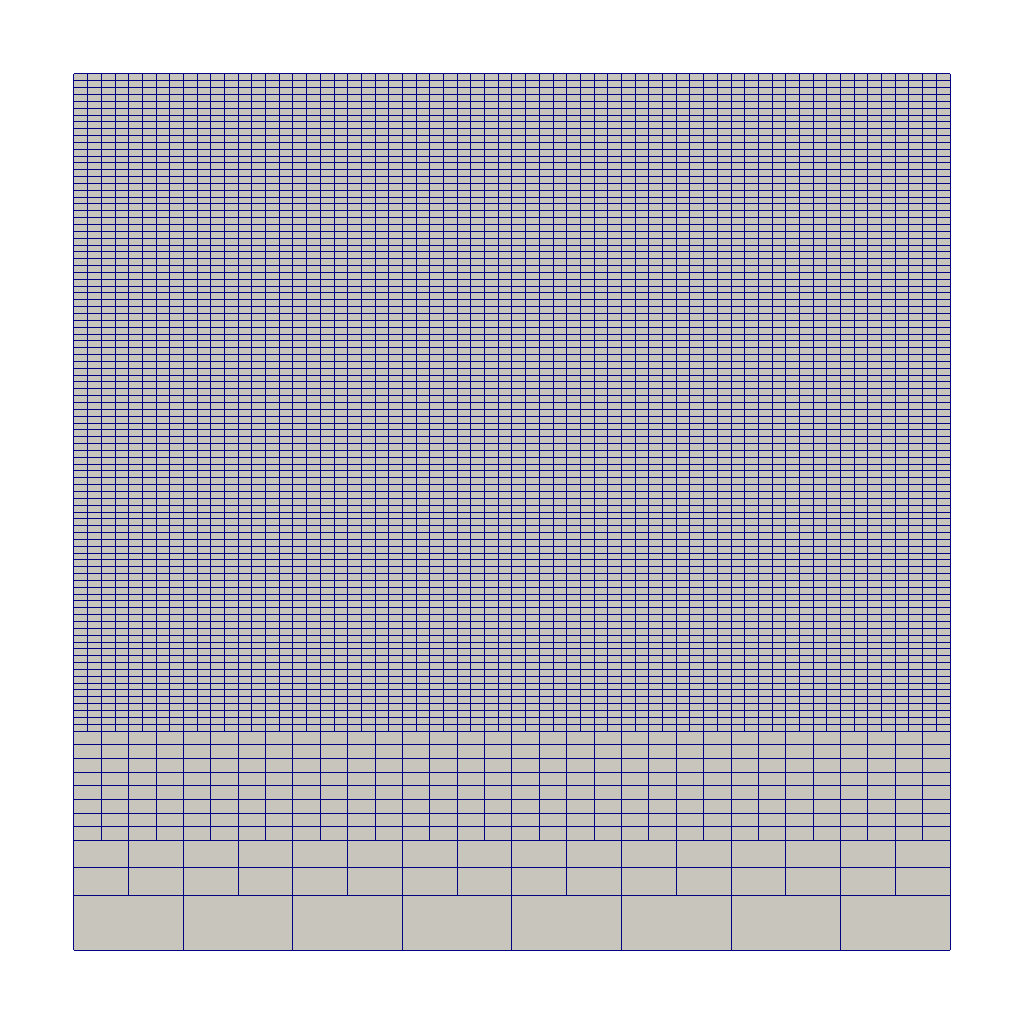}}
    \subfigure[Adapted mesh: $i= 3$.]{\includegraphics[width=0.4\textwidth]{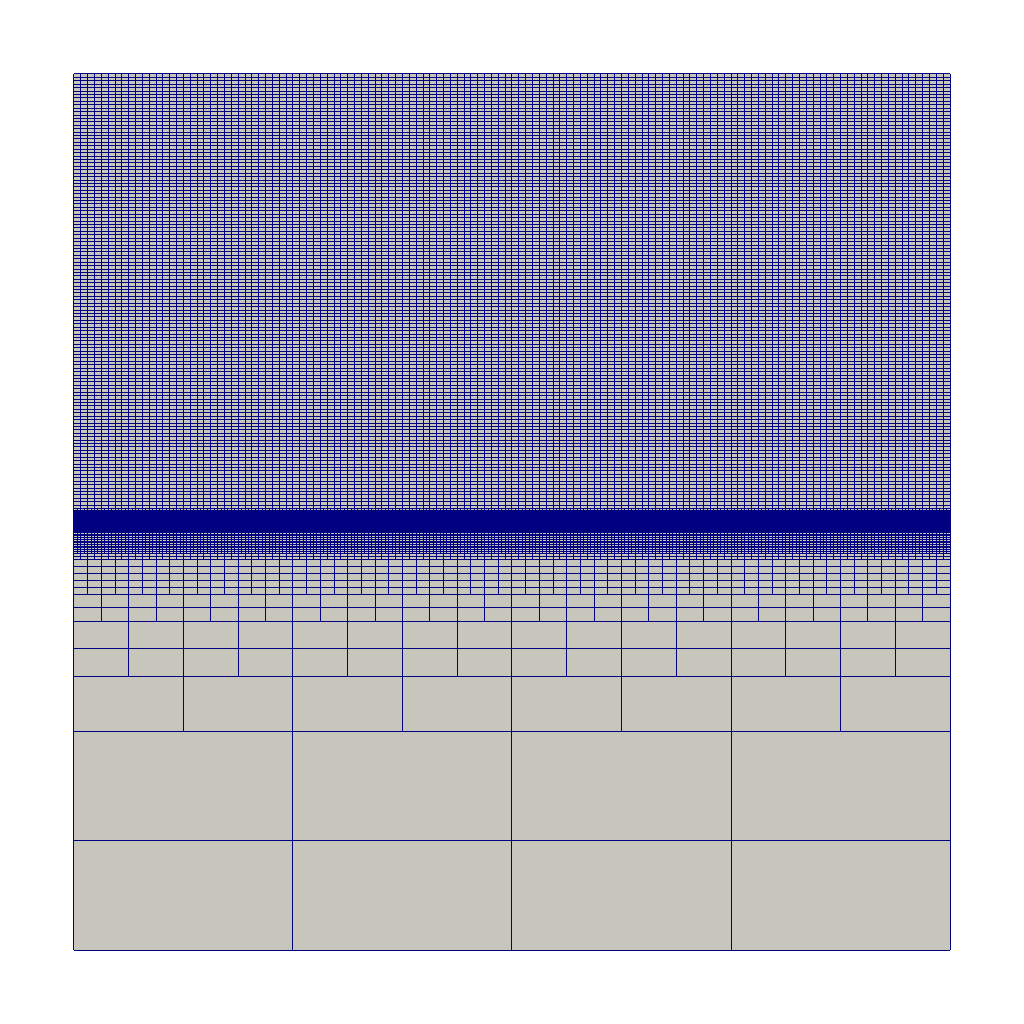}}
    \caption{Test case 2 (rectilinear discontinuity). Metric-based adaptation.}
    \label{fig:dr3_metrics}
\end{figure}

We resort to both Algorithms~\ref{algo:marking} and~\ref{algo:metrics} to 
approximate $u_\mathrm{ex}$ in order to compare the performance of the two adaptive procedures. To this aim, we pick a uniform initial mesh characterized by $32$ rectangular cells ($4$ and $8$ subdivisions along the \(x\)- and the \(y\)-direction, respectively), while assigning values $10^{-6}$ and $10$ to the input parameters \(\mathrm{tol}\) and \(i_{\mathrm{max}}\), respectively. \\
The marking-based procedure takes $9$ iterations to converge (see Figure~\ref{fig:dr3_marker}), while the metric-based algorithm breaks after only $3$ iterations by delivering the adapted mesh shown in Figure~\ref{fig:dr3_metrics} (right).
Both the adaptive procedures start by refining the upper part of the domain; successively, the refinement is further increased in correspondence with the horizontal layer, while the mesh elements are coarsened in the bottom portion of $\Omega$. In particular, we observe that, at the first iteration, the metric-based approach more massively refine the upper part of $\Omega$ with respect to the marking-based algorithm. This is compliant with the quickest convergence of Algorithm~\ref{algo:metrics}, analogously to Test case $1$. Moreover, the layer turns out to be more sharply captured by the metric-based approach.

\begin{figure}
    \centering  \includegraphics[width=0.45\textwidth]{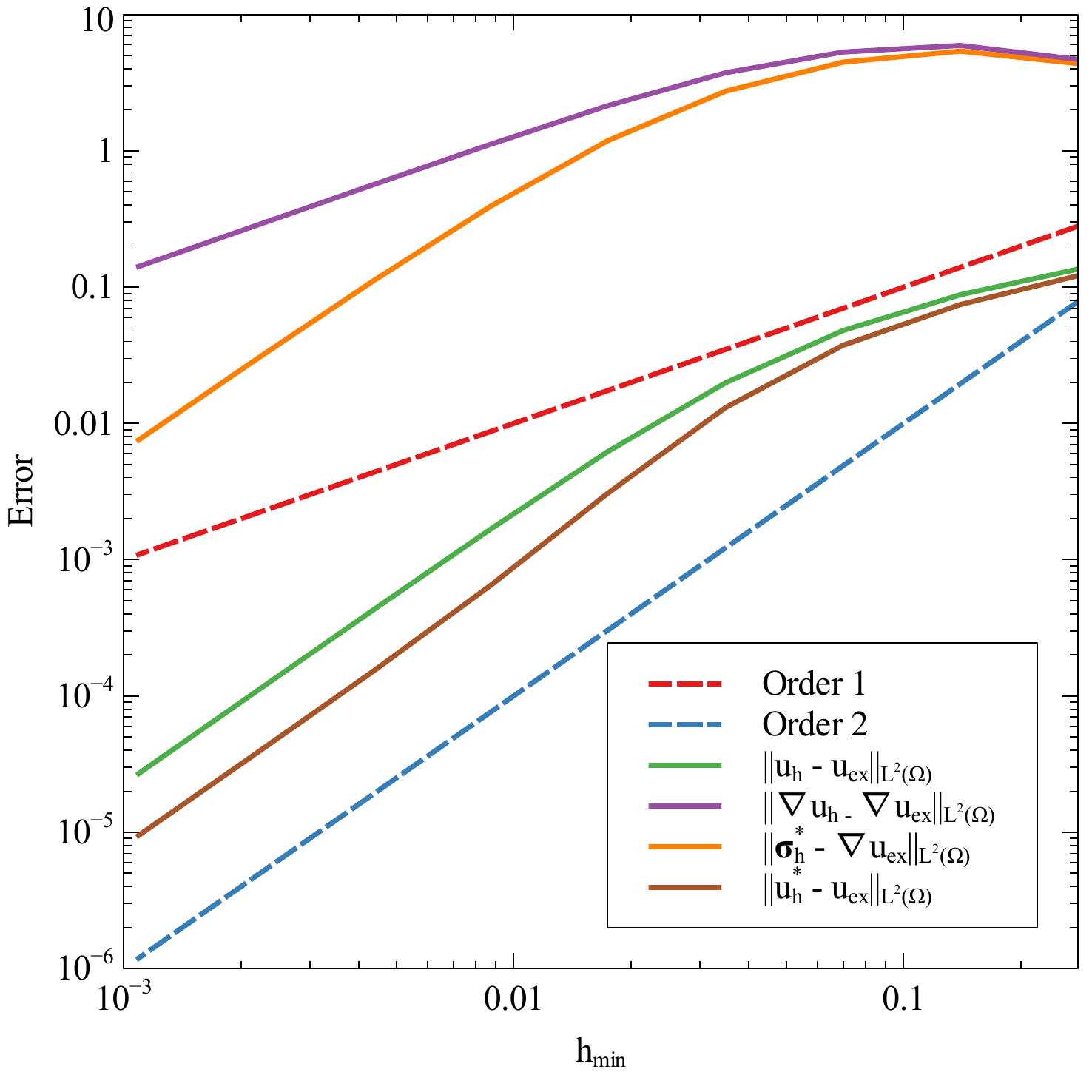}\hspace*{0.8cm}  \includegraphics[width=0.45\textwidth]{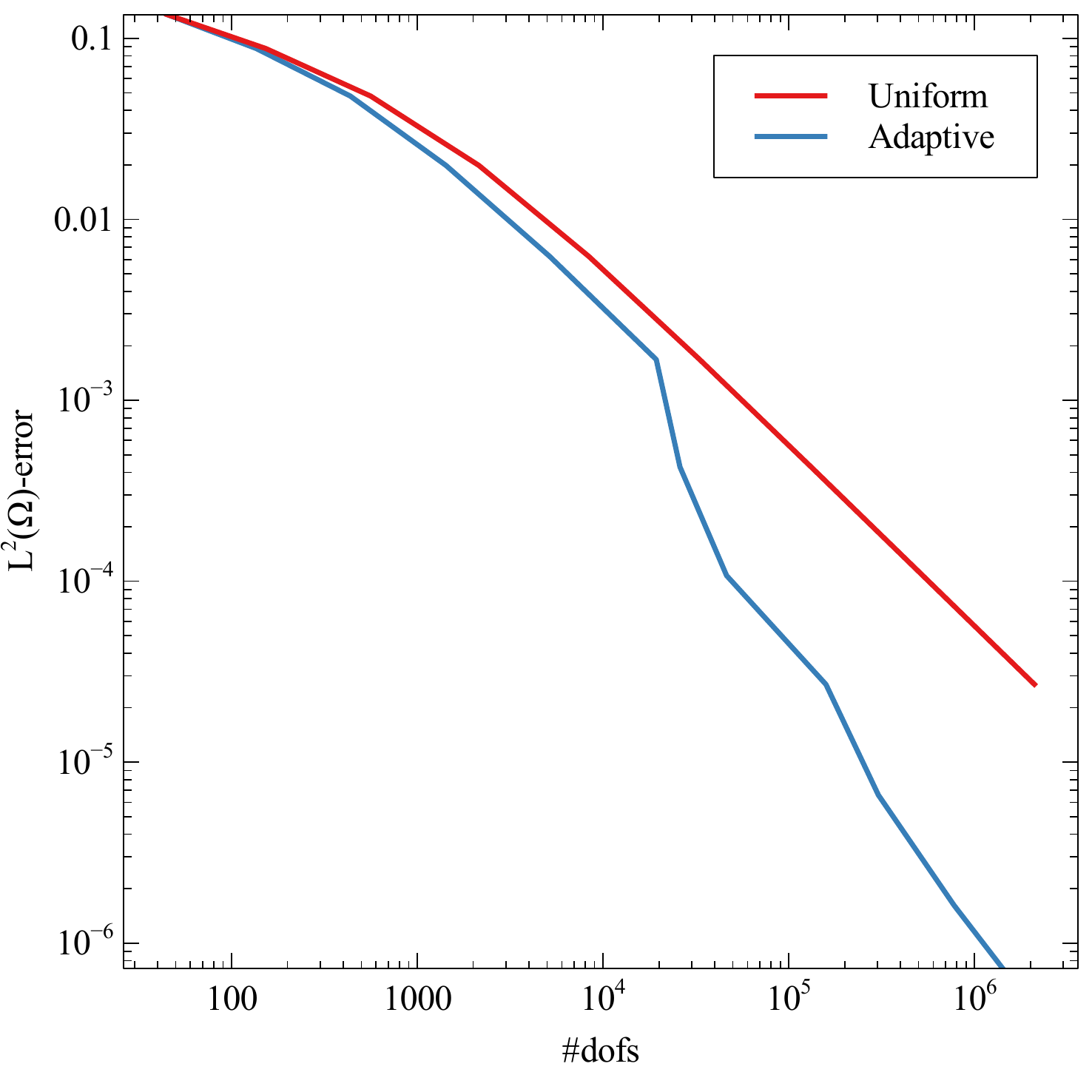}
    \caption{Test case 2 (rectilinear discontinuity). 
    Convergence analysis of the  error associated with the discrete and the recovered  solution and gradient, as a function of the mesh size $h_{\mathrm{min}}$ (left); log-log plot of the norm \(\norm{u - u_h}_{L^2(\Omega)}\) as a function of the number of dofs on a uniform and on a metric-based adapted mesh (right).}
    \label{fig:dr3_convergence}
\end{figure}

\begin{figure}
    \centering
\includegraphics[width=0.45\textwidth]{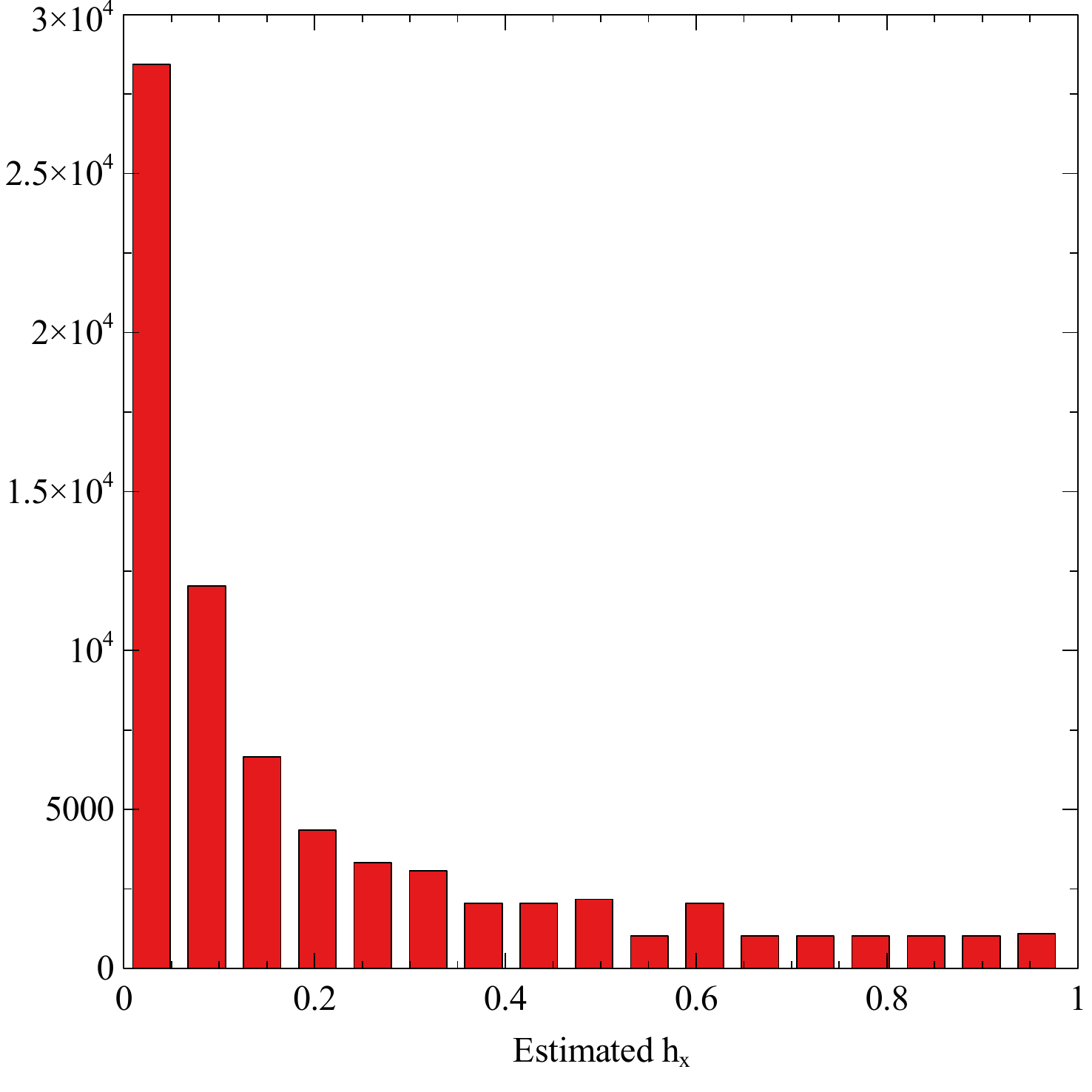}\hspace*{0.8cm}
\includegraphics[width=0.45\textwidth]{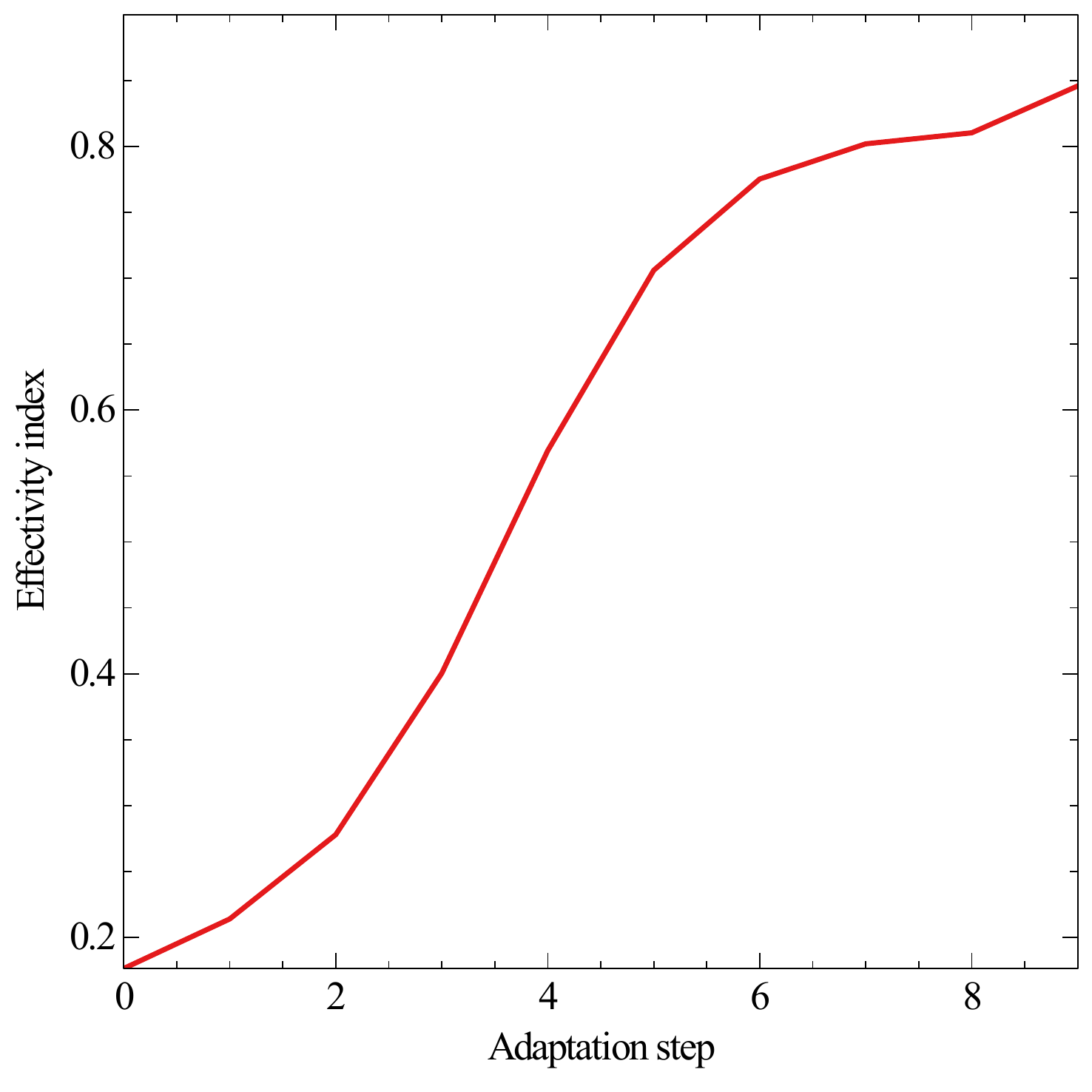}
    \caption{Test case 2 (rectilinear discontinuity). 
    Histogram of the element size computed at the last iteration of the metric-based adaptive procedure (left); trend of the effectivity index as a function of the adaptation step (right).}
    \label{fig:dr3_effectivity}
\end{figure}

The quantitative analysis in Figures~\ref{fig:dr1_convergence} and~\ref{fig:dr1_effectivity} is replicated for the configuration at hand. The convergence analysis in Figure~\ref{fig:dr3_convergence} (left) confirms the high reliability of the recovered solution, $u^*$, and gradient, $\bm{\sigma}_h^*$, which better approximate $u$ and $\nabla u$ when compared with $u_h$ and $\nabla u_h$, respectively. The orders of convergence are preserved, with a significant superconvergence exhibited by the recovered gradient. Concerning the computational advantages lead by a metric-based adapted mesh in terms of number of dofs, the discrepancy with a standard uniform refinement is more remarkable in this case with respect to what observed for Test case 1 (compare the right panels in Figures~\ref{fig:dr1_convergence} and~\ref{fig:dr3_convergence}). The same accuracy on the solution is guaranteed by more than one order of dofs less, when the number of unknowns is sufficiently large (about $10^5$ dofs). \\
Figure~\ref{fig:dr3_effectivity} (left) shows the histogram for the element size distribution predicted by Algorithm~\ref{algo:metrics}. 
By comparing this trend with the corresponding one in Figure~\ref{fig:dr1_convergence}, we observe that 
values are more concentrated around zero in this test case, the refined cells covering a larger part of the domain. Finally, the effectivity index in Figure~\ref{fig:dr3_effectivity} (right) exhibits a trend similar to the one in Figure~\ref{fig:dr1_effectivity} (right), although the stagnation trend is less evident at the last iterations.

\subsubsection{Circular discontinuity}
\label{sec:penalty}
Discontinuities across not rectilinear layers constitute a challenging task for Cartesian grids. The adoption of non-conforming meshes turns out to be instrumental in such a direction, since these grids allow to confine the refinement of the elements where it is strictly necessary. 

We assess the skills of the metric-based algorithm to deal with not straight discontinuities by considering a setting where the diffusion coefficient is discontinuous along a circular interface~\cite{raeli2018finite}. We solve a purely diffusive problem \eqref{eq:adeq} in $\Omega = (0, 1)^2$ after choosing
\begin{equation}
\label{eq:epsilon}
\varepsilon = 
\left\{
\begin{alignedat}{2}
& \varepsilon_G = 1   & \quad & \text{in } G \\
& \varepsilon_S = 100 & \quad & \text{in } S,
\end{alignedat}
\right.
\end{equation}
where \(G=\{(x-0.5)^2 + (y-0.5)^2 \geq R^2\}\) for \(R=0.25\), \(S=\Omega \setminus G\), $\beta=0$,
\(b=0\), \(f=1\). The exact solution coincides with the function
\begin{equation*}
u_\mathrm{ex}(x, y) =
\left\{
\begin{alignedat}{2}
& \frac{1}{8} - \frac{1}{4\varepsilon_G} \left((x-0.5)^2+(y-0.5)^2\right) \quad & \text{in } G \\
& \frac{1}{8} - \frac{1}{4\varepsilon_S} \left((x-0.5)^2+(y-0.5)^2\right) - \frac{R^2}{4}\left(1 - \frac{1}{\varepsilon_S}\right)\quad & \text{in } S.
\end{alignedat}
\right.
\end{equation*}

\begin{figure}
    \centering
    \subfigure[]{\includegraphics[width=0.45\textwidth]{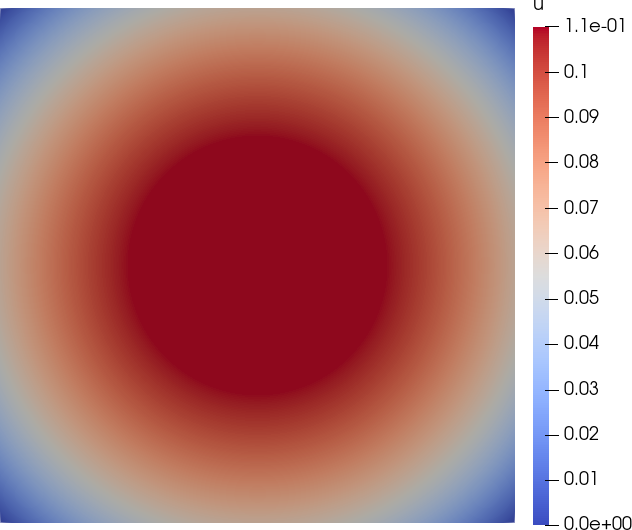}}
    \subfigure[]{\includegraphics[width=0.45\textwidth]{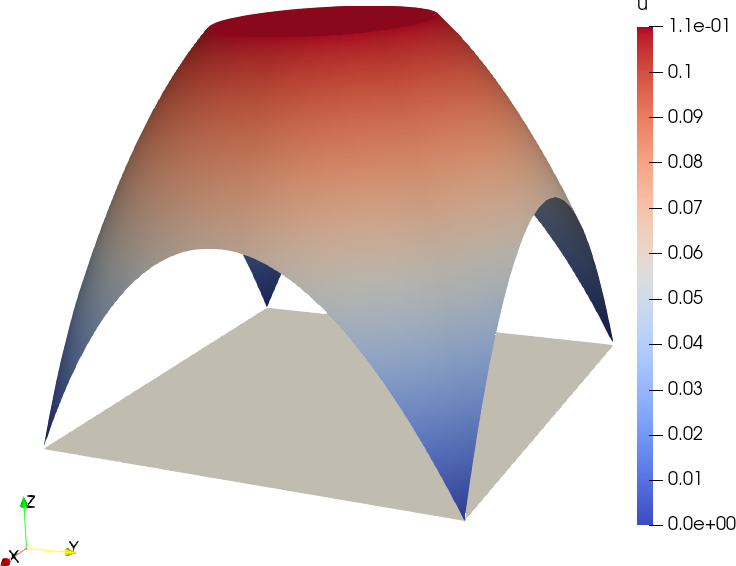}}
    \caption{Test case 2 (circular discontinuity). 2D (left) and 3D (right) view of the solution.}
    \label{fig:dr5_sol}
\end{figure}

We impose the Dirichlet boundary data \(u = u_\mathrm{ex}\) on \(\partial\Omega\) and a continuity condition at the interface \(\gamma\) between \(G\) and \(S\) on both the solution \(u\) and the associated normal flux \(\varepsilon\nabla u\).
Figure~\ref{fig:dr5_sol} provides a two- and three-dimensional view of the solution, which clearly highlights the significant variation of $u_\mathrm{ex}$ in $G$. Moreover, the penalization applied in $S$ mimics the presence of a hole centered in the domain.  

\begin{figure}
    \centering
    \hspace{-0.71cm}
    \subfigure[Adapted mesh: $i=2$.]{\includegraphics[height=0.35\textwidth]{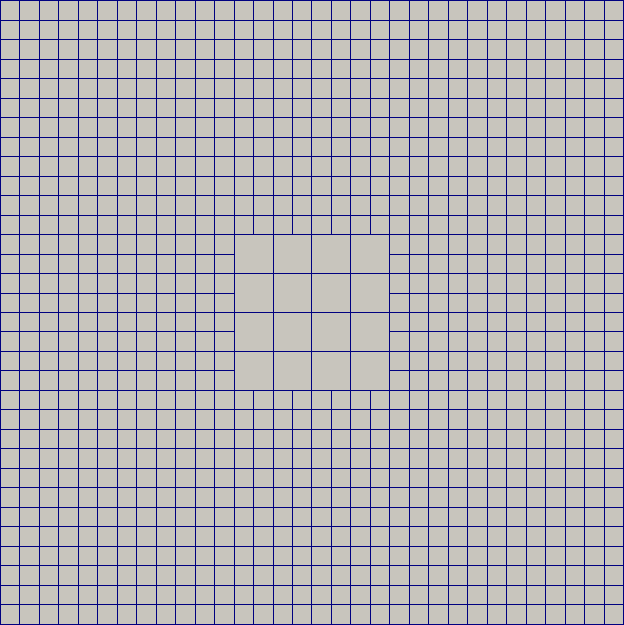}}
    \subfigure[Adapted mesh: $i=4$.]{\includegraphics[height=0.35\textwidth]{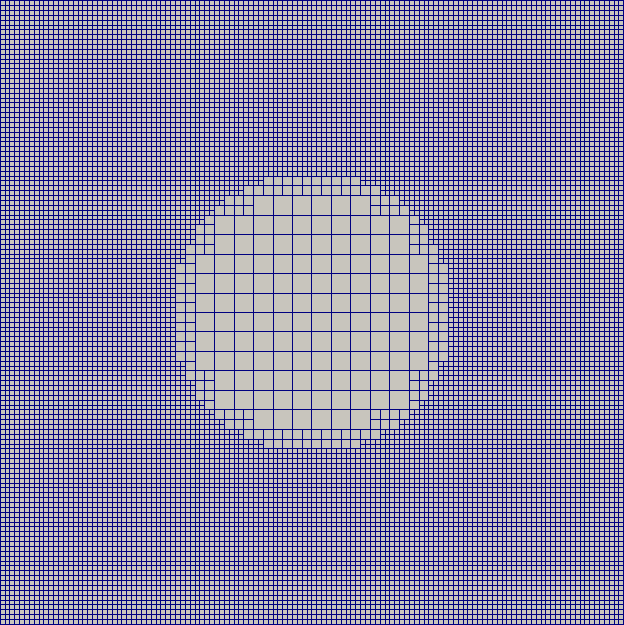}}\\
    \subfigure[Final mesh.]{\includegraphics[height=0.35\textwidth]{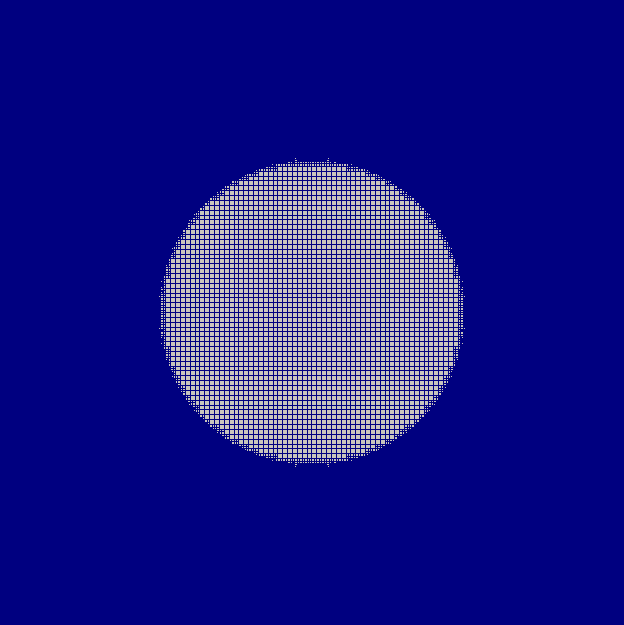}}
    \subfigure[Final mesh: detail of the interface.]{\includegraphics[height=0.35\textwidth]{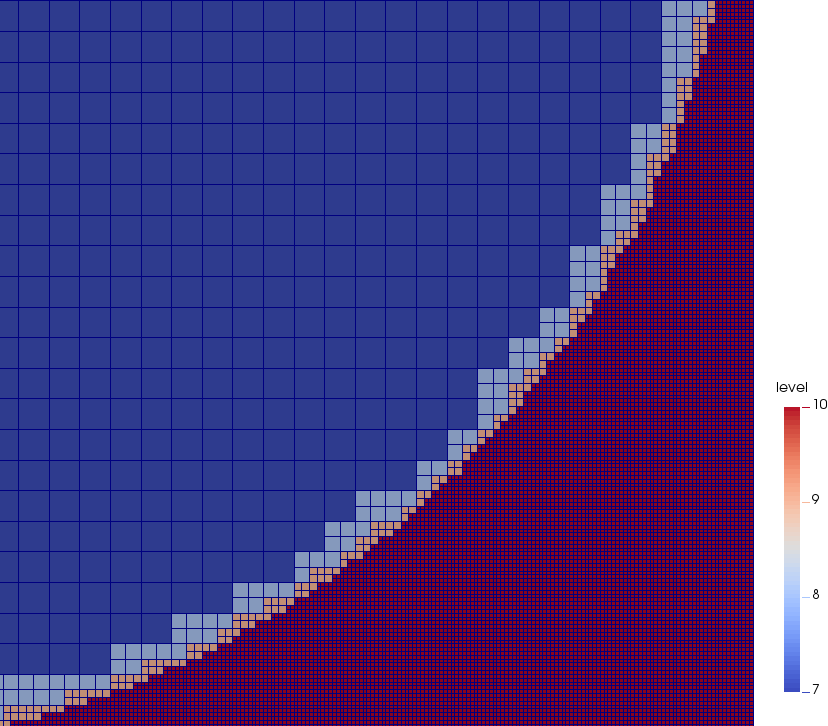}}
    \caption{Test case 2 (circular discontinuity). Metric-based adaptation.}
    \label{fig:dr5_marker}
\end{figure}

Algorithm~\ref{algo:metrics} is run on a uniform initial mesh with $8$ subdivisions along both the \(x\)- and \(y\)-direction, and by setting \(\mathrm{tol}=10^{-10}\) and \(i_{\mathrm{max}}=10\). The procedure returns a mesh with a maximum level of refinement equal to \(14\), which is shown in Figure~\ref{fig:dr5_marker} (bottom-left) together with two intermediate adapted meshes (top panels) and a zoom-in on the interface (bottom-right panel). 
We recognize a refinement of the grid along the interface between areas $S$ and $G$, consistently with the discontinuity of the diffusive coefficient. In addition, a
massive uniform refinement is applied to the region $G$; this is due to the fact that estimator $\eta$ 
tracks the solution, so that also 
a variation of $u_\mathrm{ex}$ different from a steep gradient or a layer (see Fig.~\ref{fig:dr5_sol} (right)) can be responsible for a mesh refinement (the benefits lead by such a refinement are further commented below).\\
A more quantitative analysis of the metric-based approach is provided in Table~\ref{tab:dr5}, which collects 
the number of dofs (second column), the minimum diameter of the adapted mesh (third column) together with the $L^2(\Omega)$-norm of the discretization error (fourth column) at the last eight adaptation levels (first column). From 
a geometric viewpoint, the number of dofs increases by a factor varying between $3.5$ and $4$ throughout the iterative levels, while the minimum diameter halves iteration after iteration. Concerning the error, the reduction is more significant during the first levels before stagnating around $10^{-6}$. We remark that the adaptive algorithm breaks before matching tolerance $\mathrm{tol}$ since the maximum number of iterations has been performed.
\begin{table}
    \centering
    \begin{tabular}{clll||l}
        adapt. level & dofs & \(h_{\mathrm{min}}\) & \(\|u-u_h\|_{L^2(\Omega)}\) & \(\|u-u_h\|_{L^2(\Omega)}\) \\\hline
        7 & 81 & 0.176777 & \(2.82441 \cdot 10^{-3}\) & \(1.89502 \cdot 10^{-1}\) \\
        8 & 289 & 0.0883883 & \(6.82165 \cdot 10^{-4}\) & \(1.84833 \cdot 10^{-1}\) \\
        9 & 1033 & 0.0441942 & \(1.83441 \cdot 10^{-4}\) & \(2.81516 \cdot 10^{-1}\) \\
        10 & 3829 & 0.0220971 & \(5.33397 \cdot 10^{-5}\) & \(1.36216 \cdot 10^{-1}\) \\
        11 & 14245 & 0.0110485 & \(3.05214 \cdot 10^{-5}\) & \(2.01399 \cdot 10^{-2}\) \\
        12 & 55149 & 0.00552427 & \(1.16476 \cdot 10^{-5}\) & \(3.82096 \cdot 10^{-3}\) \\
        13 & 215321 & 0.00276214 & \(4.59463 \cdot 10^{-6}\) & \(1.02780 \cdot 10^{-3}\) \\
        14 & 853511 & 0.00138107 & \(2.98502 \cdot 10^{-6}\) & \(-\)
    \end{tabular}
    \caption{Test case 2 (circular discontinuity). Quantitative information about the metric-based algorithm (left panel); comparison between the metric-based approach proposed here and the geometric procedure proposed in \cite{raeli2018finite} (fourth vs. fifth column) in terms of $L^2(\Omega)$-norm of the discretization error.}
    \label{tab:dr5}
\end{table}

We have compared the output provided by Algorithm~\ref{algo:metrics} with the adapted mesh proposed 
in~\cite{raeli2018finite} where the authors exploit a geometric indicator based on the standard concept of distance from the interface. The discrepancy between the two adaptive procedures is evident both from a qualitative and a quantitative viewpoint. The final adapted mesh in~\cite{raeli2018finite} confines the mesh refinement along the interface between $G$ and $S$. 
The quantitative comparison has been carried out by setting a maximum level of refinement as well as an upper bound on the mesh cardinality. The last column in Table~\ref{tab:dr5} is useful to cross-compare the two adaptive approaches in terms of accuracy, since gathering the value of the $L^2(\Omega)$-norm of the discretization error. It is evident that, for a fixed adaptation level, the extra-refinement in $G$ allows to reach a higher accuracy with respect to an approach where only the discontinuity of $\varepsilon$ is tracked, with, on average, a gain of two/three orders of accuracy.

\subsection{Test case 3: an advection-diffusion problem with a diagonal advective field}\label{prob:adr2}
\begin{figure}
    \centering
    \subfigure[]{\includegraphics[width=0.45\textwidth]{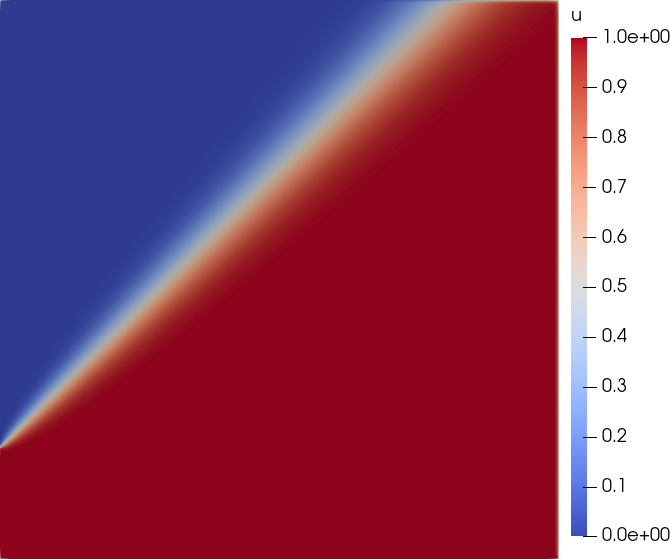}}
    \subfigure[]{\includegraphics[width=0.45\textwidth]{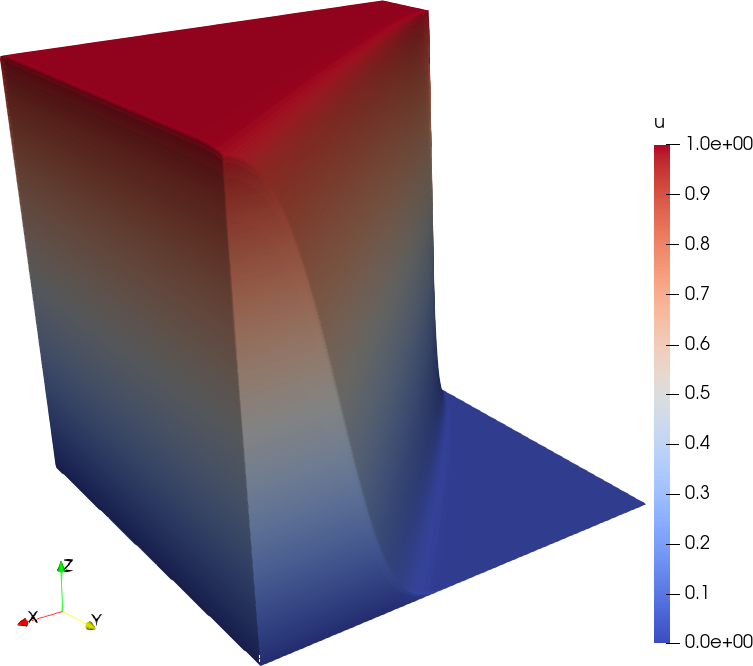}}
    \caption{Test case 3. 2D (left) and 3D (right) view of the solution.}
    \label{fig:adr2_sol}
\end{figure}

As last configuration, we consider the setting proposed in~\cite{Brooks1982199}, where problem \eqref{eq:adeq} is solved on \(\Omega = (0, 1)^2\) for \(\varepsilon = 10^{-6}\), \(\bm{\beta} = [\cos(\vartheta), \sin(\vartheta)]^T / \varepsilon\),
$b=0$ and \(f = 0\). The problem is provided with Dirichlet data, namely we impose
\begin{equation*}
u(x, y) =
\left\{
\begin{alignedat}{3}
& 1 & \quad & \text{on } \Gamma = \{(0,y)\, :\, 0 < y \leq 0.2\} \cup \{(x,0)\ :\ 0\le x\le 1\}, \\
& 0 & \quad & \text{on } \partial\Omega \setminus \Gamma.
\end{alignedat}
\right.
\end{equation*}
The solution identified by this choice of data exhibits an internal diagonal layer, with a slope equal to $\pi/4$, together with two very steep boundary layers, the former along the right side of $\Omega$, the latter in correspondence with the right portion of the top side (see Fig.~\ref{fig:adr2_sol}). We remark that the diagonal layer is characterized by a non-uniform thickness.

\begin{figure}
    \centering
    \subfigure[Adapted mesh: $i= 1$.]{\includegraphics[width=0.4\textwidth]{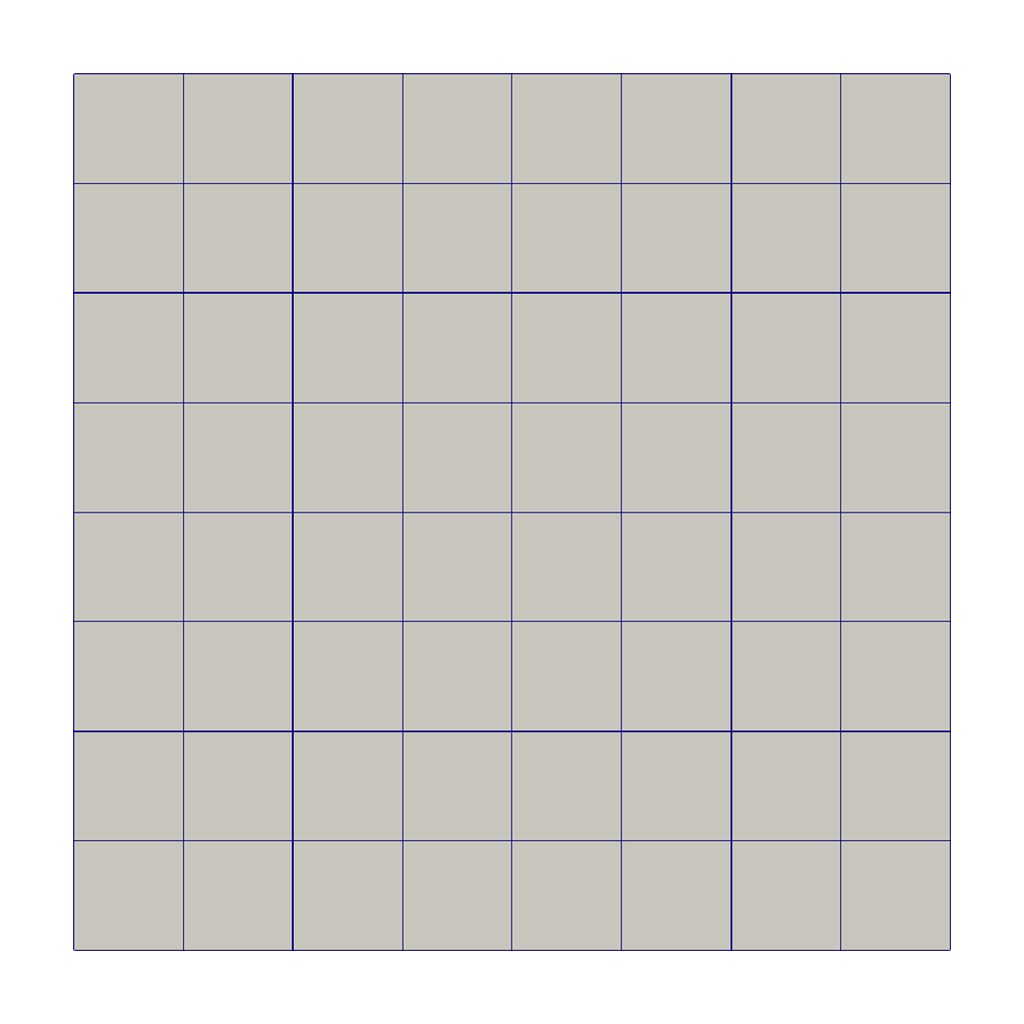}}
    \subfigure[Adapted mesh: $i= 3$.]{\includegraphics[width=0.4\textwidth]{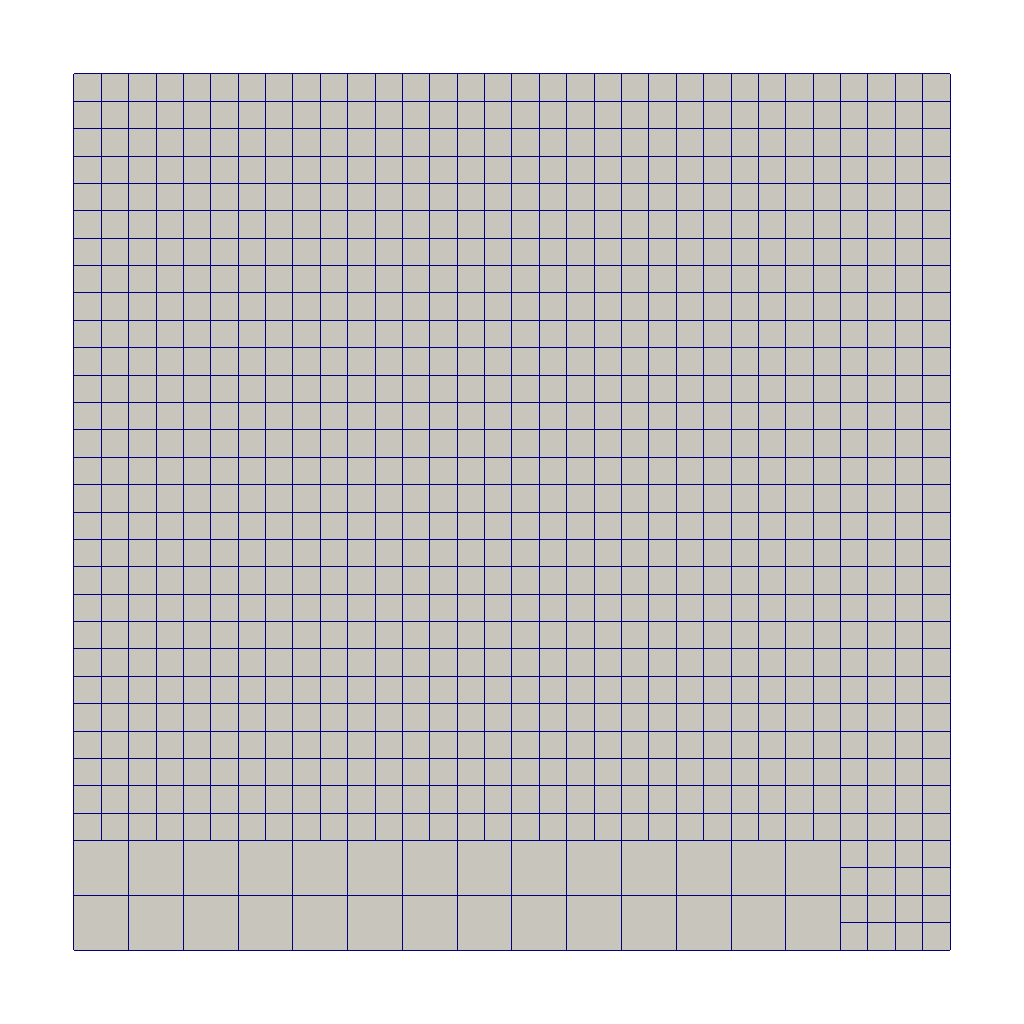}} \\
    \subfigure[Adapted mesh: $i= 6$.]{\includegraphics[width=0.4\textwidth]{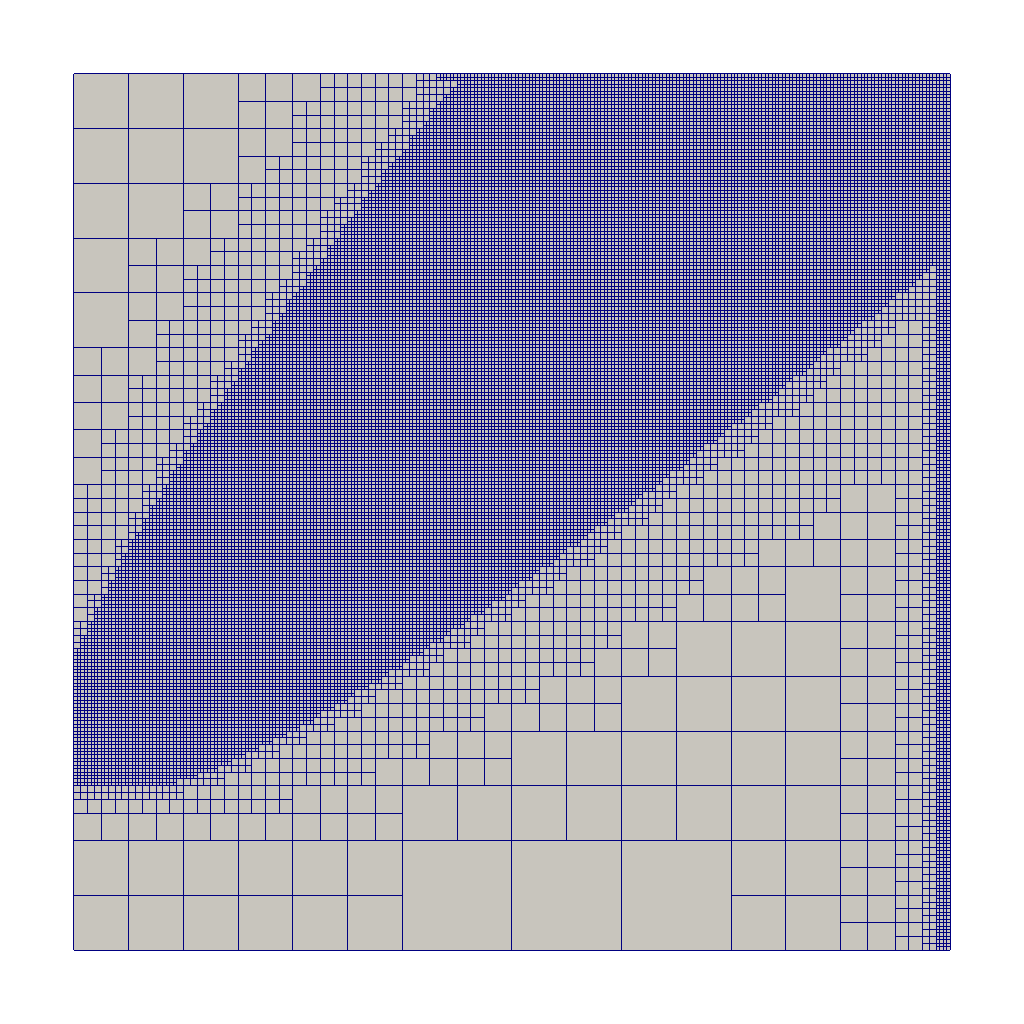}}
    \subfigure[Adapted mesh: $i= 10$]{\includegraphics[width=0.4\textwidth]{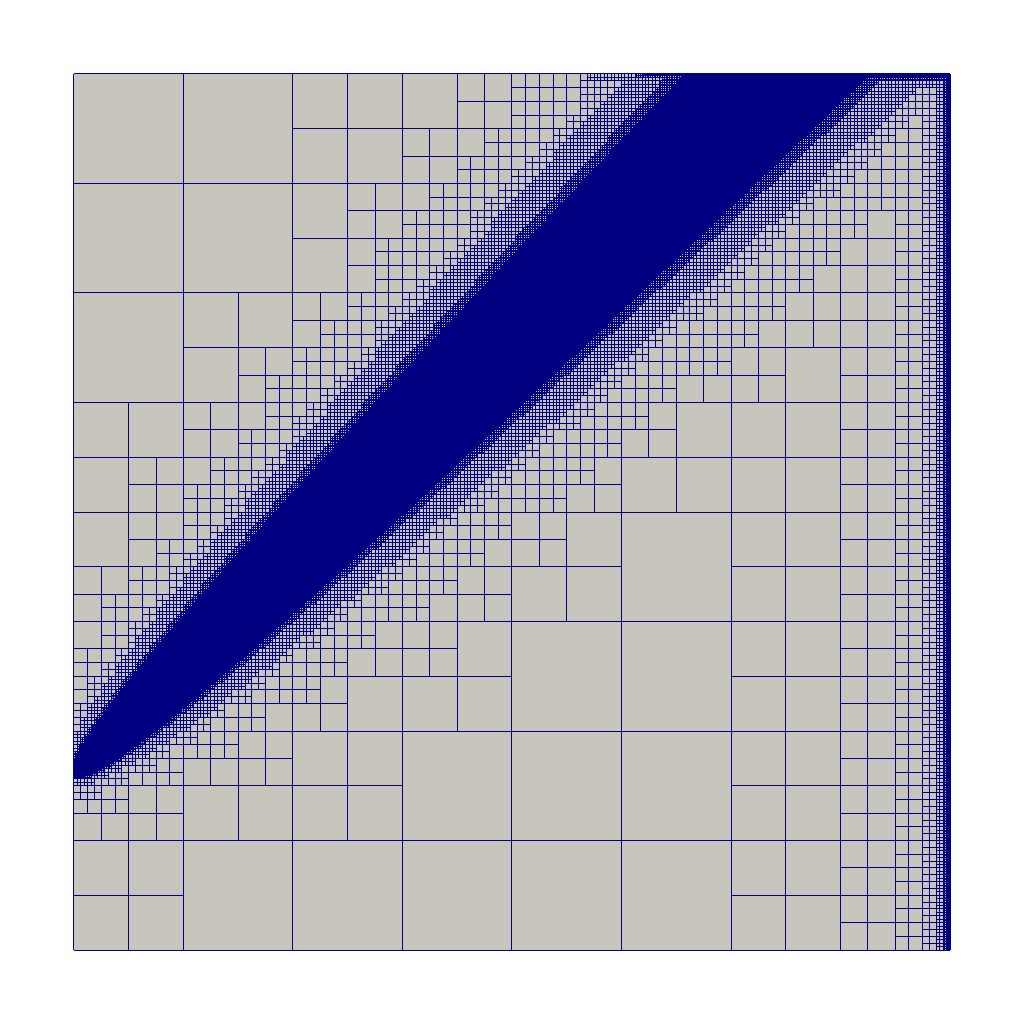}}
    \caption{Test case 3. Marking-based adaptation.}
    \label{fig:adr2_marker}
\end{figure}

\begin{figure}
    \centering
    \subfigure[Adapted mesh: $i= 2$. ]{\includegraphics[width=0.4\textwidth]{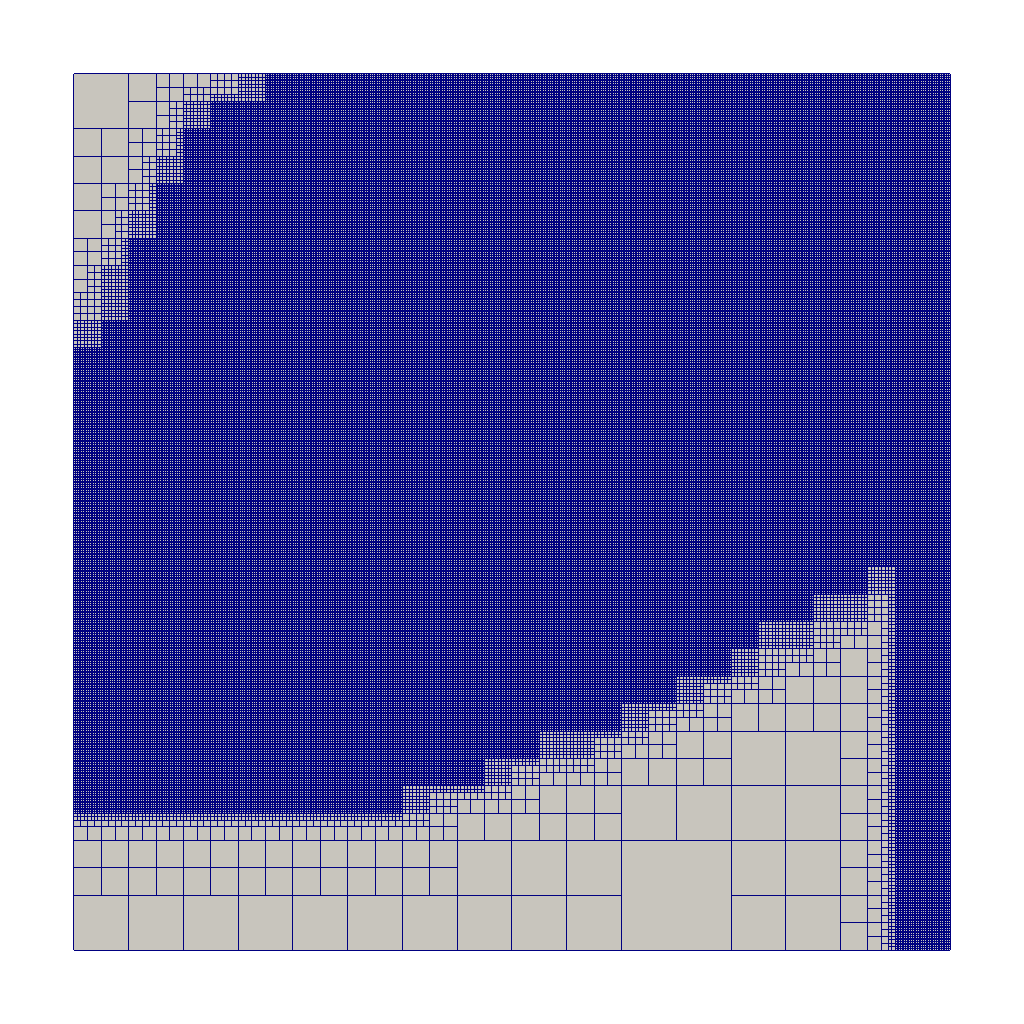}}
    \subfigure[Adapted mesh: $i= 4$. ]{\includegraphics[width=0.4\textwidth]{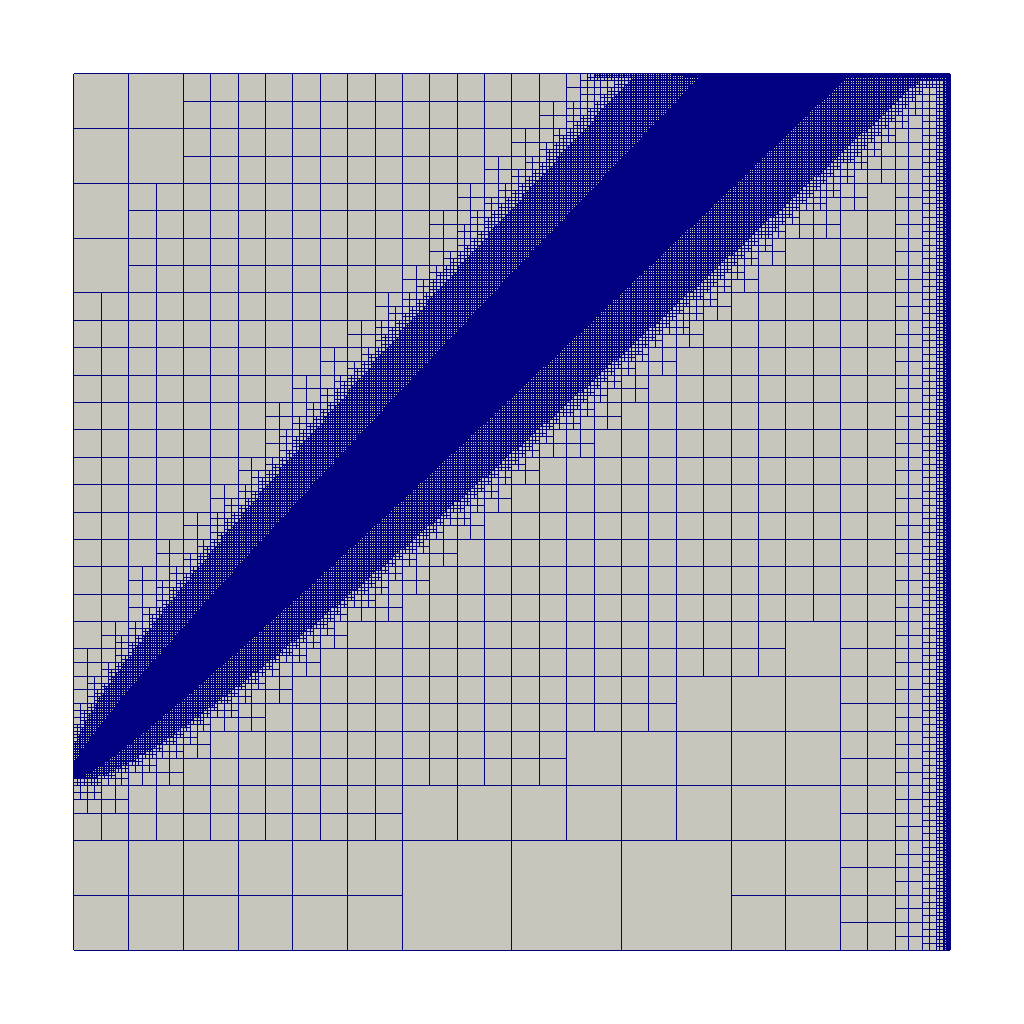}}
    \caption{Test case 3. Metric-based adaptation.}
    \label{fig:adr2_metrics}
\end{figure}

Figures~\ref{fig:adr2_marker} and~\ref{fig:adr2_metrics} show the meshes generated at certain iterations of Algorithms~\ref{algo:marking}
and~\ref{algo:metrics}, respectively. Both the adaptive procedures are run starting from an initial structured mesh characterized by $4$ subdivisions along the two axes and by setting \(\mathrm{tol}=10^{-6}\) and \(i_{\mathrm{max}}=10\).
The internal and the boundary layers are detected by both the procedures, and vary the level of refinement consistently with the different thickness of the layers. 
In particular, the boundary layers which are extremely thin are very sharply captured, while it is evident the non-uniform as well as greater thickness of the diagonal layer. 
As for the previous test cases, the metric-based approach converges in a fewer number of iterations when compared with the marking-based method. However, in contrast to Test cases 1 and 2 where the metric-based approach detects layers more sharply, here the refined areas exhibit a similar amplitude for both the procedures, although the way adopted to cluster the elements is different in the two cases (more confined and finer around the diagonal boundary for Algorithm~\ref{algo:marking}, wider and coarser for Algorithm~\ref{algo:metrics}).
In terms of efficiency, the metric-based approach outperforms the marking-based algorithm, the number of dofs being considerably lower ($150233$ in \Cref{fig:adr2_metrics}(b) versus $1371972$ in \Cref{fig:adr2_marker}(d)).

\begin{figure}
    \centering
{\includegraphics[width=0.45\textwidth]{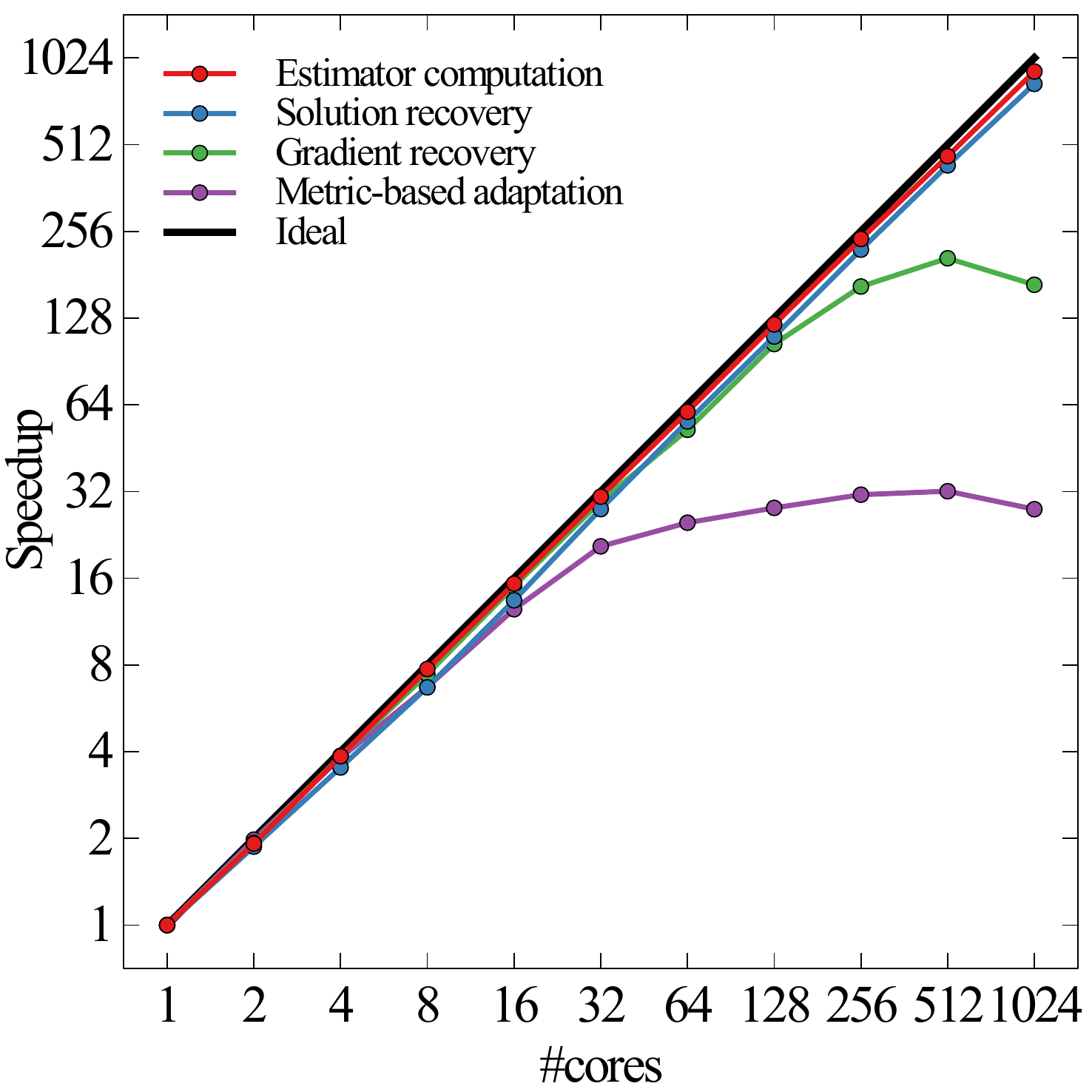}}\hspace*{0.8cm}{\includegraphics[width=0.45\textwidth]{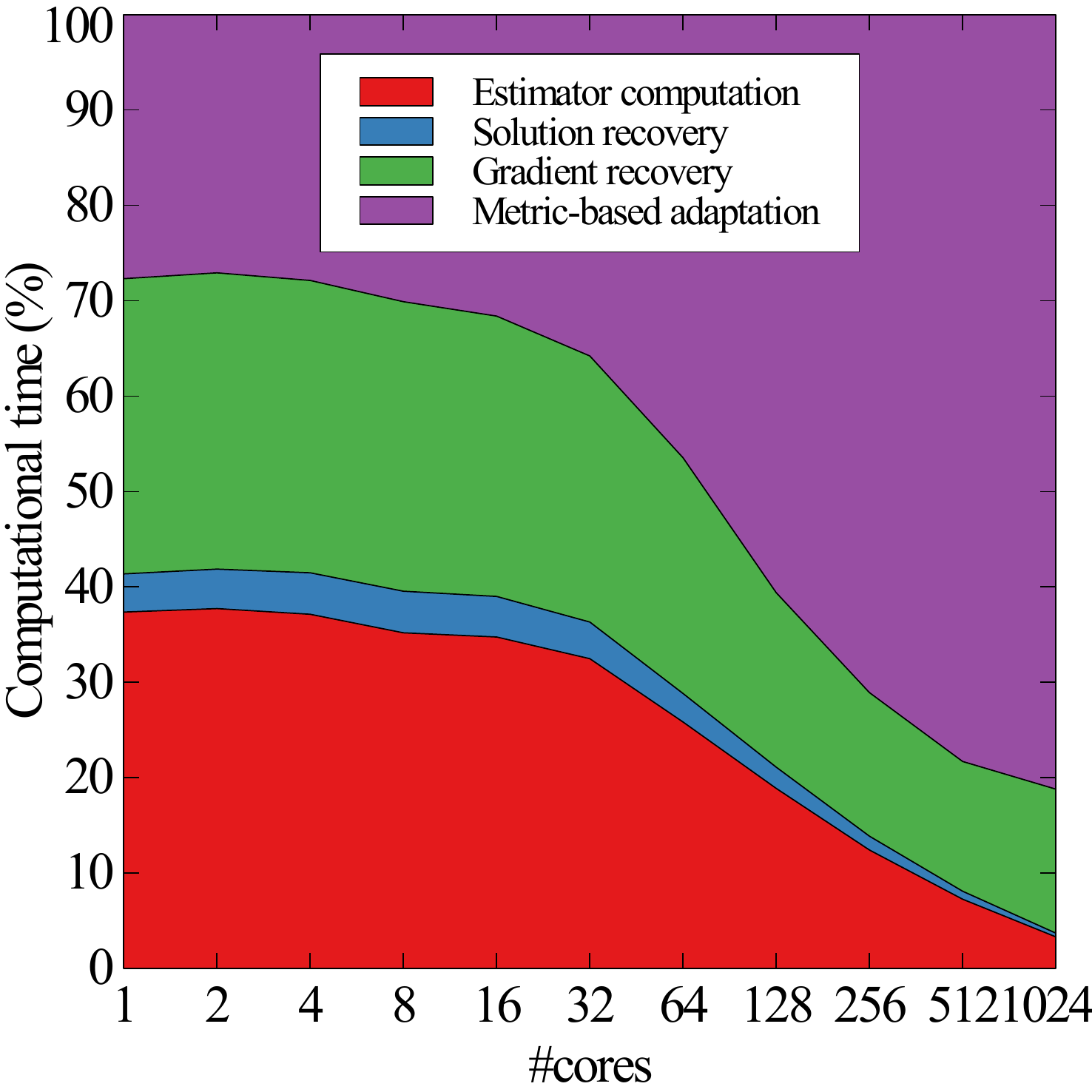}}
\caption{Test case 3. Scalabilty analysis of the
different phases of the metric-based adaptation procedure:
parallel speedup (left) and stacked barplot of the time percentage associated with the different phases of the algorithm (right) as a function of the number of parallel processes.}
    \label{fig:adr2_metrics_speedup}
\end{figure}

We exploit this test case to investigate the scalability of the metric-based adaptation procedure, by mimicking the analysis carried out in \Cref{fig:dr1_marker_speedup}.
\Cref{fig:adr2_metrics_speedup}
gathers the associated results. 
A cross-comparison between the two figures highlights a similar trend, although for this test case the bottleneck is represented by the adaptation procedure (which stagnates at about \(32\) cores). On the contrary, the gradient recovery exhibits a better speedup, approximately linear up to \(512\) cores. 
The deterioration of the mesh adaptation speedup finds a justification in the continuous change of the computational load predicted by \ref{prima_ell} for each cell as well as of the mesh partitioning among processes. 
Moreover, concerning the computational effort, running the whole algorithm requires 
a very contained time, ranging from half a second to few seconds.

\begin{figure}
    \centering
    \subfigure[Adapted mesh: $i=2$.]{\includegraphics[width=0.4\textwidth]{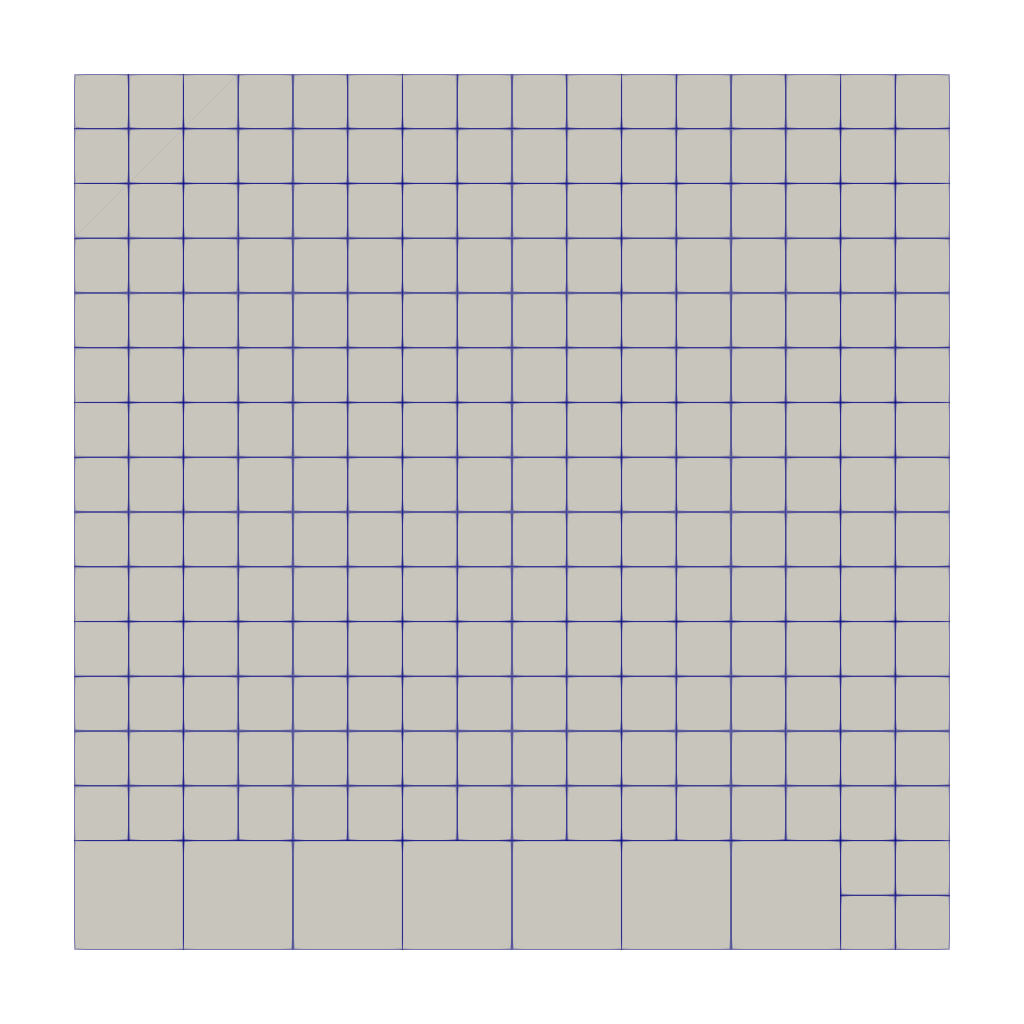}}
    \subfigure[Adapted mesh: $i=4$.]{\includegraphics[width=0.4\textwidth]{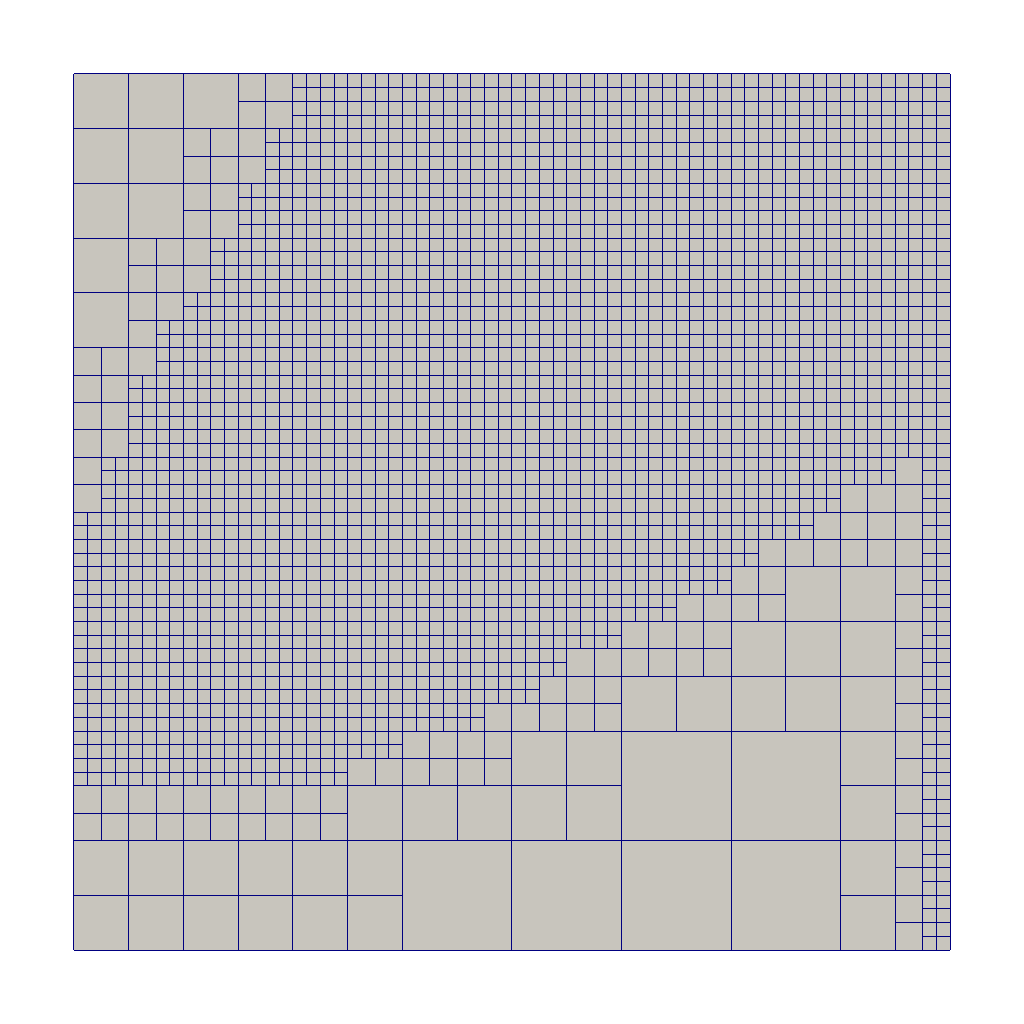}} \\
    \subfigure[Adapted mesh: $i=6$.]{\includegraphics[width=0.4\textwidth]{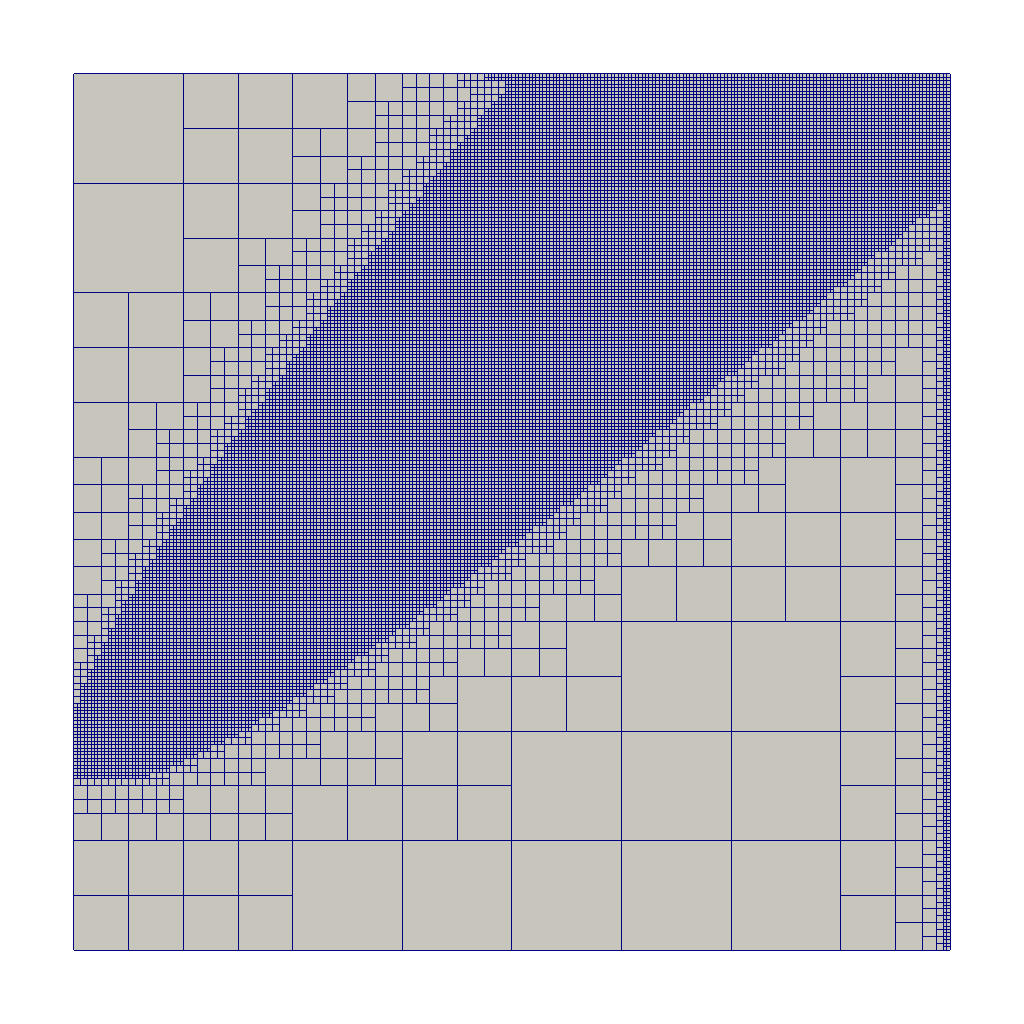}}
    \subfigure[Adapted mesh: $i=10$.]{\includegraphics[width=0.4\textwidth]{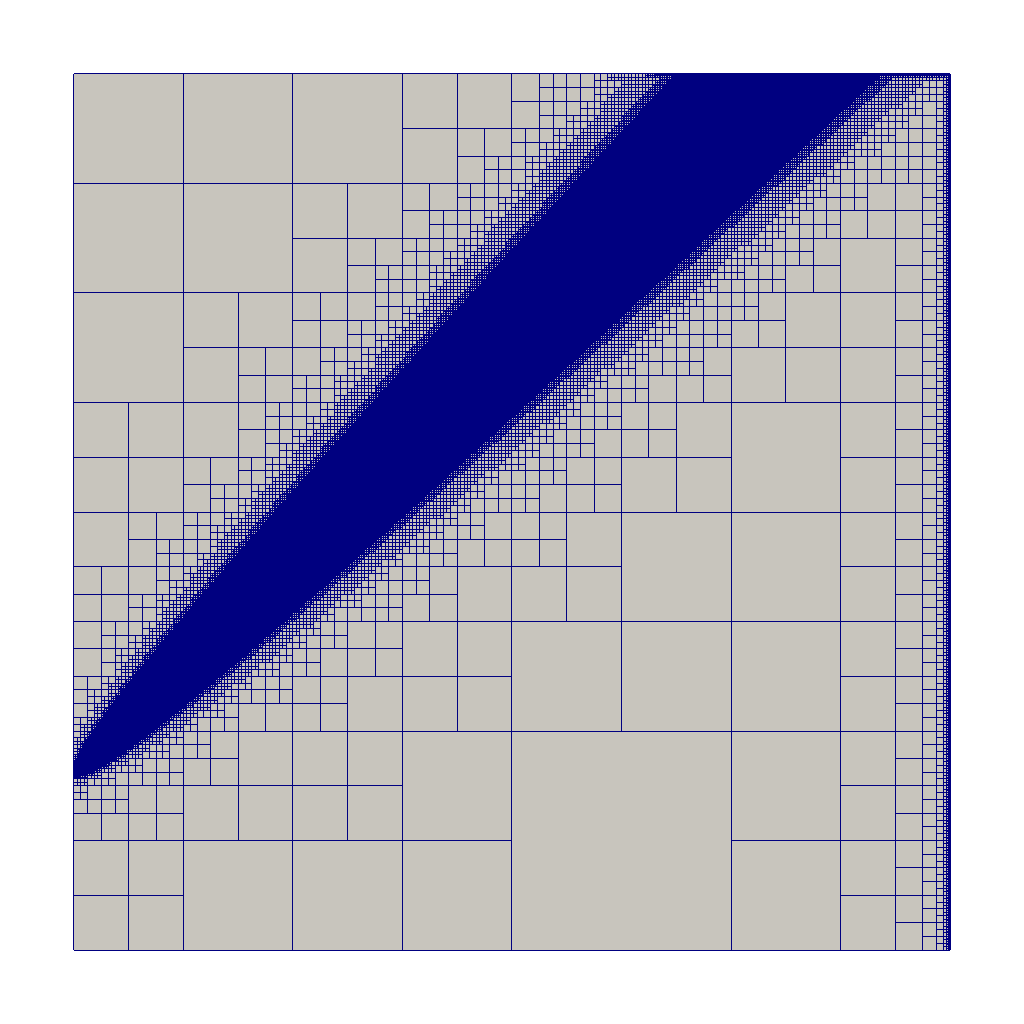}}
    \caption{Test case 3. Adaptation based on~\cite{schillinger2012isogeometric}.}
    \label{fig:adr2_isogeo}
\end{figure}

Finally, for comparison purposes, we have applied to the configuration at hand also the adaptive procedure proposed 
in~\cite{schillinger2012isogeometric}. The authors compute the area-averaged $L^2$-norm of the discrete gradient over each mesh element,
\begin{equation*}
    \gamma_k = \frac{\norm{\nabla u_h}_{L^2(\Omega^{(k)})}}{\left|\Omega^{(k)}\right|},
\end{equation*}
and then mark elements for refinement or coarsening according to condition
\begin{equation*}
    \gamma_k \geq \frac{c_1}{N_\mathrm{el}} \sum_{k=1}^{N_\mathrm{el}} \gamma_k \quad \mbox{or}\quad 
    \gamma_k \leq \frac{c_2}{N_\mathrm{el}} \sum_{k=1}^{N_\mathrm{el}} \gamma_k,
\end{equation*}
respectively with \(c_1, c_2\) user-defined parameters to be empirically tuned to confine the refinement/coarsening to the most problematic squares. Starting from the same initial mesh adopted in Figures~\ref{fig:adr2_marker}-\ref{fig:adr2_metrics}, after $10$ iterations we obtain the mesh in 
Fig.~\ref{fig:adr2_isogeo} which shows the final grid together with three preliminary adapted meshes.
In particular,  
the final tessellation is fully comparable with the one yielded by the marking-based algorithm in Fig.~\ref{fig:adr2_marker}, the main features of the solution being identified with a similar refinement pattern and after the same number ($10$) of iterations. The discrepancy between the two approaches is in terms of number of dofs, since having $1371972$ nodes in \Cref{fig:adr2_marker}(d) to be compared with $1628325$ nodes of \Cref{fig:adr2_isogeo}(d).
%
%

\section{Conclusions and perspectives}
\label{sec:conclusions}
This work formalizes a new mesh adaptation strategy for Cartesian quadtree meshes. The adaptation process is driven by the metric associated with a recovery based estimator for the $L^2(\Omega)$-norm of the dicretization error. Suitable finite difference formulas are employed to yield a higher-order approximation for the solution gradient to be, successively, employed in order to define a recovered solution. 
Finally, the discrepancy between the numerical and the recovered solution is used to settle the error estimator behind the mesh adaptation process.

The massive numerical verification in Section~\ref{sec:results} confirms the reliability of the proposed algorithm when applied to advection-diffusion-reaction problems, also in the presence of problem discontinuities not necessarily aligned with the Cartesian directions. Moreover, we have assessed the higher efficiency of the proposed metric-based approach when compared with a well-established estimate--mark--refine method and with other mesh adaptation strategies available in the literature, since guaranteeing, in general, a reduction of the number of estimation--adaptation steps and of dofs. We have also assessed the scalability of the 
metric-based algorithm. We recognize 
a remarkable performance in terms of parallel speedup, which is approximately linear for most of the phases involved in the adaptation process, as well as of
computational time, which remains of the order of few seconds also in the worst cases.

The estimator $\eta$ turns out to be fully general in terms of differential problem to be solved, dimensionality of the problem, and adopted discretization. As a consequence, 
future developments could involve the extension of the proposed algorithm to more complex three-dimensional problems, such as large scale scenarios or time-dependent applications. As a more challenging task, we aim to generalize the a posteriori setting here settled to an anisotropic context, starting from the papers already available in the literature which provide an anisotropic counterpart of the Zienkiewicz-Zhu
error estimator~\cite{enumath,farrell1,farrell2,fpc22}.

\section*{Acknowledgements}
This work was supported by Accordo Quadro ASI-POLIMI ``Attivit\`a di Ricerca e Innovazione'' n. 2018-5-HH.0, collaboration agreement between the Italian Space Agency and Politecnico di Milano.
Simona Perotto thanks the PRIN research grant n.20204LN5N5 Advanced Polyhedral Discretisations of Heterogeneous PDEs for Multiphysics Problems and the INdAM-GNCS 2022 Project Metodi di riduzione computazionale per le scienze applicate: focus su sistemi complessi.

\bibliography{references}

\end{document}